\documentclass[12pt,a4paper]{article}

\usepackage{amsmath}    
\usepackage{amssymb}
\usepackage{graphicx}   
\usepackage{verbatim}   
\usepackage{color}      

\usepackage{hyperref}   

\usepackage{enumerate}
\usepackage{mathrsfs}
\usepackage{empheq}
\usepackage{bbold}
\usepackage[british]{babel}

\usepackage[final]{showlabels}

\usepackage{amsthm}
\newtheorem{definition}{Definition}[section]

\newtheorem{theorem}{Theorem}[section]

\newtheorem{remark}{Remark}[section]

\newcommand{\del}{\partial}
\renewcommand{\theta}{\vartheta}
\renewcommand{\phi}{\varphi}
\newcommand{\vecc}[2]{\left ( \begin{array}{c}#1\\#2\\ \end{array}\right )}

\newcommand{\dd}{\mathrm{d}}

\newcommand{\id}{\mathbb{1}}

\renewcommand{\vec}{\mathbf}

\renewcommand{\title}{A semi-discrete Active Flux method for the Euler equations on Cartesian grids}

\newcommand{\authorTwo}{Wasilij Barsukow\footnote{Bordeaux Institute of Mathematics, Bordeaux University and CNRS/UMR5251, Talence, 33405 France}}
\newcommand{\authorOne}{Rémi Abgrall\footnote{Institute for Mathematics \& Computational Science, University of Zurich, Winterthurerstrasse 190, CH-8057 Zurich, Switzerland}}
\newcommand{\authorThree}{Christian Klingenberg\footnote{Institute for Mathematics, University of Wurzburg, Emil-Fischer-Strasse 40, 97074 Wurzburg, Germany}}

\begin{document}

\begin{center} \Large
\title

\vspace{1cm}

\date{}
\normalsize

\authorOne, \authorTwo, \authorThree
\end{center}

\begin{abstract}

Active Flux is an extension of the Finite Volume method and additionally incorporates point values located at cell boundaries. This gives rise to a globally continuous approximation of the solution. Originally, the Active Flux method emerged as a fully discrete method, and required an exact or approximate evolution operator for the point value update. For nonlinear problems such an operator is often difficult to obtain, in particular for multiple spatial dimensions. We demonstrate that a new semi-discrete Active Flux method (first described in \cite{abgrall22} for one space dimension) can be used to solve nonlinear hyperbolic systems in multiple dimensions without requiring evolution operators. We focus here on the compressible Euler equations of inviscid hydrodynamics and third-order accuracy. We introduce a multi-dimensional limiting strategy and demonstrate the performance of the new method on both Riemann problems and subsonic flows.



\end{abstract}

\section{Introduction}

The Active Flux method uses as its degrees of freedom both cell averages and point values at cell interfaces. The point values are shared between adjacent cells. In one spatial dimension, any choice of the reconstruction in one cell that interpolates those point values will give rise to a globally continuous reconstruction. In multiple spatial dimensions, it also seems natural to ensure global continuity, i.e. continuity all along the cell interface and not just in the point values. Traditionally, discontinuous reconstructions were favored for hyperbolic conservation laws; in \cite{roe18}, Active Flux is found to perform better than the Discontinuous Galerkin method. 

While the averages require a conservative update, the update of the point values is essentially not restricted by more than the condition that the resulting method should be stable. To this end it needs to incorporate upwinding, and the earliest version of the Active Flux method (\cite{vanleer77}, for linear advection in 1-d) traced a characteristic back to the time level $t^n$ where a reconstruction of the data was evaluated. {In \cite{kerkmann18}, approximate evolution operators were derived using a Cauchy-Kowalevskaya/ADER procedure. As the reconstruction is globally only $C^0$, Riemann problems in the derivatives were solved. This led to Active Flux methods for nonlinear systems of conservation laws in one spatial dimension.} {For scalar, nonlinear conservation laws, a fixpoint iteration can be used to systematically generate approximations to the speed of the characteristic of arbitrary order of accuracy (\cite{barsukow19activeflux,chudzik21}).} This approach was extended in \cite{barsukow19activeflux} to hyperbolic systems of conservation laws in one spatial dimension, even if they do not allow for characteristic variables. {Predictor-corrector estimates of the eigenvalues and the eigenvectors of the Jacobian led to third-order accurate approximate evolution operators.} They were used, for example, in \cite{barsukow20swaf} to solve the one-dimensional shallow water equations in presence of dry areas. 

For hyperbolic systems in multiple spatial dimensions, even if they are linear, characteristic curves in general no longer exist. Also, values in general are not transported, but the solution is a convolution of the initial data with a more or less complicated kernel. For the acoustic equations with the speed of sound $c$, for example, the solution in $\vec x$ at time $t$ depends on the initial data in a disc with radius $ct$ around $\vec x$. This disc is the interior of the intersection of the hypersurface of initial data with the cone of bicharacteristics which has its vertex at $(t, \vec x)$. In \cite{eymann13}, {\cite{fan15}}, a solution operator was given for the acoustic equations, which relied on smoothness of the initial data, and in \cite{barsukow17} a solution in the sense of distributions was obtained which could be used to solve e.g. Riemann problems. These operators can be implemented efficiently and used to update the point values in an Active Flux method for linear acoustics (\cite{barsukow18activeflux}, {\cite{chudzik23}}). An approximate evolution operator for linear acoustics based on bicharacteristics is used in \cite{chudzik23}. 
For the Euler equations, splitting and linearization is suggested in \cite{roe18}. Generally, when studying the approximation error using Taylor series for nonlinear problems, evolution operators designed using linearization are found to require an additional fix in order to achieve third-order accuracy.  For the one-dimensional case, in \cite{barsukow19activeflux} a general algorithm to achieve third order of accuracy on non-linear problems while solving only linear ones was derived. Analogous high-order approximate evolution operators for general multi-dimensional nonlinear systems of conservation laws are currently unavailable, but for the Euler equations some suggestions can be found in \cite{maeng17}.

All these Active Flux methods are fully-discrete. In \cite{zeng14,abgrall20}, a semi-discrete version of Active Flux was introduced. In order to obtain an equation for the point values, the spatial derivative in the PDE is discretized using finite difference formulae. This approach is immediately applicable to all kinds of nonlinear problems without the need to derive an evolution operator, but at the price of a reduced CFL condition. In \cite{roe21}, a link is made between the semi-discrete approach and the reduced CFL condition for high-order methods with a compact stencil in space. 

In \cite{abgrall22,abgrall22proceeding} the semi-discrete Active Flux has been applied to one-dimensional nonlinear problems, and extended to arbitrary order. The aim of the present work is, maintaining $3^\text{rd}$ order of accuracy, to extend it to the multi-dimensional Euler equations.
The paper is organized as follows: Section \ref{sec:semidiscreteAF} describes the method and Section \ref{sec:limiting} presents a novel multi-dimensional limiting strategy. Numerical results are shown in Section \ref{sec:numerical}.

\section{The semi-discrete Active Flux method} \label{sec:semidiscreteAF}

Here, we let ourselves guide by the semi-discrete approach of \cite{abgrall22} and extend it to multi-dimensional Cartesian grids, while aiming at third-order accuracy. Consider a hyperbolic $m \times m$ system of conservation laws in $d$ spatial dimensions\footnote{Boldface letters denote ``spatial'' vectors, i.e. those whose natural dimension is that of space ($d$). Other collections of scalars (such as the conserved quantities $q$) are not typeset in boldface.}
\begin{align}
 \del_t q + \nabla \cdot \vec f(q) &= 0 \qquad q \colon \mathbb R^+_0 \times \mathbb R^d \to \mathbb R^m \label{eq:conslaw}
\end{align}
For simplicity, we restrict ourselves to two spatial dimensions ($d=2$) and write $\vec f = (f^x, f^y)$, $\nabla_q f^x = J^x$, $\nabla_q f^y = J^y$ whenever convenient.

\subsection{Update of the averages}

\begin{figure}
 \centering
 \includegraphics[width=0.3\textwidth]{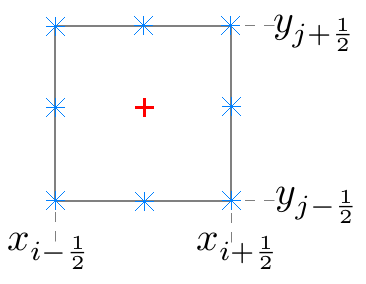}
 \caption{Degrees of freedom employed in the third-order accurate Active Flux method on Cartesian grids. Stars denote point values, the cross denotes the cell average. Every cell has access to 8 point values and one average. Point values are shared with adjacent cells; per cell there is one average and $4 \cdot \frac14 + 4 \cdot \frac12 = 3$ point values.}
 \label{fig:dofs}
\end{figure}

Integrating \eqref{eq:conslaw} over the Cartesian cell
\begin{align}
 C_{ij} := \left[ x_{i-\frac12}, x_{i+\frac12} \right ] \times \left[ y_{j-\frac12}, y_{j+\frac12} \right ]
\end{align}
and denoting the cell average by
\begin{align}
 \bar q_{ij}(t) := \frac{1}{\Delta x \Delta y} \int_{C_{ij}} q(t, \vec x) \dd \vec x
\end{align}
one finds
\begin{align}
 \frac{\dd}{\dd t} \bar q_{ij} + \frac{1}{\Delta x \Delta y} \int_{\del C_{ij}}  \vec n \cdot \vec f(q) &= 0 \label{eq:semidiscretegauss}
\end{align}
As there are degrees of freedom located at the boundary $\del C_{ij}$ of cell $C_{ij}$, we intend to use them as quadrature points for a sufficiently accurate quadrature of the integral appearing in \eqref{eq:semidiscretegauss}. Inspired by previous approaches (e.g. \cite{barsukow18activeflux,kerkmann18}) we use three Gauss-Lobatto points per edge, where the extreme points (corners) are shared (see Figure \ref{fig:dofs}). This is in contrast to \cite{zeng19} where a distribution of point values based on a Gauss-Legendre quadrature along the edge is suggested, or \cite{zeng14} where only second-order accuracy is obtained. Note also that we enforce global continuity: the point values on an edge are the same as seen from either of the adjacent cells and a value at a corner is involved in the update of four cells. This is in contrast to e.g. discontinuous Galerkin methods.

On Cartesian grids it is convenient to adopt the following notation for the 8 point values on the boundary of cell $C_{ij}$:
\begin{center}
 \begin{minipage}[l]{0.5\textwidth}
\begin{align*}
 &q_{i-\frac12,j+\frac12} &&q_{i,j+\frac12} &&q_{i+\frac12,j+\frac12}\\
 &q_{i-\frac12,j} &&                       &&q_{i+\frac12,j}\\
 &q_{i-\frac12,j-\frac12} &&q_{i,j-\frac12} &&q_{i+\frac12,j-\frac12}
\end{align*}
 \end{minipage}
\end{center}

Then, 
\begin{align}
 \frac{\dd}{\dd t} \bar q_{ij} &+ \frac{1}{\Delta x \Delta y} \int_{y_{j-\frac12}}^{y_{j+\frac12}} \dd y \Big ( f^x(q(t, x_{i+\frac12},y)) - f^x(q(t, x_{i-\frac12},y)) \Big ) \\
 &\nonumber+ \frac{1}{\Delta x \Delta y} \int_{x_{i-\frac12}}^{x_{i+\frac12}} \dd x \Big ( f^y(q(t, x, y_{j+\frac12})) - f^y(q(t, x, y_{j-\frac12})  )\Big) = 0
\end{align}
using Simpson's rule ($\omega_{-\frac12} = \omega_{\frac12} = \frac16$, $\omega_0 = \frac23$) becomes
\begin{align}
 \frac{\dd}{\dd t} \bar q_{ij}(t) &+ \frac1{\Delta x} \sum_{K=-\frac12,0,\frac12} \omega_K \left( f^x(q_{i+\frac12,j+K}(t)) - f^x(q_{i-\frac12,j+K}(t)) \right ) \\
 &\nonumber+ \frac{1}{\Delta y} \sum_{K=-\frac12,0,\frac12} \omega_K  \left ( f^y(q_{i+K,j+\frac12}(t)) - f^y(q_{i+K,j-\frac12}(t))  \right) = 0 \label{eq:averageupdatesemidiscrete}
\end{align}
This method is conservative, with e.g. the $x$-flux through the cell interface $(i+\frac12,j)$ being given by
\begin{align}
 \hat f^x_{i+\frac12,j} &= \sum_{K=-\frac12,0,\frac12} \omega_K  f^x(q_{i+\frac12,j+K})\\
 &= \frac{f^x(q_{i+\frac12,j-\frac12}) + 4f^x(q_{i+\frac12,j}) + f^x(q_{i+\frac12,j+\frac12})}{6}
\end{align}
It is also at least $3^\text{rd}$ order accurate, as it is exact for biparabolic functions.

\subsection{Update of the point values}

\newcommand{\diag}{\mathrm{diag}}

The update of the cell averages, as described above, now needs to be complemented by an update of the point values. In the one-dimensional case, it was proposed in \cite{abgrall22} to replace the spatial derivatives appearing in \eqref{eq:conslaw} by finite differences. Here, the multi-dimensional case will be addressed. Note first that hyperbolicity of \eqref{eq:conslaw} implies that it is always possible to define the positive and negative parts of the Jacobians via their eigenvalues. With $J^x = R \diag(\lambda_1, \ldots, \lambda_m)R^{-1}$ one has
 \begin{align}
 (J^x)^+ &:= R \diag(\lambda_1^+, \ldots, \lambda_m^+)R^{-1} \label{eq:jacobiansplitupwindplus}\\
 (J^x)^- &:= R \diag(\lambda_1^-, \ldots, \lambda_m^-)R^{-1} \label{eq:jacobiansplitupwindminus}
 \end{align}
 where, for scalars $a \in \mathbb R$ the positive/negative parts are simply $a^+ = \max(0, a)$, $a^- = \min(0, a)$.

The finite difference formulae are obtained by differentiating a reconstruction. Define first the unique biparabolic polynomial 
\begin{align}
q_{ij,\text{recon}} &\in P^{2,2}, \quad q_{ij,\text{recon}} \colon \left[ - \frac{\Delta x}{2}, \frac{\Delta x}{2} \right] \times \left[ - \frac{\Delta y}{2}, \frac{\Delta y}{2} \right] \to \mathbb R^m \\
P^{2,2} &= \mathrm{span}(1, x, x^2, y, xy, x^2 y, y^2, xy^2, x^2 y^2)
\end{align}
that interpolates the degrees of freedom accessible to cell $ij$:
 \begin{align*}
  q_{ij,\text{recon}}\left(-\frac{\Delta x}{2}, \frac{\Delta y}{2}  \right) &= q_{i-\frac12,j+\frac12} & q_{ij,\text{recon}}\left( \frac{\Delta x}{2}, \frac{\Delta y}{2} \right) &= q_{i+\frac12,j+\frac12}\\
  q_{ij,\text{recon}}\left( -\frac{\Delta x}{2}, -\frac{\Delta y}{2} \right) &= q_{i-\frac12,j-\frac12} & q_{ij,\text{recon}}\left( \frac{\Delta x}{2}, -\frac{\Delta y}{2} \right) &= q_{i+\frac12,j-\frac12}\\
  q_{ij,\text{recon}}\left( 0, \frac{\Delta y}{2} \right) &= q_{i,j+\frac12}  & q_{ij,\text{recon}}\left(- \frac{\Delta x}{2}, 0 \right) &= q_{i-\frac12,j} \\                     
 q_{ij,\text{recon}}\left( \frac{\Delta x}{2}, 0 \right) &= q_{i+\frac12,j} &
 q_{ij,\text{recon}}\left( 0, -\frac{\Delta y}{2} \right) &= q_{i,j-\frac12} 
 \end{align*}
and
\begin{align}
 \frac{1}{\Delta x \Delta y} \int_{-\frac{\Delta x}{2}}^{\frac{\Delta x}{2}} \int_{-\frac{\Delta y}{2}}^{\frac{\Delta y}{2}} q_{ij,\text{recon}}(x, y) \, \dd y \dd x &= \bar q_{ij}
\end{align}

 This reconstruction has already been used in \cite{barsukow18activeflux,kerkmann18} and is given there explicitly. Then we define the finite differences in the corner as
 \begin{align}
  (D^x)^+_{i+\frac12,j+\frac12}q &:= \del_x q_{ij,\text{recon}}\left.\left(x, \frac{\Delta y}{2}\right) \right|_{x = \frac{\Delta x}{2}} \label{eq:findifffirst}\\
  (D^x)^-_{i+\frac12,j+\frac12}q &:= \del_x q_{i+1,j,\text{recon}}\left.\left(x, \frac{\Delta y}{2}\right) \right|_{x = -\frac{\Delta x}{2}} \\
  (D^y)^+_{i+\frac12,j+\frac12}q &:= \del_y q_{ij,\text{recon}}\left.\left(\frac{\Delta x}{2},y\right) \right|_{y = \frac{\Delta y}{2}} \\
  (D^y)^-_{i+\frac12,j+\frac12}q &:= \del_y q_{i,j+1,\text{recon}}\left.\left(\frac{\Delta x}{2},y\right) \right|_{y = -\frac{\Delta y}{2}} 
 \end{align}
 Observe that due to continuity,
 \begin{align}
  (D^x)^+_{i+\frac12,j+\frac12} q&= \del_x q_{i,j+1,\text{recon}}\left(x, -\frac{\Delta y}{2}\right) \Big|_{x = \frac{\Delta x}{2}}
 \end{align}
 such that this would be an equivalent definition that gives the same result (and similarly for the other finite differences). Analogously, we define the finite differences on the edges
 \begin{align}
  (D^x)^+_{i+\frac12,j} q&:= \del_x q_{ij,\text{recon}}\left.\left(x, 0\right) \right|_{x = \frac{\Delta x}{2}} \\
  (D^x)^-_{i+\frac12,j}q &:= \del_x q_{i+1,j,\text{recon}}\left.\left(x, 0\right) \right|_{x = -\frac{\Delta x}{2}} \\
  (D^y)_{i+\frac12,j} q&:= \del_x q_{ij,\text{recon}}\left.\left(\frac{\Delta x}{2},y\right) \right|_{y = 0} \label{eq:findifflast}
 \end{align}
 Observe that due to continuity, there is no distinction between $(D^y)^+_{i+\frac12,j}$ and $(D^y)^-_{i+\frac12,j}$. Here, again, the symmetric definition
 \begin{align}
  (D^y)_{i+\frac12,j}q &:= \del_y q_{i+1,j,\text{recon}}\left.\left(-\frac{\Delta x}{2},y\right) \right|_{y = 0} 
 \end{align}
 yields the same result. The derivatives at $(i,j+\frac12)$ are obtained analogously. For reference we now state their explicit forms:

 \begin{align*}
\left (D^x\right )^+_{i+\frac12,j}q &= \frac{1}{4 \Delta x} \left(    4 \left (-9 \bar q_{ij} +2 \left (q_{i-\frac12,j}+2 q_{i+\frac12,j}\right )\right )+4  \left (q_{i,j-\frac12}+q_{i,j+\frac12} \right ) \right . \\   &\left .
+ q_{i-\frac12,j-\frac12}+q_{i+\frac12,j-\frac12}+q_{i-\frac12,j+\frac12}+q_{i+\frac12,j+\frac12} \right )\\
\left (D^x\right )^-_{i+\frac12,j} q&= -\frac{1}{4 \Delta x}  \left(  -36 \bar q_{i+1,j} +8  \left (2q_{i+\frac12,j}+q_{i+\frac32,j}\right ) +q_{i+\frac12,j-\frac12}  \right . \\   &\left .+4  \left (q_{i+1,j-\frac12}+q_{i+1,j+\frac12}\right )
 + q_{i+\frac32,j-\frac12} + q_{i+\frac12,j+\frac12} + q_{i+\frac32,j+\frac12} \right )\\
\left (D^y\right )_{i+\frac12,j} q&= \frac{ q_{i+\frac12,j+\frac12}-q_{i+\frac12,j-\frac12}}{\Delta y}\\
\left (D^y\right )^+_{i,j+\frac12} q&= \frac{1}{4 \Delta y }      \left( 4 \left (q_{i-\frac12,j}-9 \bar q_{ij} +q_{i+\frac12,j}\right ) 
+q_{i-\frac12,j-\frac12} + q_{i-\frac12,j+\frac12} \right . \\   &\left .+ q_{i+\frac12,j-\frac12}+ q_{i+\frac12,j+\frac12}  +8 \left (q_{i,j-\frac12} +2 q_{i,j+\frac12} \right ) \right )\\
\left (D^y\right )^-_{i,j+\frac12}q &= -\frac{1}{4 \Delta y}    \left ( 4 \left (q_{i-\frac12,j+1}-9 \bar q_{i,j+1} +q_{i+\frac12,j+1} \right )
+q_{i-\frac12,j+\frac12}  \right . \\   &\left . + q_{i+\frac12,j+\frac12}+ q_{i-\frac12,j+\frac32} + q_{i+\frac12,j+\frac32} +8  \left (2q_{i,j+\frac12}+q_{i,j+\frac32}\right ) \right )\\
\left (D^x\right )_{i,j+\frac12}q &= \frac{q_{i+\frac12,j+\frac12} -q_{i-\frac12,j+\frac12}}{\Delta x}\\
\left (D^x\right )^+_{i+\frac12,j+\frac12} q&= \frac{q_{i-\frac12,j+\frac12}-4 q_{i,j+\frac12} +3 q_{i+\frac12,j+\frac12} }{\Delta x}\\
\left (D^x\right )^-_{i+\frac12,j+\frac12} q&= \frac{4 q_{i+1,j+\frac12} - 3q_{i+\frac12,j+\frac12}-q_{i+\frac32,j+\frac12}}{\Delta x}\\
\left (D^y\right )^+_{i+\frac12,j+\frac12} q&= \frac{q_{i+\frac12,j-\frac12}-4 q_{i+\frac12,j} +3 q_{i+\frac12,j+\frac12} }{\Delta y }\\
\left (D^y\right )^-_{i+\frac12,j+\frac12} q&= \frac{4 q_{i+\frac12,j+1} - 3q_{i+\frac12,j+\frac12}-q_{i+\frac12,j+\frac32}}{\Delta y}
 \end{align*}

 However, in some situations one might be willing to employ a different reconstruction, as is, for instance, the case in Section \ref{sec:limiting} concerned with limiting. At this point one has to resort to the more general formulae \eqref{eq:findifffirst}--\eqref{eq:findifflast}.
 
Finally, the upwinding is defined as
\begin{align}
 (J^x D^x_{i+K,j+L})^\text{upw}q := (J^x)^+ (D^x)^+_{i+K,j+L}q + (J^x)^- (D^x)^-_{i+K,j+L}q
\end{align}
with $K,L \in \{ -\frac12, 0, \frac12 \}$ and an analogous definition for $J^y$. We content ourselves here with simple upwinding based on the $x$- and $y$-Jacobians separately. The multi-dimensionality advertised in works of Roe and collaborators is not entirely lost, though, because the degree of freedom at the node couples the waves from different directions in a truly multi-dimensional manner (different from DG methods, say).

We propose to update the point values as follows:
\begin{align}
 \frac{\dd}{\dd t} q_{i+\frac12,j} + (J^x D^x_{i+\frac12,j})^\text{upw} q + J^y D^y_{i+\frac12,j} q &= 0 \label{eq:edgexupdatesemidiscrete}\\
 \frac{\dd}{\dd t} q_{i,j+\frac12} + J^x D^x_{i,j+\frac12} q + (J^y D^y_{i,j+\frac12})^\text{upw} q &= 0 \label{eq:edgeyupdatesemidiscrete}\\
 \frac{\dd}{\dd t} q_{i+\frac12,j+\frac12} + (J^x D^x_{i+\frac12,j+\frac12})^\text{upw} q + (J^y D^y_{i+\frac12,j+\frac12})^\text{upw} q &= 0 \label{eq:nodeupdatesemidiscrete}
\end{align}

As the finite differences are exact for biparabolic function, one expects $3^\text{rd}$ order of accuracy.

The complete method consists of the ODEs \eqref{eq:averageupdatesemidiscrete} (average update), \eqref{eq:edgexupdatesemidiscrete}--\eqref{eq:edgeyupdatesemidiscrete} (point values at edge midpoints) and \eqref{eq:nodeupdatesemidiscrete} (point values at nodes). We propose to integrate these with an SSP-RK3 method. In \cite{abgrall22proceeding}, it was shown for linear advection in one spatial dimension that this approach leads to a stable scheme with a maximum CFL number of 0.41, a value lower than what can be achieved if an evolution operator is available (\cite{chudzik21}), but comparable e.g. to CFL numbers of DG methods. In two space dimensions we expect stability at least for half the CFL number, i.e. 0.2, or possibly higher.

\section{Limiting} \label{sec:limiting}

Existing approaches to limiting in the context of standard Finite Volume methods modify the values of the reconstruction at a cell interface. They cannot be used for Active Flux due to its global continuity and the fact that point values at cell interfaces are prescribed and cannot be modified arbitrarily. Limiting employed in \cite{kerkmann18} therefore gives up on continuity. Approaches to limiting that maintain continuity so far have only been treating the situation in which a parabolic reconstruction of monotonic discrete data (point values and average) is not monotonic, i.e. has an artificial extremum. In \cite{roe15}, a piecewise linear/parabolic reconstruction is used in this case, and in \cite{barsukow19activeflux} the same situation is handled by replacing the parabola by a power law. One can show that then the reconstruction is always monotonic whenever the discrete data are. Such modified reconstructions are effective in drastically reducing spurious oscillations, but they do not guarantee to remove them entirely. This is because the update of the averages is not limited and can itself create artificial extrema in the discrete data. However, in absence of better approaches, e.g. the power-law reconstruction is a viable limiting strategy. In particular, it is not computationally intensive.

In multiple spatial dimensions, a similar strategy is presented here for the first time. The multi-dimensional case is, however, much more complex because every cell has access to 8 point values. Even monotonicity is a vague concept, because data can be monotonic along one direction and non-monotonic along the other. We focus only on the degrees of freedom accessible to a cell. As long as the cell average is between the smallest and the largest point values we propose to reconstruct in such a way that the value of the reconstruction at any point inside the cell is also between the smallest and the largest point value. This is always possible, as will be shown. We additionally impose a monotonicity constraint along any edge, i.e. we reconstruct in a monotonic fashion if the data are monotonic. Finally, we also impose global continuity.

\begin{figure}
\centering
\includegraphics[width=0.49\textwidth]{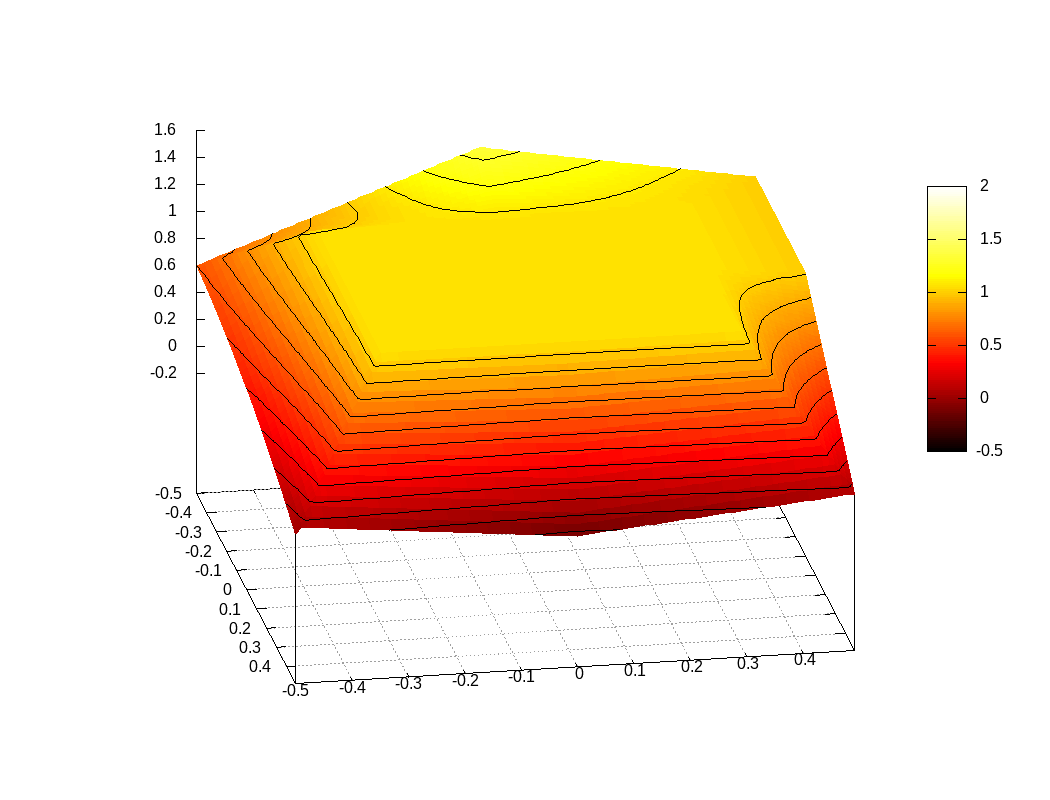} \includegraphics[width=0.49\textwidth]{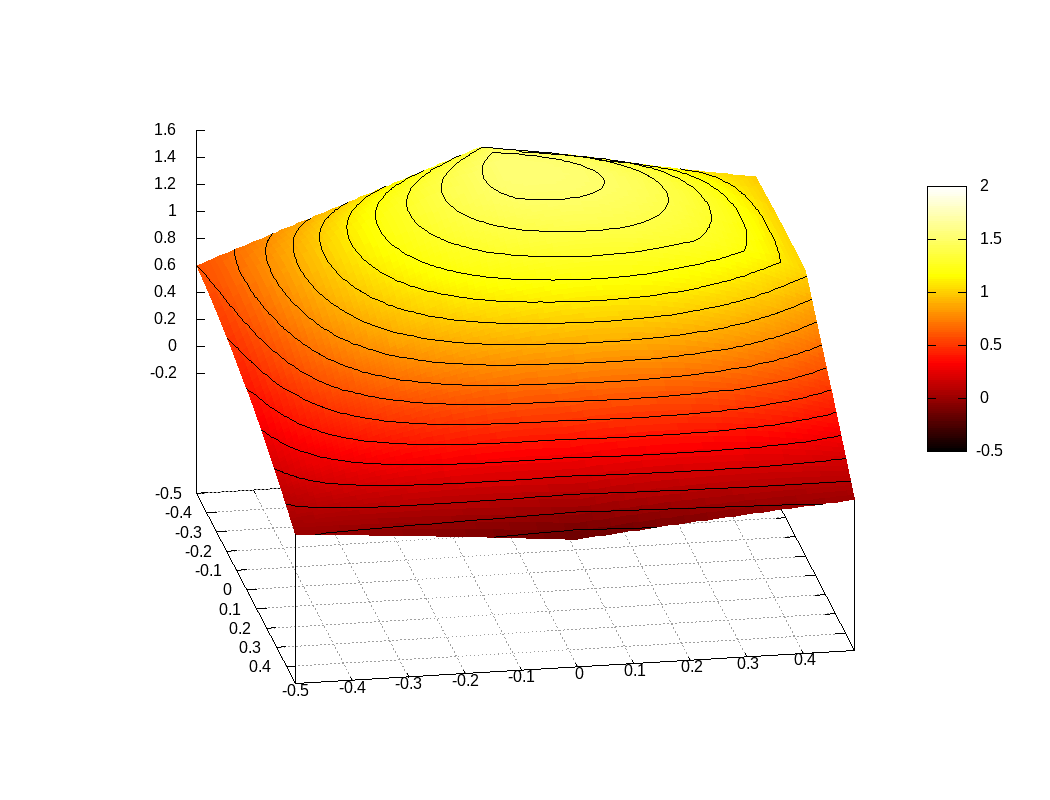} 
 \caption{ \emph{Left}: An example of a plateau reconstruction. Here, $q_\text{NW} = 1$, $q_\text{W}  = 1.35$, $q_\text{SW} = 0.6$, $q_\text{S}  = 0.4$, $q_\text{SE} = 0$, $q_\text{E}  = -0.2$, $q_\text{NE} = 0.0$, $q_\text{N}  = 1$, $ \bar q  = 0.9$ (the S-edge is on the left and $\Delta x = \Delta y = 1$). All edges but the S-edge are reconstructed as hats, the S-edge is reconstructed parabolically. \emph{Right}: A piecewise-biparabolic reconstruction of the same data; one clearly observes an overshoot. The isolines have a spacing of 0.1.}
 \label{fig:plateau}
\end{figure}

\begin{figure}
\centering
\includegraphics[width=0.49\textwidth]{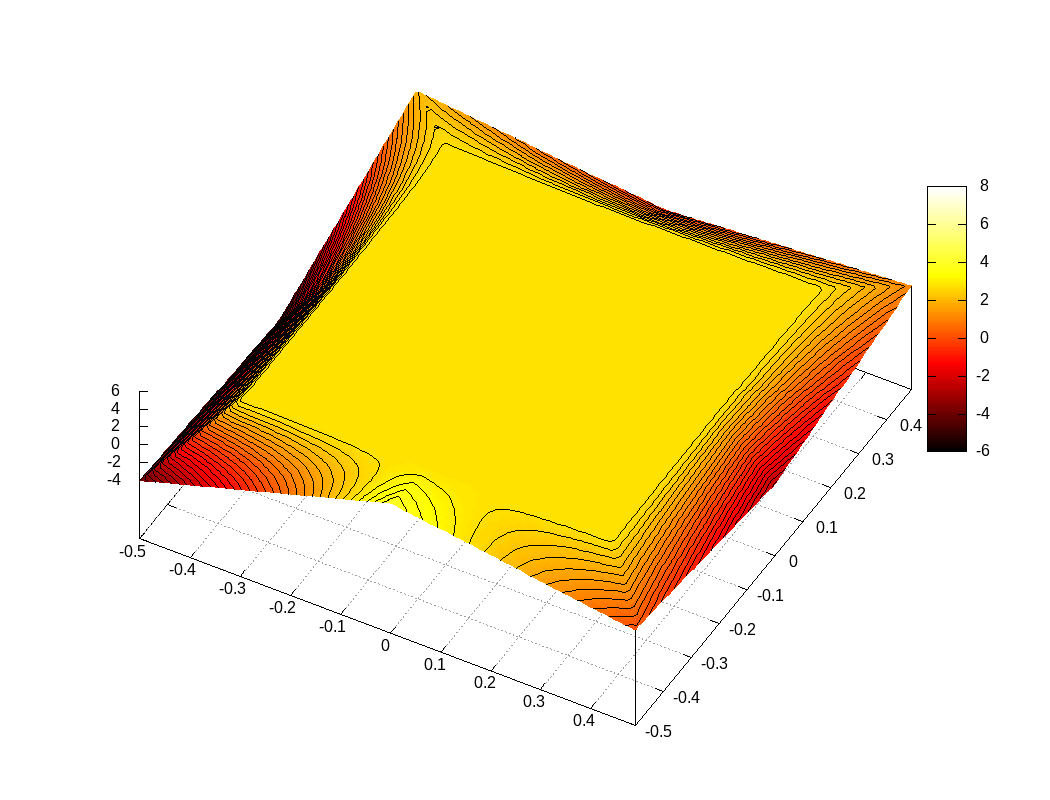} \includegraphics[width=0.49\textwidth]{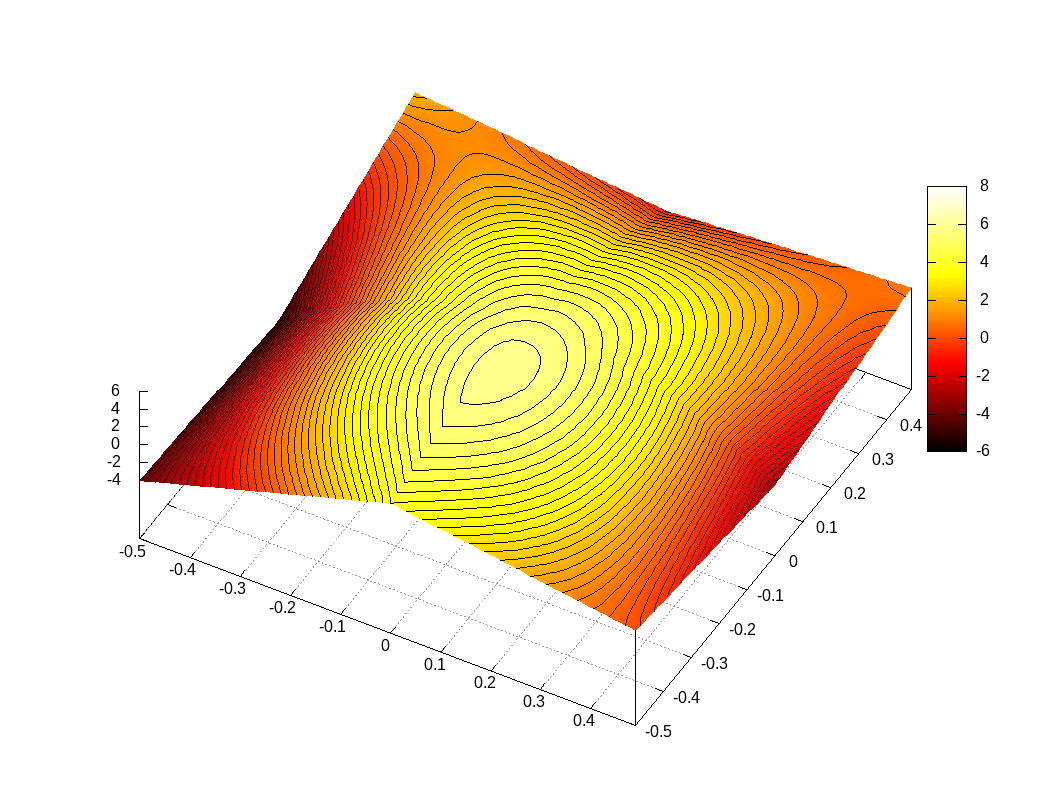} 
 \caption{ \emph{Left}: An example of a plateau reconstruction. $q_\text{NE} = 1, q_\text{NW} = 2, q_\text{SW} = -4, q_\text{SE} = 0, q_\text{N} = -1, q_\text{S} = 4, q_\text{W} = -5, q_\text{E} = -3,  \bar q  = 2$ (the W-edge is on the left and $\Delta x = \Delta y = 1$). All edges are reconstructed as hats. \emph{Right}: A piecewise-biparabolic reconstruction of the same data; one clearly observes an overshoot. The isolines have a spacing of 0.25.}
 \label{fig:plateau2}
\end{figure}

Consider point values at edge centers $q_\text{N}$, $q_\text{S}$, $ q_\text{W}$, $ q_\text{E}$ and at vertices $q_\text{NE}$, $ q_\text{SE}$, $ q_\text{NW}$, $ q_\text{SW}$ of a (reference) Cartesian cell $c = [-\frac{\Delta x}{2}, \frac{\Delta x}{2}] \times [-\frac{\Delta y}{2}, \frac{\Delta y}{2}]$ and a cell average $\bar q$ to be given. We will refer to the four edges as N-edge, S-edge, W-edge and E-edge, respectively. The reconstruction will be denoted by $q_\text{recon} \colon c \to \mathbb R$ for simplicity.
There exist two ways how the reconstruction can create new extrema, which can occur independently of each other:

\begin{enumerate}[1.]

\item It can happen that the parabolic reconstruction along an edge (as part of a biparabolic reconstruction in the cell) overshoots/undershoots the three point values along the edge in question. For the example of an N-edge, this happens if either
\begin{itemize}
 \item the point values $q_\text{NW}$, $q_\text{N}$, $q_\text{NE}$ are not monotonic and $q_\text{NW} \neq q_\text{NE}$, or if
 \item they are monotonic (i.e. either $q_\text{NW} < q_\text{N} < q_\text{NE}$ or $q_\text{NW} > q_\text{N} > q_\text{NE}$), but
 \begin{align}
 \left |q_\text{N} - \frac{q_\text{NE} + q_\text{NW}}{2} \right | > \frac{|q_\text{NE} - q_\text{NW}|}{4} \label{eq:conditionhatreconedge}
\end{align}
such that the parabolic reconstruction has an extremum not present in the discrete data.
\end{itemize}

In this case the reconstruction along the edge will be chosen continuous piecewise linear (``hat''). We will say that the \textbf{reconstruction along the edge is limited}, or just that the ``edge is limited''. To ensure continuity, the reconstruction in any cell with a limited edge can then no longer be biparabolic, but needs to be modified as detailed below and in Section \ref{sec:piecewisebiparabolic}.

\item Define 
\begin{align}
 m := \min(q_\text{N}, q_\text{S}, q_\text{W}, q_\text{E}, q_\text{NE}, q_\text{SE}, q_\text{NW}, q_\text{SW} ) \\
 M := \max(q_\text{N}, q_\text{S}, q_\text{W}, q_\text{E}, q_\text{NE}, q_\text{SE}, q_\text{NW}, q_\text{SW} ) 
\end{align}
It can happen that despite
\begin{align}
 m < \bar q < M \label{eq:quasimonotone}
\end{align}
the reconstruction $q_\text{recon}$ inside the cell $c$ has an extremum with no counterpart in the discrete data, i.e.
\begin{align}
 \exists \vec x \in c \text{ such that either } q_\text{recon}(\vec x) < m \text{ or }q_\text{recon}(\vec x) > M
\end{align}

This situation will be improved by introducing a piecewise defined reconstruction with a central region where the function is constant (``plateau''), and connecting the plateau to the (parabolic or hat) reconstructions along the edges in a continuous fashion. More details are given below and in Section \ref{sec:plateau}; Figures \ref{fig:plateau} and \ref{fig:plateau2} show examples. This new reconstruction fulfills \begin{align}
 m < q_\text{recon}(\vec x) < M \qquad \forall \vec x \in c
\end{align}
We will say that the \textbf{reconstruction inside the cell is limited}, or just that the ``cell is limited''.

This situation appears already in 1-d, in which case it has been suggested in \cite{barsukow19activeflux} to replace the parabolic reconstruction in the cell by a power law. A multi-dimensional analogue of the power law seems unfeasible, though, and we resort here to a piecewise defined, but easier function.

\end{enumerate}

The two situations are independent: any number of edges along the boundary of a cell might require limiting, and this will not generally imply anything about whether the cell itself is to be limited. The possible presence of hat functions along the boundary requires the reconstruction inside the cell to flexibly adapt to the different combinations of edge-reconstructions in order to be continuous. For instance, the plateau reconstruction needs to connect the plateau continuously to either a parabola, or a hat function (see Section \ref{sec:plateau}). Also, if there exists at least one edge that is reconstructed as a hat function, then one cannot use a biparabolic reconstruction inside the cell any longer, even if the cell is not limited. If at least one edge is limited, but the cell is not, a piecewise-biparabolic reconstruction will be used, detailed in Section \ref{sec:piecewisebiparabolic}.

As we are aiming at a globally continuous reconstruction, that is computed locally from merely the cell average and the point values of the cell, the reconstruction along an edge can only depend on the three values associated to this edge, and cannot depend on other values in the cell. Indeed, if edge-reconstruction of one of edges of $c$ were to depend on, say, the average in the cell $c$, then the reconstruction in the neighbouring cell $c'$ would also need to know about the average in $c$.

Due to the particular choice of degrees of freedom for Active Flux the reconstruction has to fulfill two types of conditions: It is supposed to interpolate the point values at cell interfaces and its average is supposed to be equal to the given one. The latter condition -- merely to simplify the calculations -- will be replaced by a (yet unknown) point value $q_\text{C}$ at cell center which is kept as a variable in the formulae. Once the type of reconstruction in all regions of the cell has been determined, their integrals over the respective domains of definition can easily be found as functions of $q_\text{C}$, and $q_\text{C}$ is then determined by imposing the average of the reconstruction over the entire cell. This is a linear equation in $q_\text{C}$ due to linearity of the interpolation problem which makes $q_\text{C}$ enter linearly everywhere. The explicit formulae below therefore also depend on $q_\text{C}$, but the reconstruction in a cell in the end only depends on the point values along its boundary and on its average. This detour does not change the result but simplifies the algorithm.

The overall structure of the reconstruction algorithm is:
\begin{enumerate}
 \item Decide for every edge of the cell whether it is reconstructed parabolically, or as a hat function.
 \item \label{it:piecewbipara} Assume as hypothesis that the cell does not require limiting (i.e. that it is reconstructed in a piecewise biparabolic fashion) and compute the value of $q_\text{C}$ that ensures that the average of the reconstruction agrees with the given cell average.
 \item Check \eqref{eq:quasimonotone} and if true, decide whether the piecewise-biparabolic reconstruction obtained in \ref{it:piecewbipara} has an artificial extremum in the sense of what has been described above\footnote{This happens numerically by testing a given number of locations.}
 \item If this is the case, the cell needs to be limited with a plateau reconstruction. Compute the parameters $\eta, q_\text{p}$ (see below) of the plateau reconstruction that ensure maximum principle preservation and the correct value of the average of the reconstruction.
\end{enumerate}

A pedagogical derivation of the reconstruction algorithm is given in Section \ref{app:recon}. Here, we only state all the relevant results in a concise way.

\begin{theorem} \label{thm:limitedrecon}
 The following reconstruction $q_\text{recon} \colon \left[ -\frac{\Delta x}{2}, \frac{\Delta x}{2} \right] \times \left[ -\frac{\Delta y}{2}, \frac{\Delta y}{2} \right] \to \mathbb R$ is continuous, interpolates all the point values along the boundary of the cell, its average agrees with the given cell average and the reconstruction has the following properties:
 \begin{enumerate}[(i)]
   \item If Condition \eqref{eq:quasimonotone}, i.e. $m < \bar q < M $ is fulfilled, then $m \leq q_\text{recon}(\vec x) \leq M$ for all $\vec x$ inside the cell.
   \item If $q_\text{NW} < q_\text{N} < q_\text{NE}$, then $q_\text{NW} \leq q_\text{recon}(\vec x) \leq q_\text{NE}$ for all $\vec x$ along the N-edge, and similarly for all the other edges.
\end{enumerate}

The definition of the reconstruction is as follows:
If $m < \bar q < M$ is not fulfilled, or if additionally $m < q^\text{pw. biparab.}_\text{recon}(x, y) < M$ for all $(x, y) \in c$, then
\begin{align}
 q_\text{recon}(x, y) := q^\text{pw. biparab.}_\text{recon}(x, y)
\end{align}
otherwise
\begin{align}
 q_\text{recon}(x, y) := q^\text{plateau}_\text{recon}(x, y),
\end{align}
the two types of reconstruction being defined as follows:

\begin{enumerate}[1.]
 \item The piecewise-barabolic reconstruction:
\begin{align}
 q_\text{recon}^\text{pw. biparab.}(x, y) &:= \\
 &\!\!\!\!\!\!\!\!\!\!\!\!\!\!\!\!\!\!\!\!\!\!\!q_\text{recon}^\text{W}\left(\frac{q_\text{SW}-\bar Q}{2}, q_\text{W}-\bar Q, \frac{q_\text{NW}-\bar Q}{2}, x, y, \text{S}, \text{N}, \text{W}, \frac{\bar q-\bar Q}{4}\right) \\&\!\!\!\!\!\!\!\!\!\!\!\!\!\!\!\!\!\!\!\!\!\!\! \nonumber
+q_\text{recon}^\text{S}\left(\frac{q_\text{SE}-\bar Q}{2}, q_\text{S}-\bar Q, \frac{q_\text{SW}-\bar Q}{2}, x, y, \text{E}, \text{W}, \text{S}, \frac{ \bar q-\bar Q }{4}\right)\\
&\!\!\!\!\!\!\!\!\!\!\!\!\!\!\!\!\!\!\!\!\!\!\!\nonumber+q_\text{recon}^\text{N}\left(\frac{q_\text{NW}-\bar Q}{2}, q_\text{N}-\bar Q, \frac{q_\text{NE}-\bar Q}{2}, x, y, \text{W}, \text{E}, \text{N}, \frac{\bar q-\bar Q}{4}\right)\\&\!\!\!\!\!\!\!\!\!\!\!\!\!\!\!\!\!\!\!\!\!\!\! \nonumber
+q_\text{recon}^\text{E}\left(\frac{q_\text{NE}-\bar Q}{2}, q_\text{E}-\bar Q, \frac{q_\text{SE}-\bar Q}{2}, x, y, \text{N}, \text{S}, \text{E}, \frac{ \bar q-\bar Q}{4}\right)\\
				&\!\!\!\!\!\!\!\!\!\!\!\!\!\!\!\!\!\!\!\!\!\!\!\nonumber+ \bar Q
\end{align}
with $\bar Q := \frac{q_\text{SW} + q_\text{W} + q_\text{NW} + q_\text{N} + q_\text{NE} + q_\text{E} + q_\text{SE} + q_\text{S}}{8}$ (other choices are possible) and
\begin{align}
    q_\text{recon}^\text{S}(q_\text{SE}, q_\text{S}, q_\text{SW}, x, y, \text{E}, \text{W}, \text{S},  \bar q ) &= q_\text{recon}^\text{W}(q_\text{SE}, q_\text{S}, q_\text{SW}, y, -x, \text{E}, \text{W}, \text{S},  \bar q )  \label{eq:rotationreconS}\\
    q_\text{recon}^\text{N}(q_\text{NW}, q_\text{N}, q_\text{NE}, x, y, \text{W}, \text{E}, \text{N},  \bar q ) &= q_\text{recon}^\text{W}(q_\text{NW}, q_\text{N}, q_\text{NE}, -y, x, \text{W}, \text{E}, \text{N},  \bar q ) \\
    q_\text{recon}^\text{E}(q_\text{NE}, q_\text{E}, q_\text{SE}, x, y, \text{N}, \text{S}, \text{E},  \bar q ) &= q_\text{recon}^\text{W}(q_\text{NE}, q_\text{E}, q_\text{SE}, -x, -y, \text{N}, \text{S}, \text{E},  \bar q ) \label{eq:rotationreconE}
\end{align}
and
\begin{align}
 &q_\text{recon}^\text{W}(q_\text{SW}, q_\text{W}, q_\text{NW}, x, y, \text{S}, \text{N}, \text{W},  \bar q ) \\\nonumber &\,\,\,\,\,\,\,\,\,\,\,\,\,\,\,\,\,\,\,\,\,\,\,\,\,\,=
 \begin{cases}
  \text{\eqref{eq:bipararecon}} & \text{N,S,W parabolic} \\
\text{\eqref{eq:parahathatleft}--\eqref{eq:parahathatright}} & \text{W parabolic, N,S hat} \\
\text{\eqref{eq:parahatparaleft}--\eqref{eq:parahatpararight}} & \text{W,S parabolic, N hat}\\
\text{\eqref{eq:paraparahatleft}--\eqref{eq:paraparahatright}} & \text{W,N parabolic, S hat}\\
\text{\eqref{eq:hatWparaS} and \eqref{eq:hatWparaN}} & \text{W hat, N,S parabolic}\\
\text{\eqref{eq:hatWparaS} and \eqref{eq:hatWhatNleft}--\eqref{eq:hatWhatNright}} & \text{W,N hat, S parabolic}\\
\text{\eqref{eq:hatWparaN} and \eqref{eq:hathatbottomleft}--\eqref{eq:hathatbottomright}} & \text{W,S hat, N parabolic}\\
\text{\eqref{eq:hatWhatNleft}--\eqref{eq:hatWhatNright} and \eqref{eq:hathatbottomleft}--\eqref{eq:hathatbottomright}} & \text{W,S,N hat}
 \end{cases}
\end{align}
Here, N/S/E/W denote the edges of the cell. $q_\text{C}$ fulfills {\footnotesize
\begin{align*} 
  \begin{cases}
 q_\text{C} = \frac{1}{16}(36  \bar q -q_\text{NW}-q_\text{SW}-4 q_\text{W}) & \text{N,S,W parabolic} \\
 q_\text{C} =\frac1{32} (72  \bar q -3 q_\text{NW}-3q_\text{SW}-8 q_\text{W}) & \text{W parabolic, N,S hat} \\
q_\text{C} =\frac{1}{32} (72  \bar q -3 q_\text{NW}-2 q_\text{SW}-8 q_\text{W}) & \text{W,S parabolic, N hat}\\
q_\text{C} =\frac1{32} (72  \bar q -2 q_\text{NW}-3 q_\text{SW}-8 q_\text{W}) & \text{W,N parabolic, S hat}\\
\bar q = \frac{2 q_\text{C}}9+ \frac{q_\text{SW}+q_\text{W}}{24} + \frac{2 q_\text{C}}9+\frac{q_\text{NW}+q_\text{W}}{24} & \text{W hat, N,S parabolic}\\
\bar q = \frac{2 q_\text{C}}9+ \frac{q_\text{SW}+q_\text{W}}{24} + \frac{2 q_\text{C}}9+\frac1{576} (35 q_\text{NW}+q_\text{SW}+22 q_\text{W}) & \text{W,N hat, S parabolic}\\
\bar q = \frac{2 q_\text{C}}9+\frac{q_\text{NW}+q_\text{W}}{24} + \frac{2 q_\text{C}}9+ \frac1{576} (q_\text{NW}+35 q_\text{SW}+22 q_\text{W}) & \text{W,S hat, N parabolic}\\
\bar q = \frac{2 q_\text{C}}9+\frac1{576} (35 q_\text{NW}+q_\text{SW}+22 q_\text{W}) +  \frac{2 q_\text{C}}9+ \frac1{576} (q_\text{NW}+35 q_\text{SW}+22 q_\text{W}) & \text{W,S,N hat}
 \end{cases}
\end{align*}}

\item The plateau reconstruction:

{\footnotesize

\begin{align*}
 q_\text{recon}^\text{plateau}(x, y) &:= \begin{cases}
                                          q_\text{p} \text{ if } (x, y) \in \left[\Delta x \left(\eta-\frac12 \right), \Delta x \left(\frac12 - \eta \right)\right] \times \left[\Delta y \left(\eta-\frac12 \right), \Delta y \left(\frac12 - \eta \right)\right] \\ 
                                          q_\text{recon}^\text{trapeze,W}(q_\text{SW}, q_\text{W}, q_\text{NW}, x, y, \text{W},  \eta, q_\text{p}) \text{ if }(x,y)\in \text{W-trapeze} \\
                                          q_\text{recon}^\text{trapeze,S}(q_\text{SE}, q_\text{S}, q_\text{SW}, x, y,  \text{S},   \eta, q_\text{p} ) \text{ if }(x,y)\in \text{S-trapeze}\\
                                            q_\text{recon}^\text{trapeze,N}(q_\text{NW}, q_\text{N}, q_\text{NE}, x, y, \text{N},  \eta, q_\text{p}) \text{ if }(x,y)\in \text{N-trapeze} \\
                        q_\text{recon}^\text{trapeze,E}(q_\text{NE}, q_\text{E}, q_\text{SE}, x, y,  \text{E},   \eta, q_\text{p} )\text{ if }(x,y)\in \text{E-trapeze}
                                         \end{cases}
\end{align*}}
 with
\begin{align}
    q_\text{recon}^\text{trapeze,W}(q_\text{SW}, q_\text{W}, q_\text{NW}, x, y, \text{W},  \eta, q_\text{p}) = \begin{cases}
                                                                                                                    \eqref{eq:trapezeWparabolic} & (x, y) \in \text{W parabolic}\\
                                          \text{\eqref{eq:trapezeWhattop}--\eqref{eq:trapezeWhatbottom}} & (x, y) \in \text{W hat}
                                                                                                                   \end{cases}
\end{align}
defined only in {\footnotesize
\begin{align}
\text{W-trapeze} = \left\{(x, y) \text{ s.t. } x \in \left[-\frac{\Delta x}{2}, -\Delta x\left(\frac12- \eta\right) \right] \text{ and } y \in \left[\frac{x}{\Delta y} \Delta x , -\frac{x}{\Delta y} \Delta x  \right ] \right \}
\end{align}}

The reconstructions of the other trapezes are
\begin{align}
    q_\text{recon}^\text{trapeze,S}(q_\text{SE}, q_\text{S}, q_\text{SW}, x, y,  \text{S},   \eta, q_\text{p} ) &= q_\text{recon}^\text{trapeze,W}(q_\text{SE}, q_\text{S}, q_\text{SW}, y, -x, \text{S},   \eta, q_\text{p})  \\
    q_\text{recon}^\text{trapeze,N}(q_\text{NW}, q_\text{N}, q_\text{NE}, x, y, \text{N},  \eta, q_\text{p}) &= q_\text{recon}^\text{trapeze,W}(q_\text{NW}, q_\text{N}, q_\text{NE}, -y, x,  \text{N},  \eta, q_\text{p} ) \\
    q_\text{recon}^\text{trapeze,E}(q_\text{NE}, q_\text{E}, q_\text{SE}, x, y,  \text{E},   \eta, q_\text{p} ) &= q_\text{recon}^\text{trapeze,W}(q_\text{NE}, q_\text{E}, q_\text{SE}, -x, -y,  \text{E},   \eta, q_\text{p} ) 
\end{align}

The parameters $q_\text{p}$ and $\eta$ are found according to the procedure descrobed in Section \ref{ssec:plateaumaximumprinciple}.

\end{enumerate}

\end{theorem}
\begin{proof}
 Continuity is a consequence of Theorems \ref{thm:continuitypara} and \ref{thm:continuityhat} and of the fact that the reconstruction along any edge only involves the points on this edge. 
 The pointwise and average interpolation property follows from Theorems \ref{thm:reconinterpol} and \ref{thm:reconplateauinterpol}.
 Preservation of the maximum principle along the edges is clear from \eqref{eq:conditionhatreconedge} and the idea of reconstructing a hat function along the edge; preservation of the maximum principle for the cell follows from Theorem \ref{thm:reconplateauinterpol}.
\end{proof}

\begin{theorem} 
 The usage of the reconstruction from Theorem \ref{thm:limitedrecon} in every cell leads to a globally continuous reconstruction.
\end{theorem}
\begin{proof}
 It follows by construction (see Sections \ref{ssec:continuitypara}, \ref{ssec:plateauinterpolation}) that the reconstruction in a cell $c$ continuously turns into the reconstruction along the edge as $c \ni (x,y) \to s \in \del c$. The reconstructions along the edges only depend on the three point values located on the edge, and thus the limit as $(x,y)$ approaches the same edge from the other cell is the same.
\end{proof}

\section{Numerical results} \label{sec:numerical}

Here, the Euler equations with $q = (\rho, \rho u, \rho v, e)$,
\begin{align}
 f^x &= (\rho u, \rho u^2 + p, \rho u v, u (e+p)) \\
 f^y &= (\rho v, \rho u v, \rho v^2 + p , v (e+p)) \\
 e &= \frac{p}{\gamma-1} + \frac12 \rho (u^2 + v^2)
\end{align}
and $\gamma=1.4$ are solved using the Active Flux method described above. Initial data are denoted by $\rho_0$, $u_0$, $v_0$, $p_0$.

\subsection{Convergence studies}

\begin{figure}
 \centering
 \includegraphics[width=.49\textwidth]{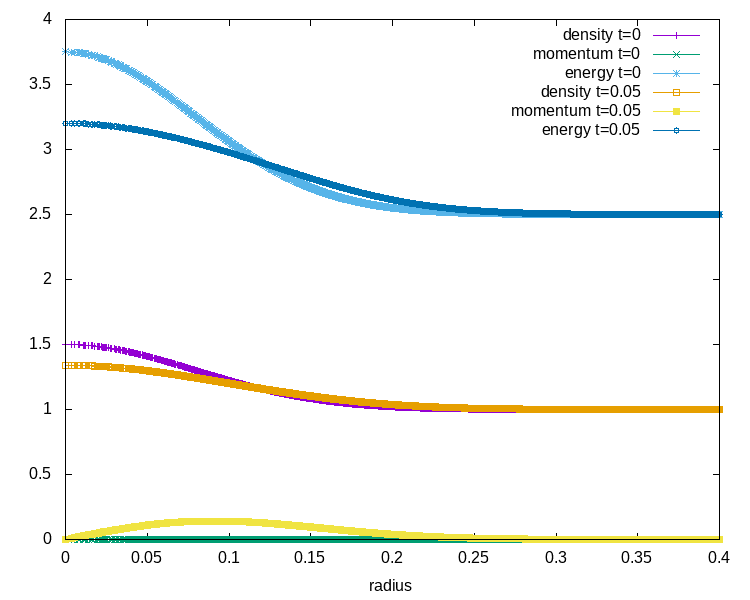}\hfill \includegraphics[width=.49\textwidth]{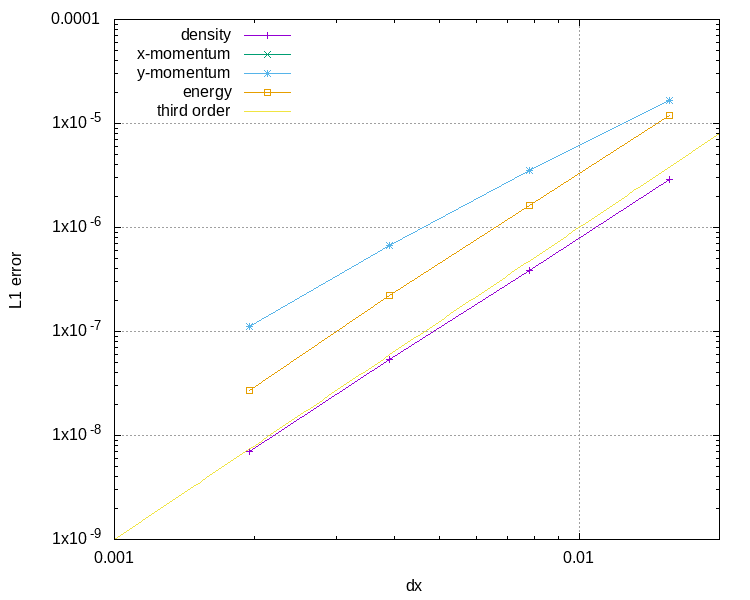} \\
 \includegraphics[width=.49\textwidth]{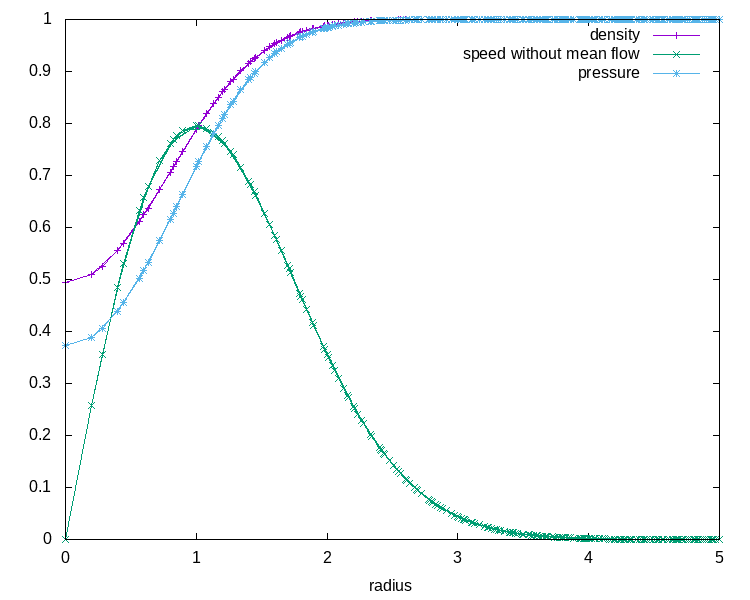}\hfill \includegraphics[width=.49\textwidth]{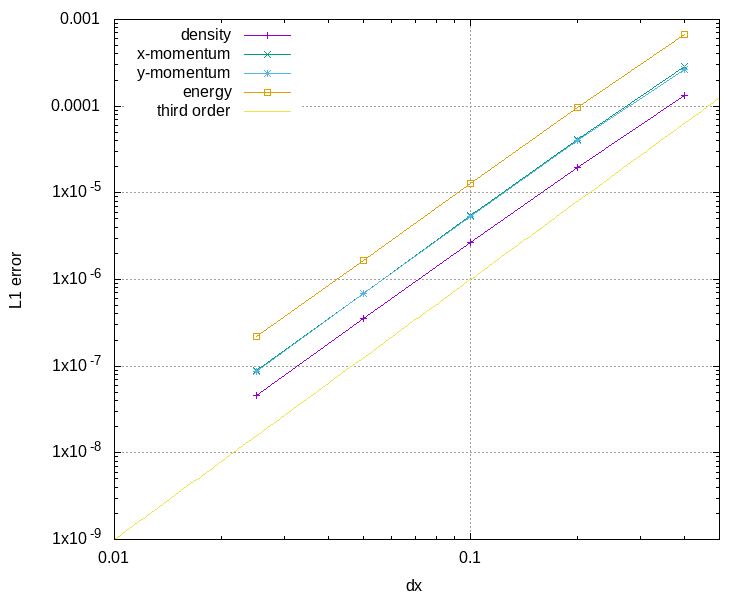} 
 \caption{Convergence study. \emph{Top left}: Setup \eqref{eq:convergence1u}--\eqref{eq:convergence1p} at initial time and at $t=0.05$, shown as a scatter plot as a function of radius, computed on a $256\times 256$ grid. \emph{Bottom left}: Setup \eqref{eq:convergence2u}--\eqref{eq:convergence2p} at initial time on a grid of $100 \times 100$, shown as a scatter plot. \emph{Right top}: $\ell^1$ error of the numerical solution of the point values, i.e. the average of $\frac{1}{|\Omega|} \sum_{ij} |q_{i+\frac12,j+\frac12}(t^n) - q_\text{ref}(t^n, x_{i+\frac12}, y_{j+\frac12})| \Delta x \Delta y$ for the setup \eqref{eq:convergence1u}--\eqref{eq:convergence1p}. The analogous error of the averages has virtually the same values and is not shown. \emph{Right bottom}: Same for the setup \eqref{eq:convergence2u}--\eqref{eq:convergence2p}.}
 \label{fig:convergence}
\end{figure}

For a convergence analysis, the following initial data (similar to those used in \cite{kerkmann18,barsukow19activeflux}) are solved until $t=0.05$ on grids of different resolution:
\begin{align}
 u_0(x,y) &= v_0(x,y) = 0 \qquad r := \sqrt{(x - \frac12)^2 + (y - \frac12)^2} \label{eq:convergence1u}\\
 \rho_0(x,y) &= p_0(x,y) = 1 + \frac12 \exp\left(-80 r^2\right) \label{eq:convergence1p}
\end{align}
Figure \ref{fig:convergence} shows the setup and the error, computed with respect to a reference solution obtained on a grid of $1024 \times 1024$ covering $[0,1]^2$. 

To facilitate comparison with other methods, we have also performed a convergence analysis using the setup of an isentropic ($p = \rho^{\gamma}$), traveling vortex as in \cite{hu99}, Example 3.3:
\begin{align}
\vecc{u_0(x, y)}{v_0(x, y)} &= \vecc{1}{1} + \frac{\Gamma }{2\pi}\exp\left(\frac{1 - r^2}2 \right) \vecc{-y}{x} \label{eq:convergence2u}\\
T_0(x,y) := \frac{p_0(x,y)}{\rho_0(x,y)} &= 1 - \frac{(\gamma-1) \Gamma^2}{8 \gamma \pi^2} \exp(1 - r^2) \label{eq:convergence2p}\\ 
r &:= \sqrt{(x - 10)^2 + (y - 10)^2}
\end{align}
The exact solution is pure advection with speed $(1, 1)$. We use $\Gamma = 5$ and compute the errors at time $t=2$. This vortex is not compactly supported, but its rotational velocity decays rapidly. We choose a larger domain than \cite{hu99}, namely $[0,20]^2$, to ensure that the rotational velocity of the vortex has decreased to the level of machine error at the boundary of the domain. 

In both cases limiting is not used and we use a CFL number of 0.2. In Figure \ref{fig:convergence} one observes third order accuracy in agreement with the expectation.

\subsection{Spherical shock tube}

\begin{figure}
 \centering
 \includegraphics[width=0.49\textwidth]{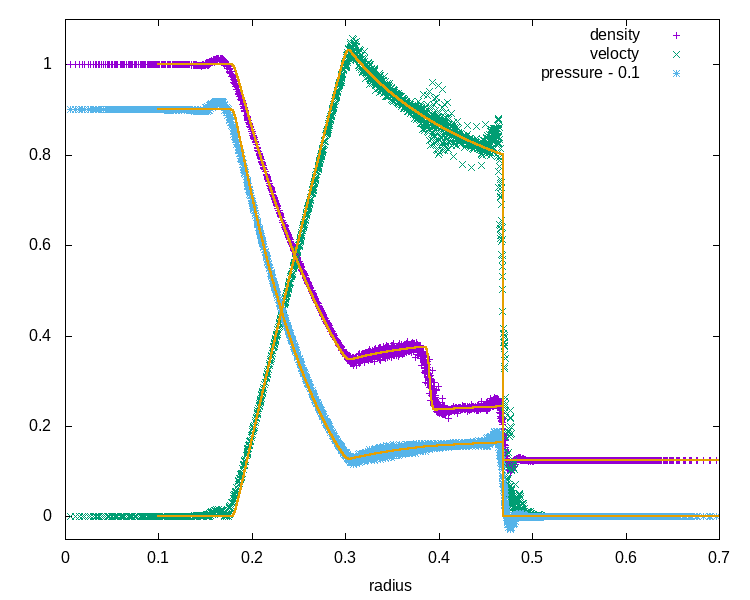} \hfill  \includegraphics[width=0.49\textwidth]{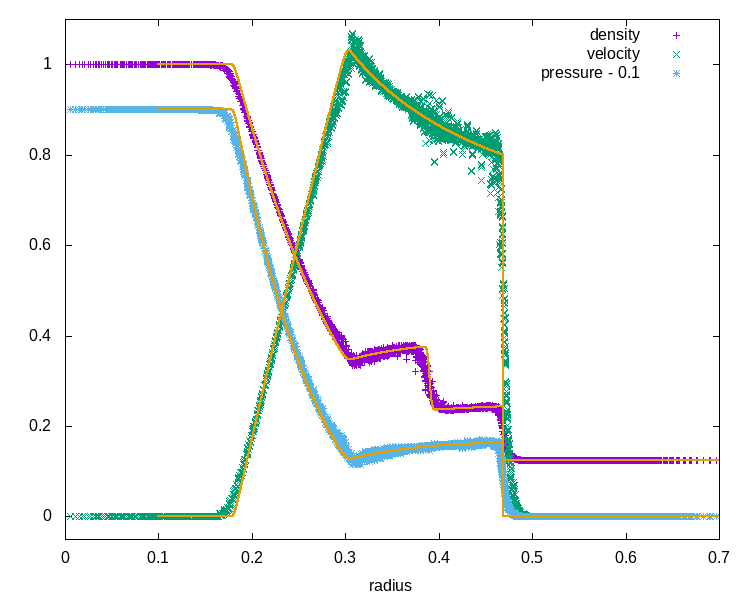}
 \caption{Radial scatter plot of the two-dimensional version of Sod's shock tube solved on a $100 \times 100$ grid. The solid line shows a finely resolved solution of the one-dimensional, radial Euler equations obtained with a standard Finite Volume method. We show the pressure offset by 0.1 for better readability of the plot. \emph{Left}: No limiting. \emph{Right}: Limiting used. The limiting is successful at suppressing oscillations at the shock and around the rarefaction. As in a radial scatter plot points from different locations end up shown in the same location, around the contact wave one rather observes scatter than oscillations, i.e. a deviation from radial symmetry.}
 \label{fig:sod}
\end{figure}

As a first test with discontinuities, Figure \ref{fig:sod} shows a 2-dimensional version of Sod's shock tube:
\begin{align}
 \rho_0(x, y) &= \begin{cases} 1 & r < 0.3 \\ 0.125 & \text{else} \end{cases} &
 p_0(x, y) &= \begin{cases} 1 & r < 0.3 \\ 0.1 & \text{else} \end{cases}  \\
 u_0(x, y) &= v_0(x,y) = 0
\end{align}
with $r = \sqrt{(x-\frac12)^2 + (y-\frac12)^2}$. We use a CFL number of 0.05, a reduction we find helpful to avoid negative pressure. One observes that the limiting is successful in suppressing oscillations at the shock and around the rarefaction. However, both the limited and the unlimited versions show some scatter around the contact wave, the origin of which will be studied as future work. Global continuity does not impede Active Flux from converging to weak solutions, because the update of the averages is conservative and fulfills a version of the Lax-Wendroff-theorem (\cite{abgrall20}).

\subsection{Multi-dimensional Riemann problems}

\begin{figure}
 \centering
 \includegraphics[width=0.49\textwidth]{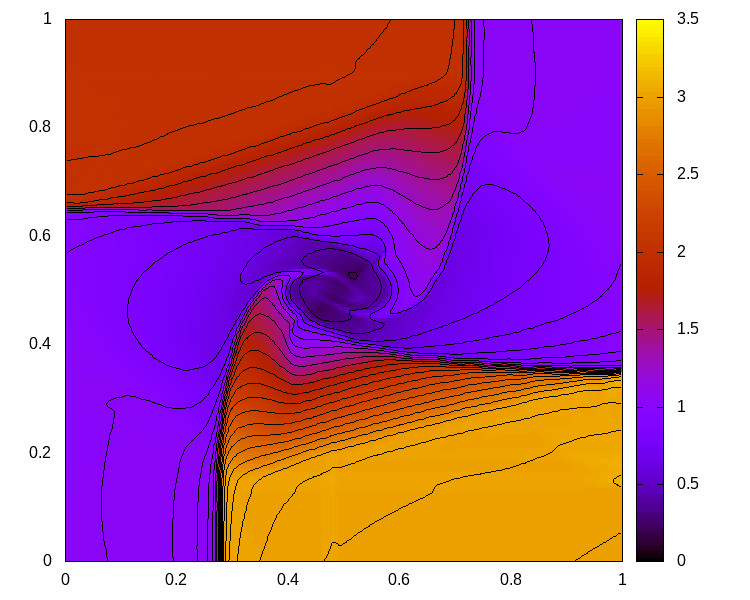} \hfill \includegraphics[width=0.49\textwidth]{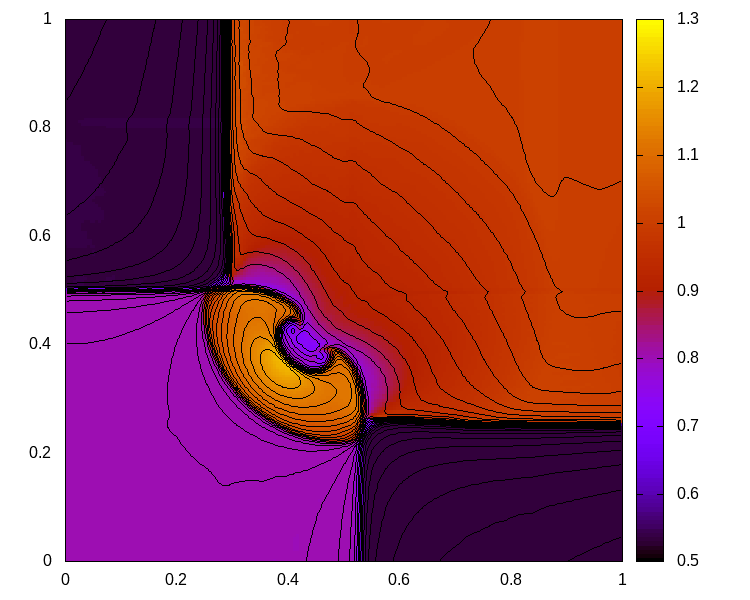}
 \includegraphics[width=0.49\textwidth]{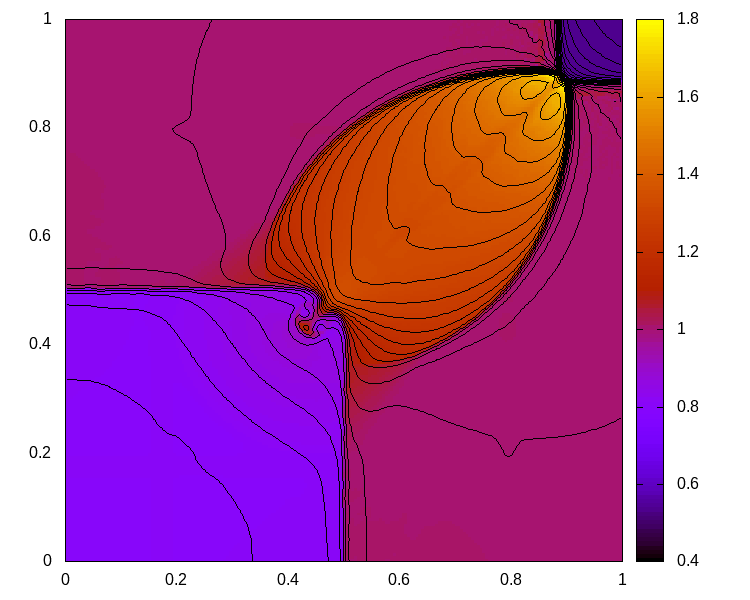} \hfill \includegraphics[width=0.49\textwidth]{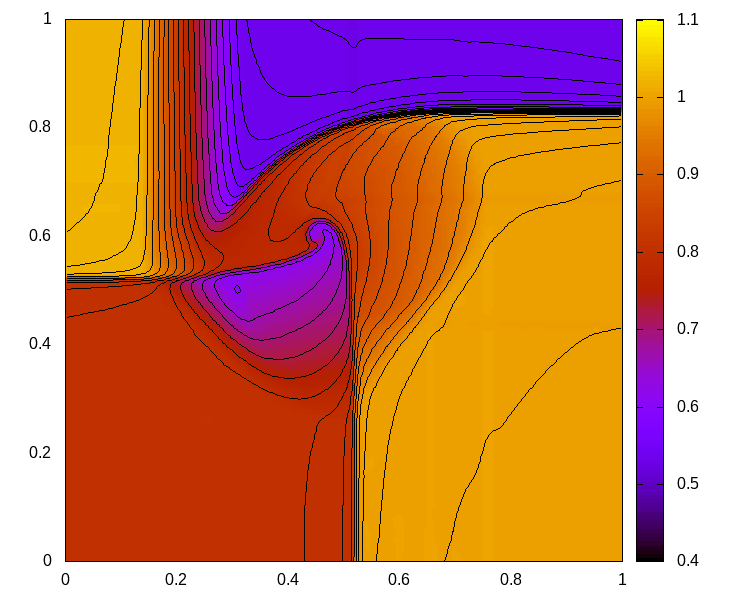}
 \caption{Multi-dimensional Riemann problems solved on a grid with $\Delta x = \Delta y = \frac{1}{200}$ using limiting as described in Section \ref{sec:limiting}. Configurations 6 (\emph{top left}), 11 (\emph{top right}), 12 (\emph{bottom left}) and 16 (\emph{bottom right}) from \cite{lax98} are shown. Density is shown in color and contour.}
 \label{fig:laxliu}
\end{figure}

\begin{figure}
 \centering
 \includegraphics[width=0.49\textwidth]{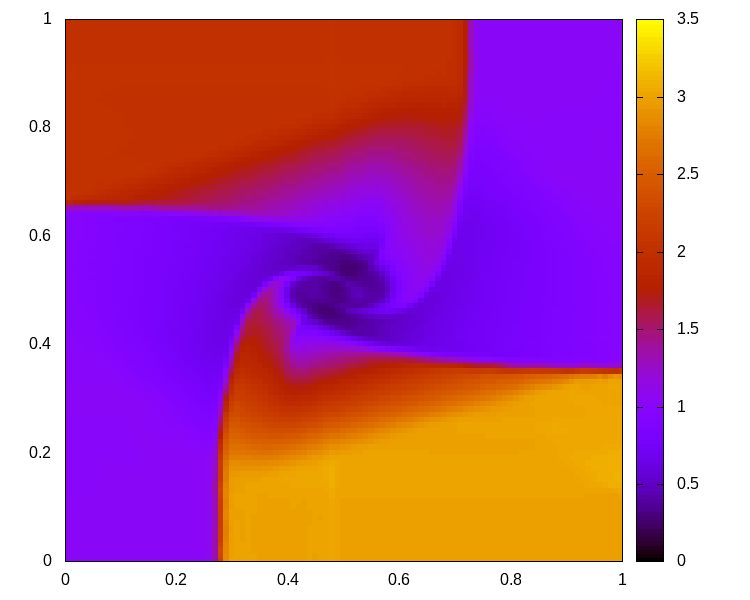} \hfill \includegraphics[width=0.49\textwidth]{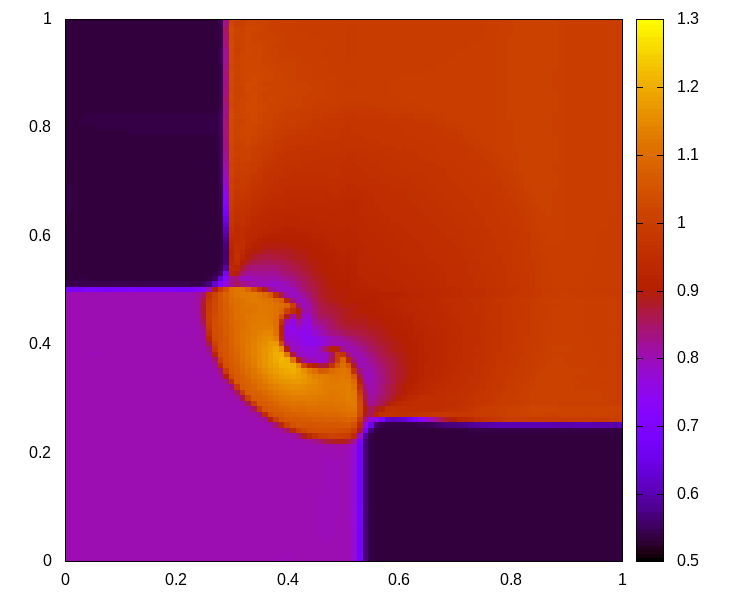}
 \includegraphics[width=0.49\textwidth]{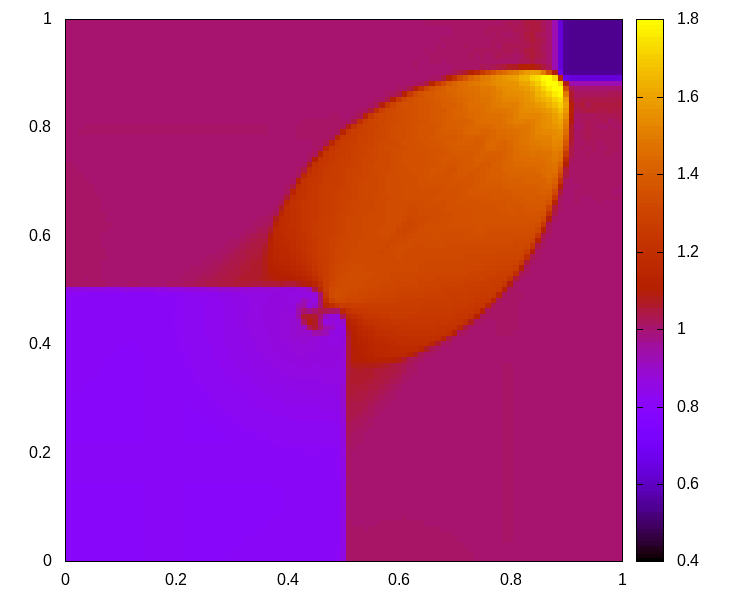} \hfill \includegraphics[width=0.49\textwidth]{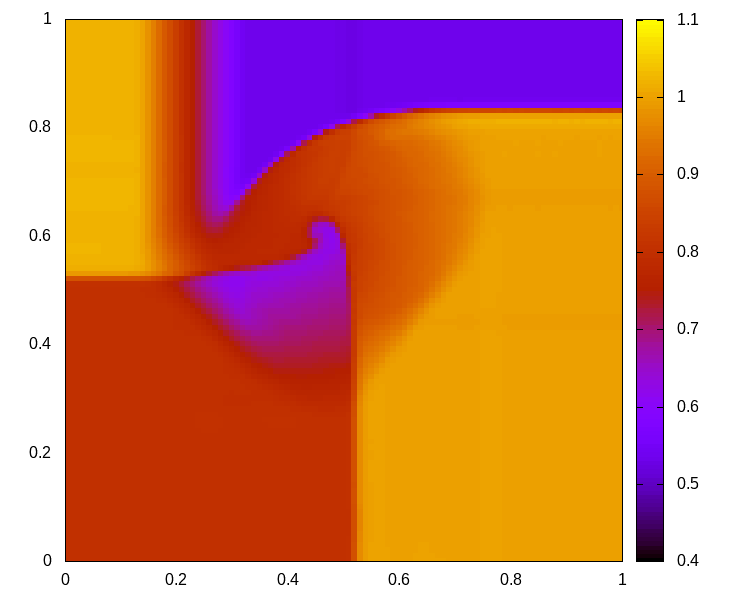}
 \caption{Setup as in Figure \ref{fig:laxliu}, but solved on a grid with $\Delta x = \Delta y = \frac{1}{100}$.}
 \label{fig:laxliucoarse}
\end{figure}

\begin{figure}
 \centering
 \includegraphics[width=0.49\textwidth]{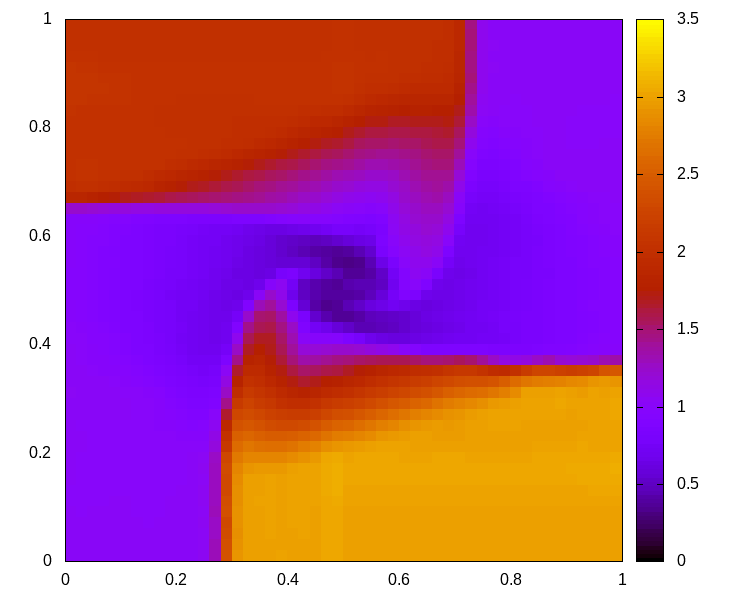} \hfill \includegraphics[width=0.49\textwidth]{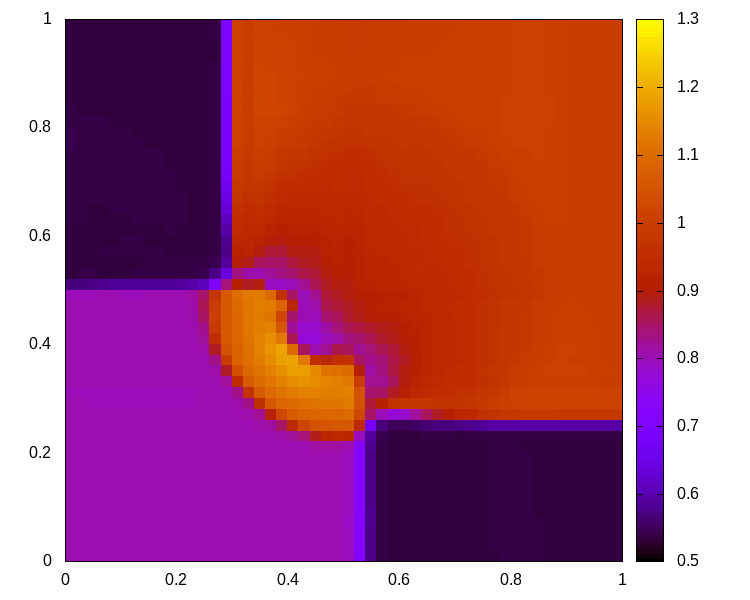}
 \includegraphics[width=0.49\textwidth]{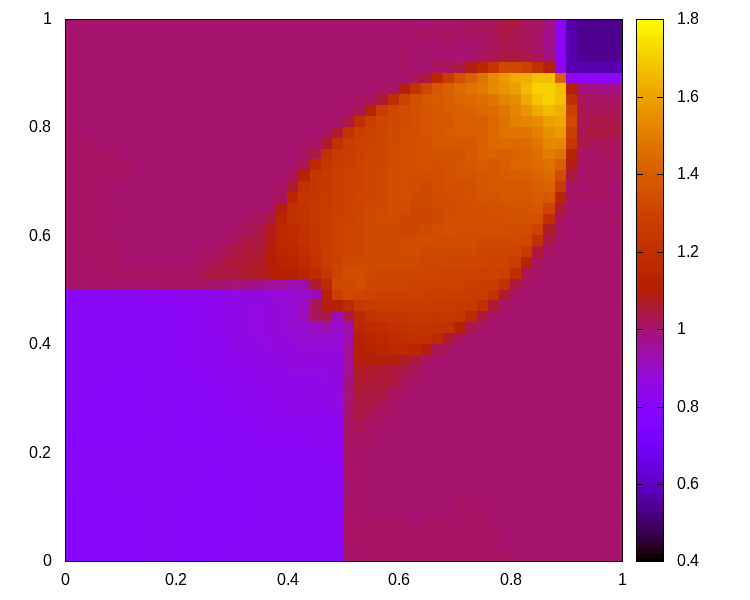} \hfill \includegraphics[width=0.49\textwidth]{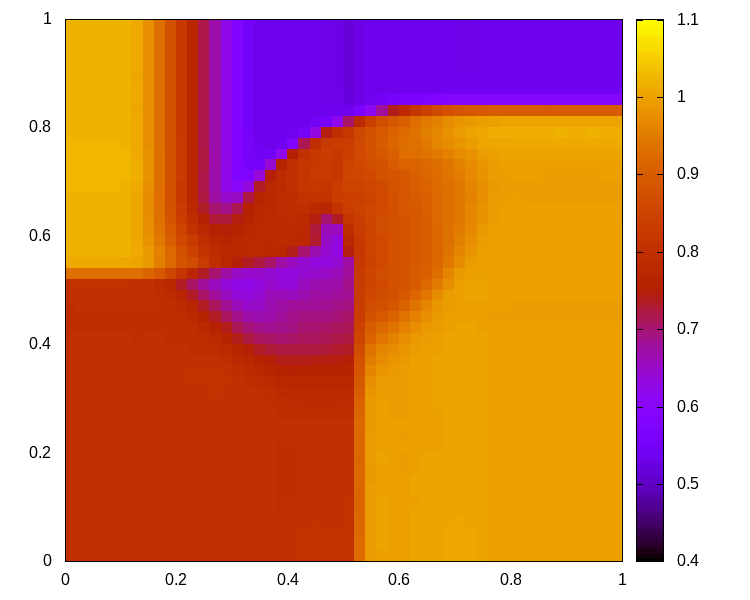}
 \caption{Setup as in Figure \ref{fig:laxliu}, but solved on a grid with $\Delta x = \Delta y = \frac{1}{50}$.}
 \label{fig:laxliucoarse2}
\end{figure}

\begin{figure}
 \centering
 \includegraphics[width=0.49\textwidth]{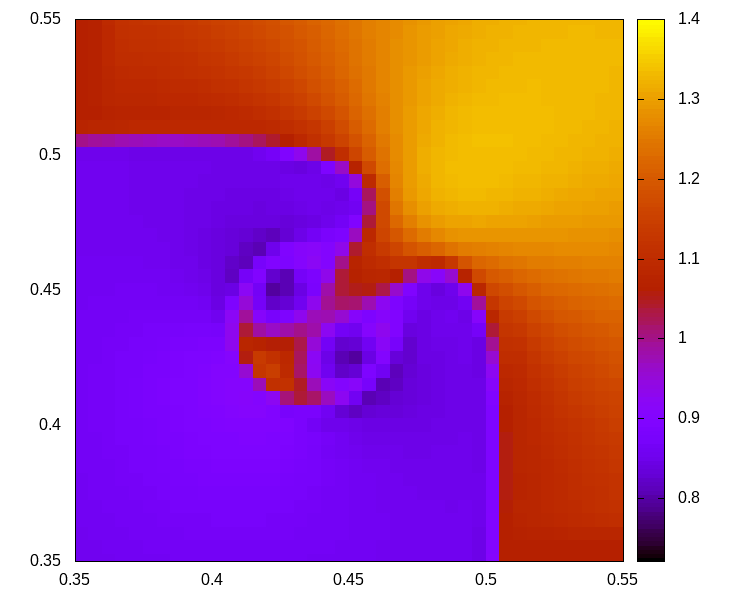} \hfill \includegraphics[width=0.49\textwidth]{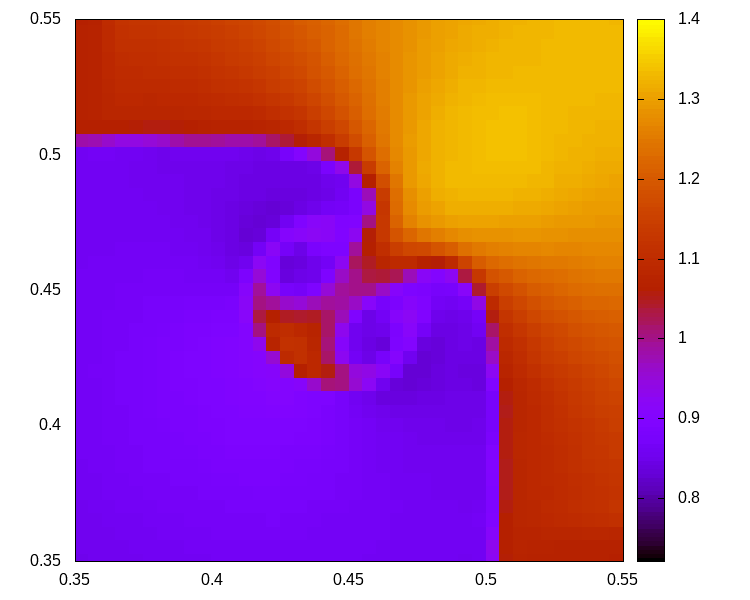}
 \caption{Influence of limiting on the central region in Configuration 12. \emph{Left}: Limiting off. \emph{Right}: Limiting on. Without limiting one observes some undershoots in the vicinity of the vortices. The structure of the solution feature is, however, not degraded by applying the limiter.}
 \label{fig:laxliulimdiff1}
\end{figure}

\begin{figure}
 \centering
 \includegraphics[width=\textwidth]{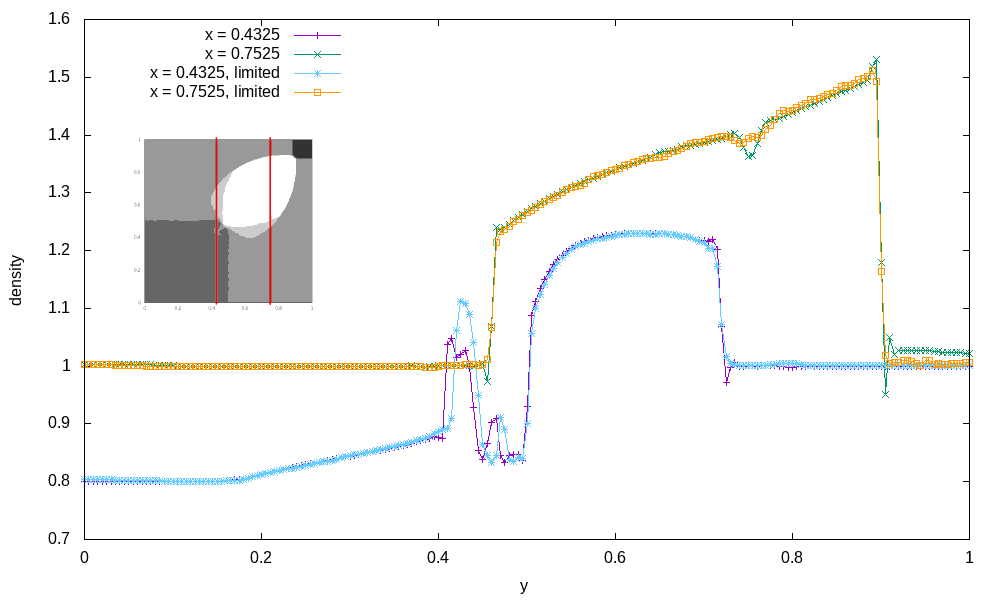} 
 \caption{Influence of limiting on Configuration 12. Density is shown along the lines $x = 0.4325$ and $x = 0.7525$ (as indicated in the inset). One observes that limiting successfully removes spurious oscillations in the vicinity of discontinuities. However, it also gently shifts the location of the central double-vortex and smears out the feature along the $x=y$ diagonal in the first quadrant.}
 \label{fig:laxliulimdiff2}
\end{figure}

In \cite{lax98}, particular multi-dimensional Riemann problems were studied, designed such that the one-dimensional Riemann problems outside the central interaction region result in elementary waves. Inside the interaction region these Riemann problems display a lot of sophisticated structure. They will illustrate the ability of the proposed method to solve complex interactions of shocks, rarefactions and slip lines. All the Riemann problems shown in Figure \ref{fig:laxliu} are solved on grids with $\Delta x = \Delta y = \frac{1}{200}$ (the original publication used $\frac{1}{400}$) with a domain slightly larger than the one shown (to exclude the influence of boundary conditions). A CFL number of 0.05 was used, as well as limiting, as described in Section \ref{sec:limiting}. Figures \ref{fig:laxliucoarse}--\ref{fig:laxliucoarse2} show results on even coarser meshes. Figures \ref{fig:laxliulimdiff1}--\ref{fig:laxliulimdiff2} show a comparison between results obtained with and without limiting. It seems that the intricate structures in the interaction region are not significantly smeared out by the limiting while oscillations at shocks are very efficiently suppressed.

\subsection{Low Mach number vortex}

\begin{figure}
 \centering
 \includegraphics[width=0.32\textwidth]{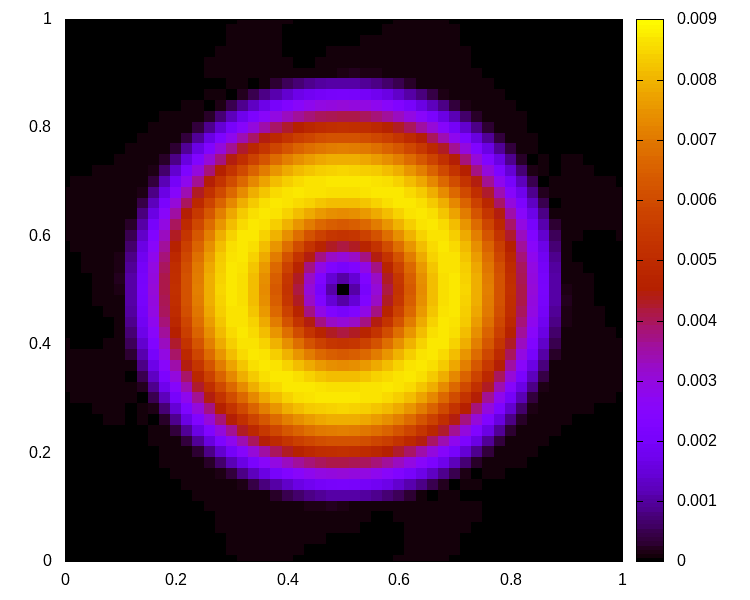} \hfill \includegraphics[width=0.32\textwidth]{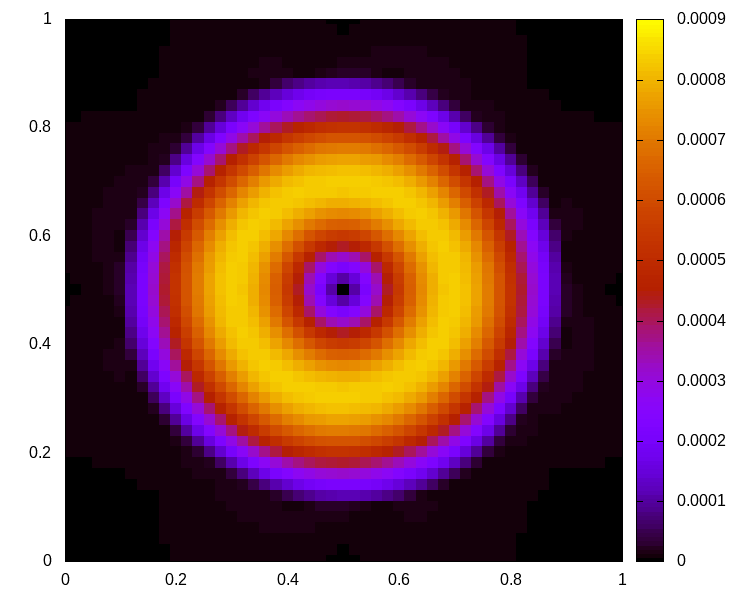} \hfill
 \includegraphics[width=0.32\textwidth]{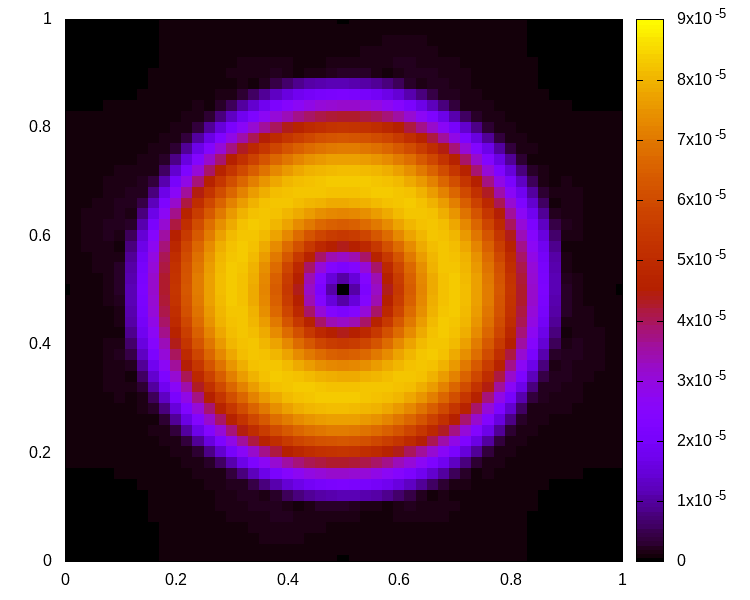} \\
 \includegraphics[width=0.32\textwidth]{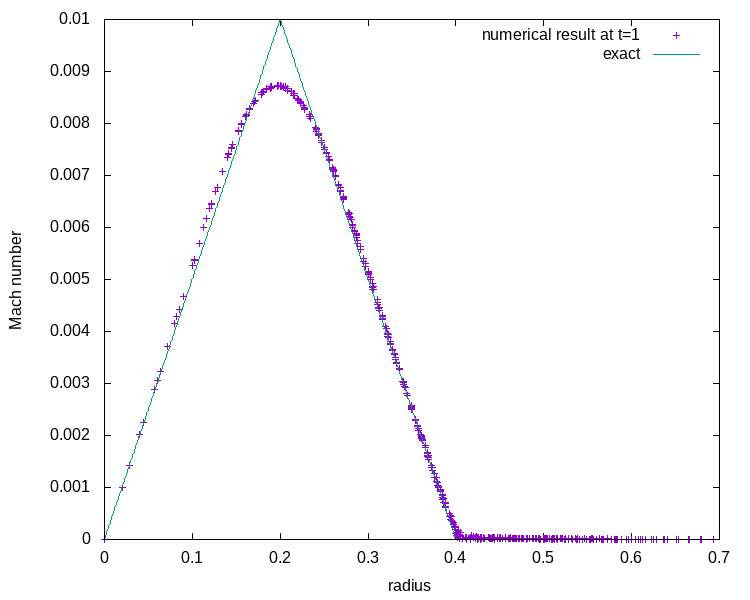} \hfill \includegraphics[width=0.32\textwidth]{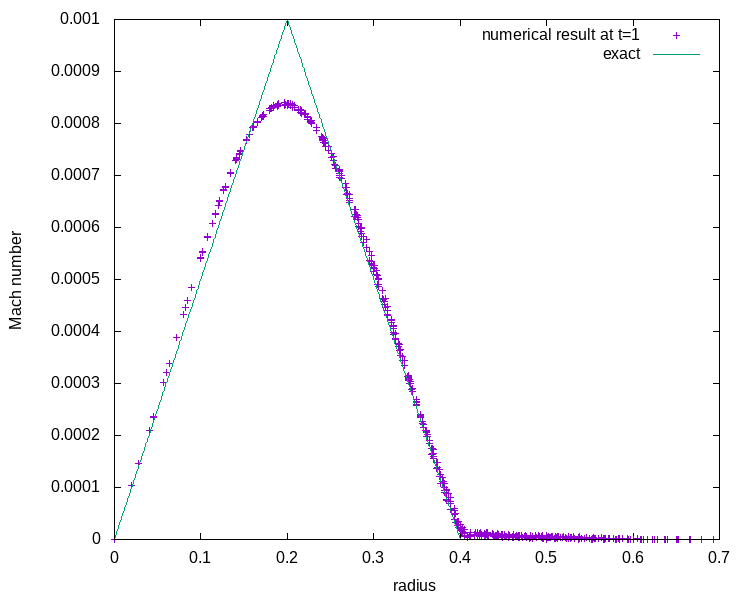} \hfill \includegraphics[width=0.32\textwidth]{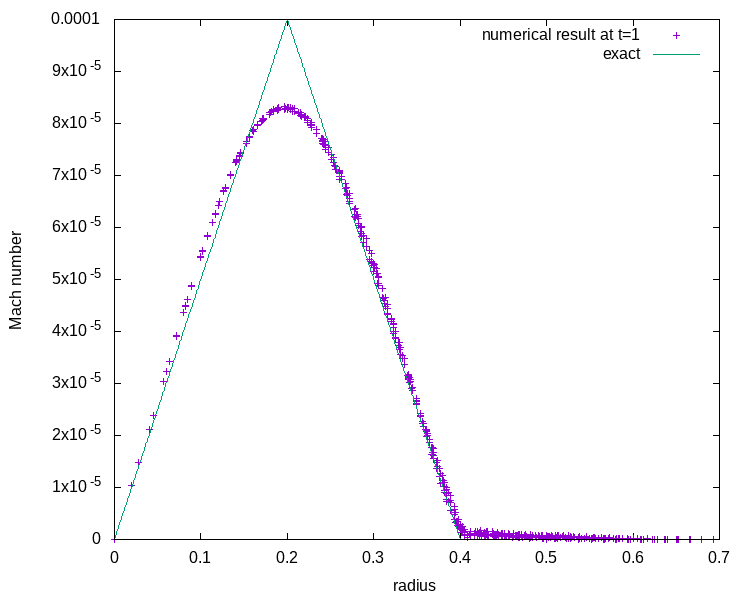}
 \caption{The Gresho vortex at various maximal Mach numbers ($10^{-2}$, $10^{-3}$, $10^{-4}$ \emph{from left to right}). While the Mach number decreases by a factor of 10 from one plot to the other, the effect of numerical diffusion remains the same asymptotically, i.e. the observed diffusion is purely advective, not acoustic.}
 \label{fig:gresho}
\end{figure}

In e.g. \cite{barsukow16} a subsonic vortex from \cite{gresho90} has been used to assess the low Mach number properties of the numerical method. The velocity field associated with this stationary solution of the Euler equations is divergence-free. The Mach number of the flow is modified by changing the background pressure, see \cite{barsukow16} for setup details. While numerical methods that are not low Mach compliant diffuse the vortex on acoustic time scales, low Mach number compliant methods apply only advective diffusion, i.e. the rate of decay of the numerical solution is asymptotically independent of the Mach number. Figure \ref{fig:gresho} shows the vortex on a $50 \times 50$ grid for different values of the Mach number until $t=1$ with a CFL number of 0.2. One observes the independence of the results on the Mach number leading us to the conclusion that Active Flux is well-suited for subsonic/low Mach number flows. This parallels the finding in \cite{barsukow18activeflux} where is was shown that Active Flux (upon usage of the exact evolution operator) is stationarity preserving and thus low Mach compliant for linear acoustics. We also observe that the advective diffusion leads to a (Mach number independent) smearing of the vortex, to levels that would require solving this setup on a $200 \times 200$ grid with a first-order low Mach number compliant method (\cite{barsukow20cgk}).

\begin{figure}
 \centering
 \includegraphics[width=0.32\textwidth]{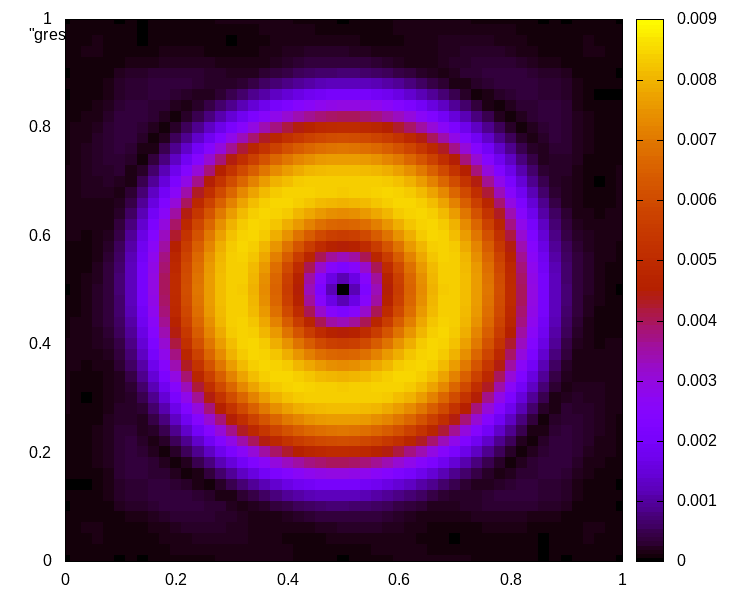} \hfill \includegraphics[width=0.32\textwidth]{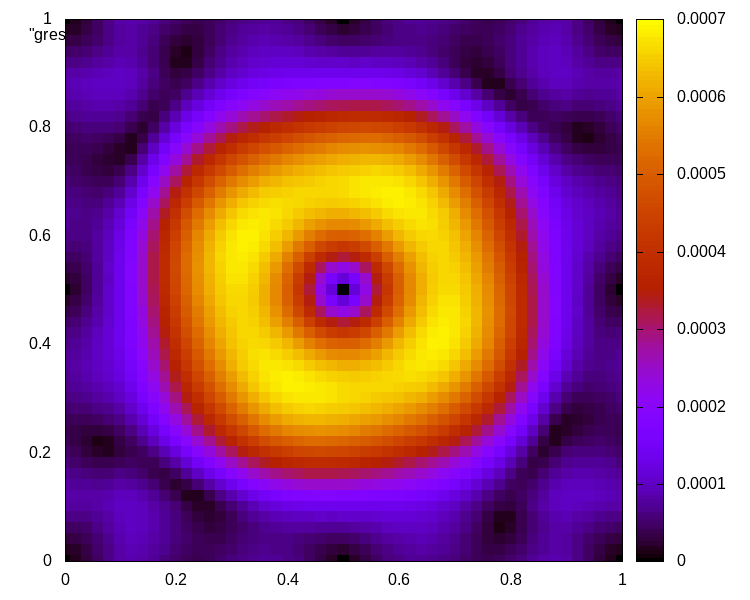} \hfill  \includegraphics[width=0.32\textwidth]{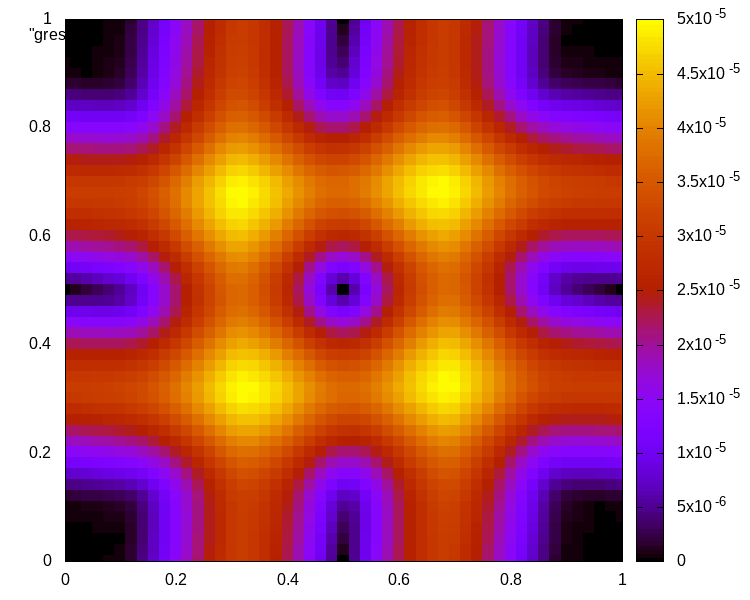} \hfill \\
 \includegraphics[width=0.32\textwidth]{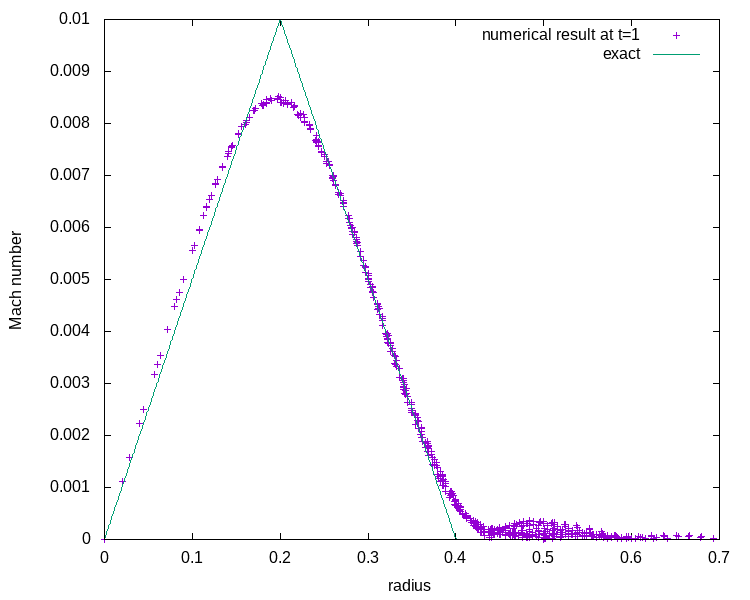} \hfill \includegraphics[width=0.32\textwidth]{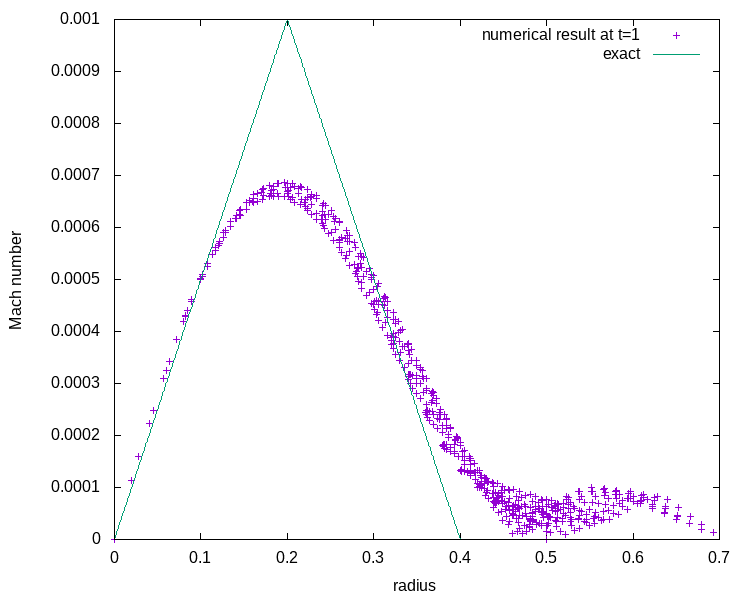} \hfill \includegraphics[width=0.32\textwidth]{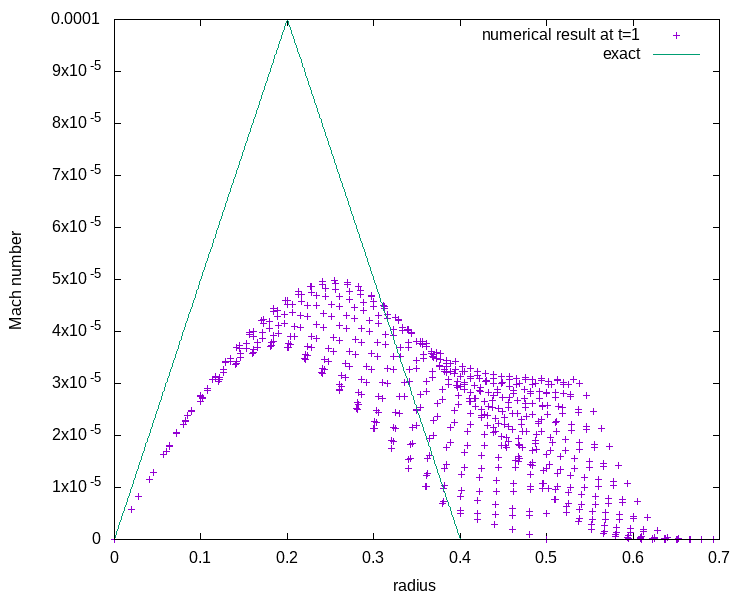}
 \caption{The Gresho vortex as in Figure \ref{fig:gresho} solved using \eqref{eq:jacobiansplitrusanovplus}--\eqref{eq:jacobiansplitrusanovminus} instead of \eqref{eq:jacobiansplitupwindplus}--\eqref{eq:jacobiansplitupwindminus}. The presence of numerical diffusion in the terms associated to acoustics makes the choice \eqref{eq:jacobiansplitupwindplus}--\eqref{eq:jacobiansplitupwindminus} unsuitable for low Mach number computations.}
 \label{fig:greshorusanov}
\end{figure}

The good low Mach number behaviour of the proposed method is due to the way how the point values are evolved, and thus not trivial.
For comparison, in Figure \ref{fig:greshorusanov} the results of a simulation are shown which uses
\begin{align}
 (J_x)^+ &:= \frac12 (J_x + s \id) \label{eq:jacobiansplitrusanovplus}\\
 (J_x)^- &:= \frac12(J_x - s \id) \qquad s := \max_k |\lambda_k| \label{eq:jacobiansplitrusanovminus}
\end{align}
instead of \eqref{eq:jacobiansplitupwindplus}--\eqref{eq:jacobiansplitupwindminus}. One observes artefacts typical of numerical diffusion associated to the acoustic terms, and thus a method unsuitable for the low Mach number regime.

\subsection{Kelvin-Helmholtz instability}

A special kind of a Kelvin-Helmholtz instability triggered by the passage of an acoustic wave has been used in \cite{munz03} to assess the properties of a numerical method for subsonic flow. This setup is of interest because acoustic phenomena are present on top of a low Mach number flow. On a domain $\left[-\frac1{\mathcal M}, \frac1{\mathcal M}\right] \times \left[0, \frac{2}{5 \mathcal M}\right]$, the initial data consist of those of a right-running sound wave $q^\text{s}$ on top of a vertically stratified background $q^\text{bg}$
\begin{align}
 q_0(x, y) = q^\text{bg}(y) + q^\text{s}(x)
\end{align}
with 
\begin{align}
    \rho_0^\text{bg}(y) &= 1 + \phi(y) &\rho_0^\text{s}(x) &= \frac{\mathcal M}5 \psi(x) \\
    u_0^\text{bg}(y) &= 0 & u_0^\text{s}(x) &= \sqrt{\gamma}  \psi(x) \\
    v_0^\text{bg}(y) &= 0 & p_0^\text{s}(x) &= \frac{1}{\mathcal M} \gamma  \psi(x)\\
    p_0^\text{bg}(y) &= \frac{1}{\mathcal M^2}   &v_0^\text{s}(x) &= 0
\end{align}
and
\begin{align}
 \phi(y) &:= \begin{cases}
            2\mathcal M y & y < 4 \\
            2\mathcal M(y - 4) - 0.4 & \text{else}
           \end{cases} &
 \psi(x) &:= 1 + \cos(\pi \mathcal M x)
\end{align}

For the linearized Euler equations
\begin{align}
 \del_t \rho^\text{s}(t,x) + \bar \rho \del_x u^\text{s}(t,x) &= 0\\
 \del_t u^\text{s}(t,x) + \frac{1}{\bar \rho} \del_x p^\text{s}(t,x) &= 0\\
 \del_t p^\text{s}(t,x) + \bar \rho c^2 \del_x u^\text{s}(t,x) &= 0
\end{align}
with $c^2 = \frac{5 \gamma}{\mathcal M^2}$ and $\bar \rho = \frac{1}{\sqrt{5}}$ the initial data $q_0^\text{s}(x)$ evolve as follows
\begin{align}
 \rho^\text{s}(t,x) &= \rho_0^\text{s}(x - ct) & u^\text{s}(t,x) &= u_0^\text{s}(x - ct) & p^\text{s}(t,x) &= p_0^\text{s}(x - ct).
\end{align}
This justifies referring to $q_0^\text{s}$ as the initial data of a right-running sound wave. The non-linearity of the full Euler equations leads to its self-steepening over time. 

\begin{figure}
 \centering
 \includegraphics[width=\textwidth]{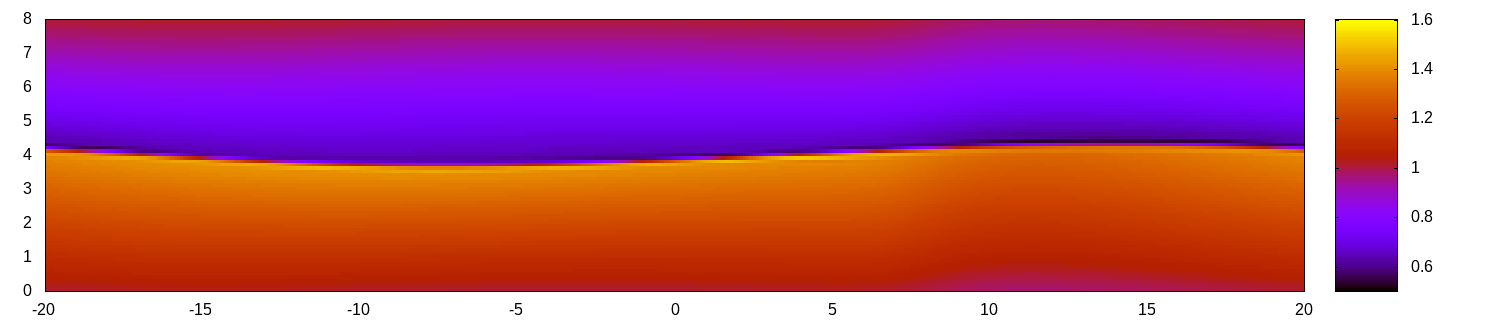}\\
 \includegraphics[width=\textwidth]{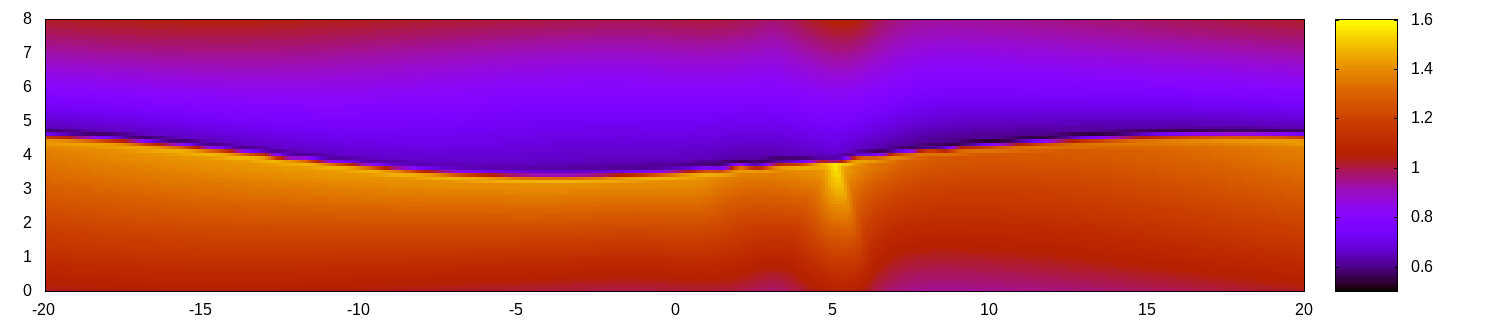}\\
 \includegraphics[width=\textwidth]{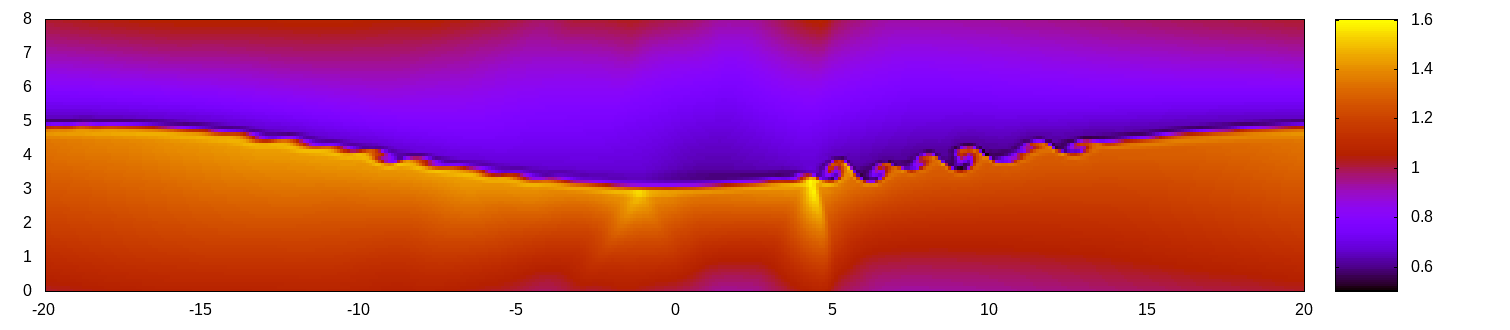}\\
 \includegraphics[width=\textwidth]{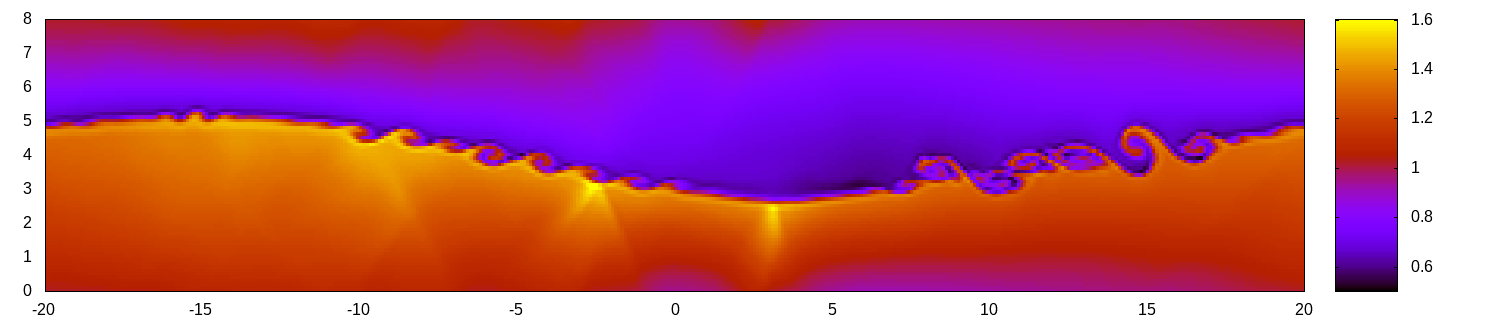}
 \caption{A Kelvin-Helmholtz instability is triggered by the passage of acoustic waves (visible as weak shocks). The setup is computed on a $400 \times 80$ grid without using limiting. Density is shown at times $t=3, 6, 9, 12$.}
 \label{fig:kh1}
\end{figure}

\begin{figure}
 \centering
 \includegraphics[width=\textwidth]{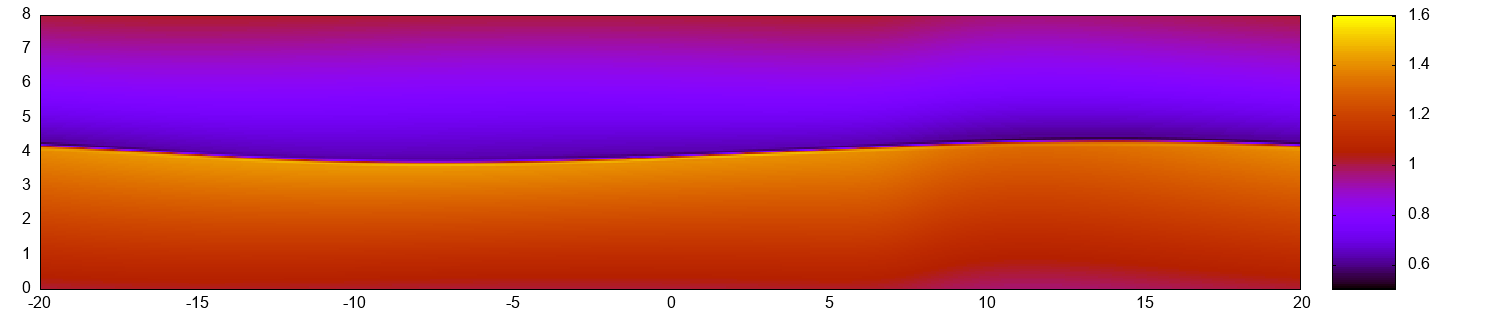}\\
 \includegraphics[width=\textwidth]{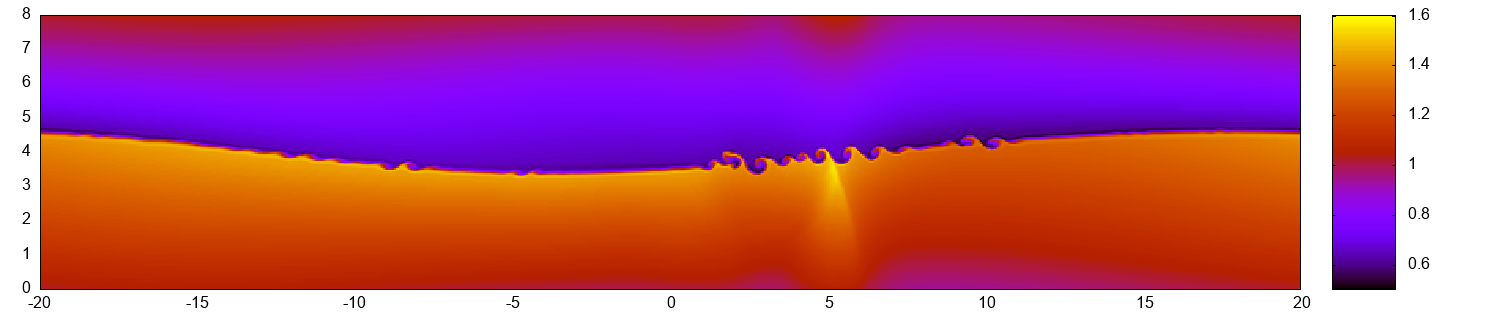}\\
 \includegraphics[width=\textwidth]{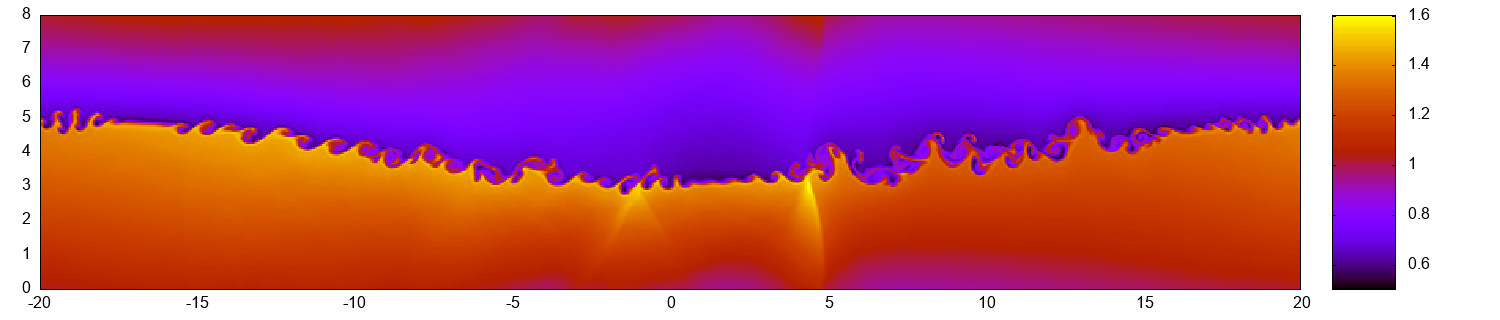}\\
 \includegraphics[width=\textwidth]{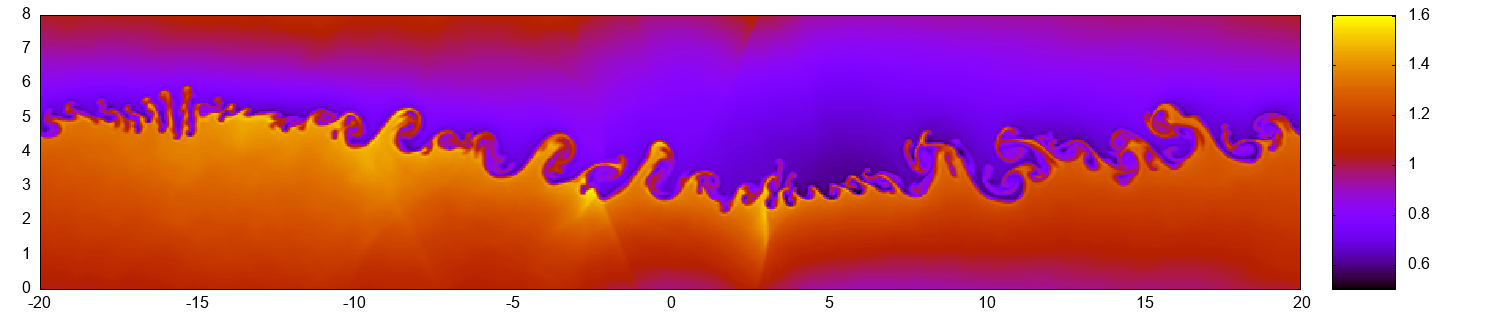}
 \caption{Same setup as Figure \ref{fig:kh1}, but on a grid of $800\times 160$.}
 \label{fig:kh2}
\end{figure}

\begin{figure}
 \centering
 \includegraphics[width=\textwidth]{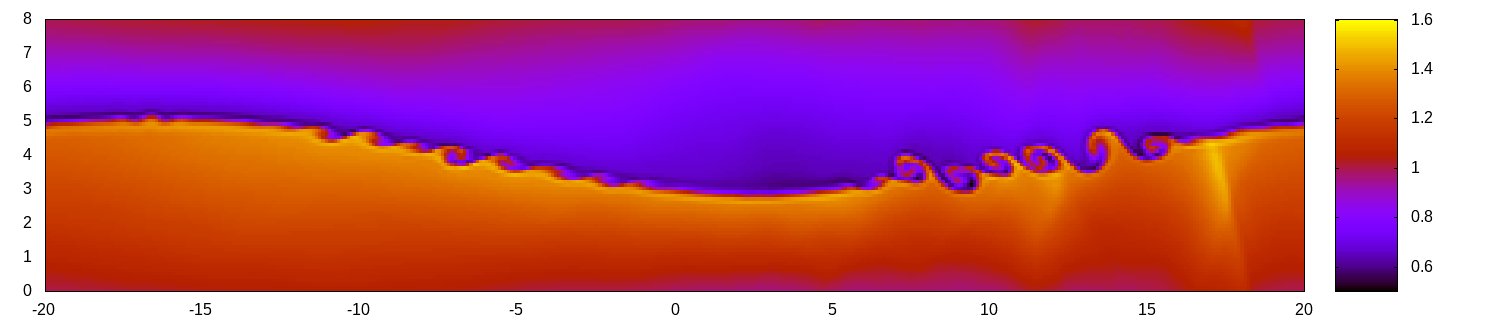}\\
 \includegraphics[width=\textwidth]{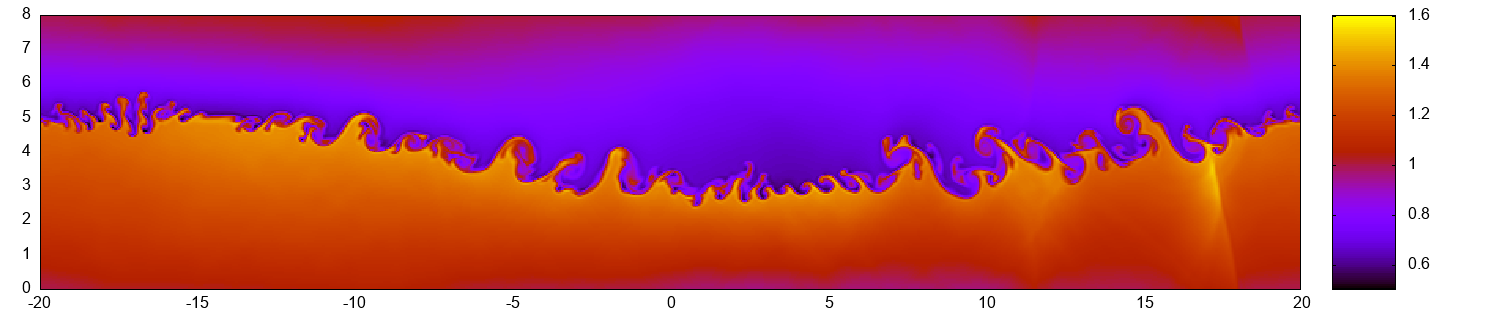}\\
 \includegraphics[width=\textwidth]{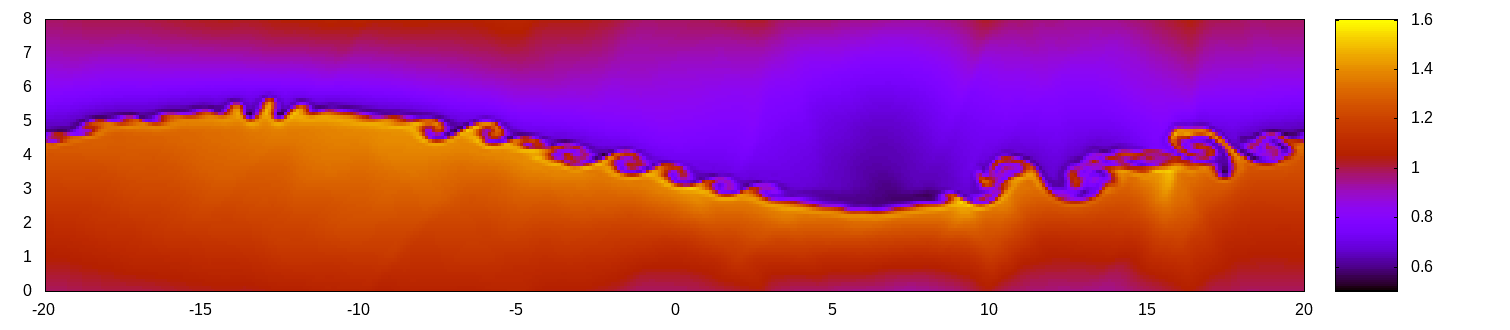}
 \includegraphics[width=\textwidth]{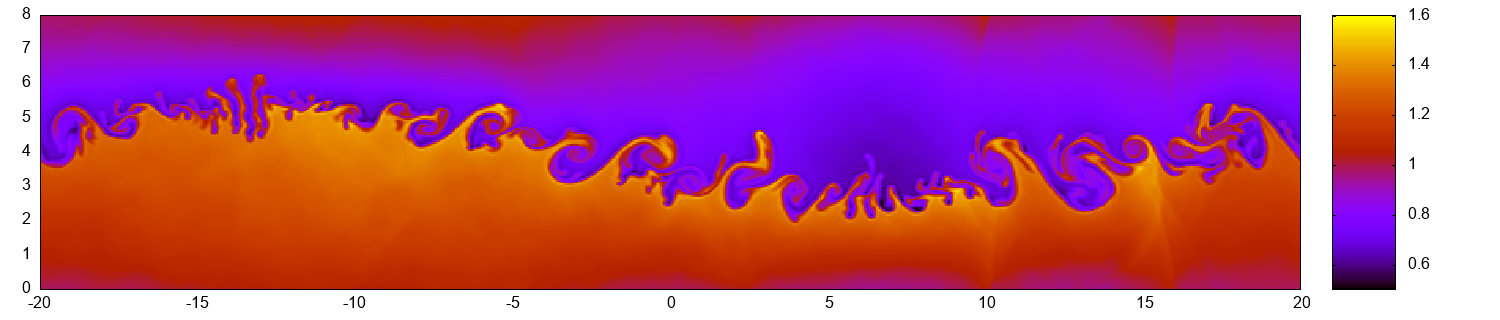}
 \caption{From top to bottom: $400 \times 80$, $t=11$; $800\times 160$, $t=11$; $400 \times 80$, $t=14$; $800\times 160$, $t=14$. Density is shown. These times are the ones shown in \cite{munz03}.}
 \label{fig:kh3}
\end{figure}

Additionally, due to the saw-tooth-shaped change in $y$-direction of the background density, the sound wave is not moving the same way above and below the interface located initially at $y=4$. A shear flow is induced, which causes a Kelvin-Helmholtz-type instability. Here we show this setup on grids of $400 \times 80$ (Figure \ref{fig:kh1}) and $800 \times 160$ (Figure \ref{fig:kh2}) with a CFL of 0.15, periodic boundaries and $\mathcal M = \frac{1}{20}$. No limiting was used. Figure \ref{fig:kh3} shows the results at additional times to facilitate comparison to \cite{munz03}. Numerical methods that are not low Mach number compliant would add so much diffusion to the subsonic background flow that the instability would be artificially stabilized (unless the grid is refined excessively). One observes that here, the method is able to adequately\footnote{without claiming that the solution is physically correct, of course, as we are solving only the Euler equations and are neglecting viscosity and other effects crucial for such setups.} evolve both the instability (comparable to results in \cite{peraire11} obtained using a high-order DG method) and the sound waves (which have become weak shocks) passing through the domain. This leads to the hypothesis that the proposed Active Flux method is well-suited for all-speed regimes.

\section{Conclusions}

Active Flux combines aspects of Finite Volume and Finite Element methods. The evolution of cell averages ensures shock-capturing properties, while the presence of point values at cell interfaces leads to a globally continuous reconstruction, which is a great difference to Godunov methods. The incorporation of additional degrees of freedom and thus the compact nature of the stencil makes the method high-order, but efficient for parallelization or the implementation of boundary conditions. The shared degrees of freedom imply less memory cost than DG methods. Finally, point values do not need to be expressed in conservative variables, i.e. Active Flux offers more freedom than conventional approaches.

For Active Flux, an emphasis on multi-dimensional thinking was put from the beginning, e.g. in \cite{eymann13}. This, again, is different from Godunov methods which are rooted in one-dimensional thinking even when applied in multiple dimensions. Interestingly, in the framework of Godunov methods, structure preservation is not obtained even if all the multi-dimensional Riemann problems are solved exactly (see \cite{barsukow17}). As has been shown in \cite{barsukow18activeflux} for linear acoustics, using an exact evolution operator for Active Flux results in a method that is stationarity preserving. 

Due to the inherent high-order nature of Active Flux, sufficiently high-order evolution operators for nonlinear, multi-dimensional systems, that can be used in fully discrete Active Flux methods are non-trivial to derive. For a semi-discrete Active Flux method (\cite{abgrall20,abgrall22}), it is easier to find spatial discretizations that use the degrees of freedom of Active Flux and to immediately write down evolution equations for the point values. The semi-discrete problem can then be integrated in time using standard methods, at the price of a reduced CFL number. The incorporation of truly multi-dimensional information is also more difficult.

To show that, however, such an Active Flux method can be successfully used to solve the multi-dimensional Euler equations is the aim of the present work. One finds that, endowed with a limiting strategy, Active Flux is indeed able to solve complex flow problems. This has been demonstrated here for examples of multi-dimensional Riemann problems and for subsonic flows, highlighting in particular the ability of Active Flux to cope with low Mach number flows. This leads us to conjecture that also the semi-discrete Active Flux method possesses structure-preserving properties. A comparison between the semi-discrete and the fully discrete, classical approach to Active Flux in the setting of linear acoustics is subject of a manuscript currently in preparation.

The semi-discrete approach can in principle be applied to other hyperbolic systems of conservation laws. Further research is necessary to understand the theoretical aspects of this method, such as entropy inequalities, or to construct appropriate boundary conditions. The presented approach to limiting is rather a proof of concept than a final choice. Future work will be directed towards better limiting procedures, a study of other choices of the point value update and towards preservation of physical conditions such as positivity of the pressure. A comparison between the classical and the semi-discrete approach for nonlinear problems, and to this end further development of high-order evolution operators for multi-dimensional nonlinear systems also need to be addressed in future.

\section*{Acknowledgement}

CK and WB acknowledge funding by the Deutsche Forschungsgemeinschaft (DFG, German Research Foundation) within \emph{SPP 2410 Hyperbolic Balance Laws in Fluid Mechanics: Complexity, Scales, Randomness (CoScaRa)}, project number 525941602.

\newcommand{\etalchar}[1]{$^{#1}$}

\appendix

\section{Detailed derivation of the multi-dimensional limiting} \label{app:recon}

\subsection{Piecewise-biparabolic reconstruction} \label{sec:piecewisebiparabolic}

If none of the edges needs to be limited, then the natural choice of the reconstruction is biparabolic, as it has been used already since \cite{barsukow18activeflux,kerkmann18}. 
However, if one of the edges is reconstructed as a hat, then something else needs to be done inside the cell in order to ensure continuity. We generally choose to subdivide the cell into regions (quadrants or halves, depending on the situation) and to reconstruct biparabolically in every such region while maintaining global continuity.

Then, if the discrete data fulfill condition \eqref{eq:quasimonotone}, the reconstruction is tested (in an approximate way) to check whether $m \leq q_\text{recon}(x,y) \leq M$ holds inside the cell. In case not, the reconstruction is discarded and replaced by the plateau reconstruction of Section \ref{sec:plateau}.

Linearity of the problem (in the point values and the average) will be exploited by considering all point values apart from $q_\text{SW}, q_\text{W}, q_\text{NW}$ to vanish:

\begin{definition}
 Consider all point values apart from $q_\text{SW}, q_\text{W}, q_\text{NW}$ to vanish. Then a reconstruction of the cell that interpolates these values pointwise and whose average agrees with $\bar q$ is called the \textbf{edge-basis-function} $q_\text{recon}^\text{W}$ of the W-edge:
 \begin{align}
  q_\text{recon}^\text{W}\left(-\frac{\Delta x}{2}, \frac{\Delta y}{2}\right) &= q_\text{NW} &
  q_\text{recon}^\text{W}\left(-\frac{\Delta x}{2}, 0\right) &= q_\text{W}\\
  q_\text{recon}^\text{W}\left(-\frac{\Delta x}{2}, -\frac{\Delta y}{2}\right) &= q_\text{SW} &
  \frac{1}{\Delta x \Delta y} \int_{-\frac{\Delta x}{2}}^{\frac{\Delta x}{2}} \int_{-\frac{\Delta y}{2}}^{\frac{\Delta y}{2}} q_\text{recon}^\text{W}(x, y) \dd x \dd y  &= \bar q
 \end{align}

 Similar notions will be used for the other edges.
\end{definition}

Observe that an edge-basis-function is a reconstruction of the entire cell. In the following, only the edge-basis-functions for the W-edge will be given explicitly, as those for the other edges can be obtained by rotation, as long as $\Delta y = \Delta x$ (otherwise some rescaling is necessary).

\begin{theorem}
If the edge-basis-function for the W-edge is
\begin{align}
 q_\text{recon}^\text{W}(q_\text{SW}, q_\text{W}, q_\text{NW}, x, y, \text{S}, \text{N}, \text{W},  \bar q )
\end{align}
then the other basis functions are
\begin{align}
    q_\text{recon}^\text{S}(q_\text{SE}, q_\text{S}, q_\text{SW}, x, y, \text{E}, \text{W}, \text{S},  \bar q ) &= q_\text{recon}^\text{W}(q_\text{SE}, q_\text{S}, q_\text{SW}, y, -x, \text{E}, \text{W}, \text{S},  \bar q )  \\
    q_\text{recon}^\text{N}(q_\text{NW}, q_\text{N}, q_\text{NE}, x, y, \text{W}, \text{E}, \text{N},  \bar q ) &= q_\text{recon}^\text{W}(q_\text{NW}, q_\text{N}, q_\text{NE}, -y, x, \text{W}, \text{E}, \text{N},  \bar q ) \\
    q_\text{recon}^\text{E}(q_\text{NE}, q_\text{E}, q_\text{SE}, x, y, \text{N}, \text{S}, \text{E},  \bar q ) &= q_\text{recon}^\text{W}(q_\text{NE}, q_\text{E}, q_\text{SE}, -x, -y, \text{N}, \text{S}, \text{E},  \bar q ) 
\end{align}
\end{theorem}

The edge-basis-function depends on $q_\text{SW}, q_\text{W}, q_\text{NW}$, on whether the reconstruction of the W-edge is parabolic or hat, and -- this complicates things a little -- on whether the neighbouring edges (S and N) are reconstructed as hats or as parabolae. This is necessary due to global continuity and because the corner values $q_\text{SW}, q_\text{NW}$ are shared with the S- and N-edges. 

The final reconstruction is obtained through summation:
\begin{theorem}\label{thm:reconinterpol}
 The following reconstruction $q_\text{recon}$ interpolates all the point values along the boundary of the cell and its average agrees with the given cell average:
 \begin{align}
  q_\text{recon}(x,y) &:= q_\text{recon}^\text{W}\left(\frac{q_\text{SW}-\bar Q}{2}, q_\text{W}-\bar Q, \frac{q_\text{NW}-\bar Q}{2}, x, y, \text{S}, \text{N}, \text{W}, \frac{ \bar q-\bar Q}{4}\right) \\&\nonumber
+q_\text{recon}^\text{S}\left(\frac{q_\text{SE}-\bar Q}{2}, q_\text{S}-\bar Q, \frac{q_\text{SW}-\bar Q}{2}, x, y, \text{E}, \text{W}, \text{S}, \frac{ \bar q-\bar Q}{4}\right)\\
&\nonumber+q_\text{recon}^\text{N}\left(\frac{q_\text{NW}-\bar Q}{2}, q_\text{N}-\bar Q, \frac{q_\text{NE}-\bar Q}{2}, x, y, \text{W}, \text{E}, \text{N}, \frac{ \bar q-\bar Q}{4}\right) \\&\nonumber
+q_\text{recon}^\text{E}\left(\frac{q_\text{NE}-\bar Q}{2}, q_\text{E}-\bar Q, \frac{q_\text{SE}-\bar Q}{2}, x, y, \text{N}, \text{S}, \text{E}, \frac{\bar q-\bar Q}{4}\right)\\
				&\nonumber+ \bar Q
 \end{align}
 where 
 $\bar Q := \frac{q_\text{SW} + q_\text{W} + q_\text{NW} + q_\text{N} + q_\text{NE} + q_\text{E} + q_\text{SE} + q_\text{S}}{8}$.
 Moreover, as all the point values tend to $\bar q$,
 \begin{align}
  q_\text{recon}(x, y) \to \bar q \label{eq:reconcontinuity}
 \end{align}
 for all $x, y$.
\end{theorem}
\begin{proof}
 The pointwise interpolation property is clear because, for example,
 \begin{align}
  q_\text{recon}\left( \frac{\Delta x}{2}, \frac{\Delta y}{2} \right) &= q_\text{recon}^\text{N}\left(\frac{q_\text{NW}-\bar Q}{2}, q_\text{N}-\bar Q, \frac{q_\text{NE}-\bar Q}{2}, \frac{\Delta x}{2}, \frac{\Delta y}{2}, \text{W}, \text{E}, \text{N}, \frac{\bar q-\bar Q}{4}\right)\\
+\nonumber&q_\text{recon}^\text{E}\left(\frac{q_\text{NE}-\bar Q}{2}, q_\text{E}-\bar Q, \frac{q_\text{SE}-\bar Q}{2}, \frac{\Delta x}{2}, \frac{\Delta y}{2}, \text{N}, \text{S}, \text{E}, \frac{\bar q-\bar Q}{4}\right) + \bar Q\\
&= \frac{q_\text{NE}-\bar Q}{2} + \frac{q_\text{NE}-\bar Q}{2}+\bar Q = q_\text{NE}
 \end{align}
 The correctness of the average follows from
 \begin{align}
  \frac{1}{\Delta x \Delta y} \int_{-\frac{\Delta x}{2}}^{\frac{\Delta x}{2}} \int_{-\frac{\Delta y}{2}}^{\frac{\Delta y}{2}} q_\text{recon}^\text{W}(x, y) \dd x \dd y = 4 \cdot \frac{\bar q-\bar Q}{4} + \bar Q = \bar q
 \end{align}
 
 Finally, property \eqref{eq:reconcontinuity} is trivial if $\bar q = 0$, because the reconstruction is linear in all the point values and in the average, and thus $q_\text{recon}(x, y) \to 0$ uniformly in this case. If the point values tend to $\bar q \neq 0$, then so does $\bar Q \to \bar q$ and thus
 \begin{align}
  q_\text{recon}(x,y) \to 0 + \bar Q   \to \bar q
 \end{align}
\end{proof} 
\begin{remark}: One might think that it would be sufficient to define the reconstruction as
\begin{align}
  &q_\text{recon}^\text{W}\left(\frac{q_\text{SW}}{2}, q_\text{W}, \frac{q_\text{NW}}{2}, x, y, \text{S}, \text{N}, \text{W}, \frac{\bar q}{4}\right)
+q_\text{recon}^\text{S}\left(\frac{q_\text{SE}}{2}, q_\text{S}, \frac{q_\text{SW}}{2}, x, y, \text{E}, \text{W}, \text{S}, \frac{\bar q}{4}\right)\\
+&q_\text{recon}^\text{N}\left(\frac{q_\text{NW}}{2}, q_\text{N}, \frac{q_\text{NE}}{2}, x, y, \text{W}, \text{E}, \text{N}, \frac{\bar q}{4}\right)
+q_\text{recon}^\text{E}\left(\frac{q_\text{NE}}{2}, q_\text{E}, \frac{q_\text{SE}}{2}, x, y, \text{N}, \text{S}, \text{E}, \frac{\bar q}{4}\right)
 \end{align}
This function also has the interpolation properties in Theorem \ref{thm:reconinterpol}. However, in the limit of all the point values converging to $\bar q$, property \eqref{eq:reconcontinuity} is not guaranteed. Linearity merely implies that in the limit, $q_\text{recon}$ will be proportional to $\bar q$, but it can still have a non-trivial dependence on $x, y$. The only case where one can be sure of obtaining a uniform constant is when $\bar q=0$. Thus, one can first subtract a uniform constant $k$ from all the discrete data, for instance the average of all the point values, reconstruct, and add it back:
\begin{align}
  &q_\text{recon}^\text{W}\left(\frac{q_\text{SW}}{2}-k, q_\text{W}-k, \frac{q_\text{NW}}{2}-k, x, y, \text{S}, \text{N}, \text{W}, \frac{\bar q-k}{4}\right)\\
&+q_\text{recon}^\text{S}\left(\frac{q_\text{SE}}{2}-k, q_\text{S}-k, \frac{q_\text{SW}}{2}-k, x, y, \text{E}, \text{W}, \text{S}, \frac{\bar q-k}{4}\right)\\
+&q_\text{recon}^\text{N}\left(\frac{q_\text{NW}}{2}-k, q_\text{N}-k, \frac{q_\text{NE}}{2}-k, x, y, \text{W}, \text{E}, \text{N}, \frac{\bar q-k}{4}\right)\\
&+q_\text{recon}^\text{E}\left(\frac{q_\text{NE}}{2}, q_\text{E}, \frac{q_\text{SE}}{2}, x, y, \text{N}, \text{S}, \text{E}, \frac{\bar q-k}{4}\right) + k
 \end{align}
 Now, if all point values tend to $\bar q$, $k$ also tends to $\bar q$, and all the reconstructions tend to uniform 0. The choice of $k = \bar Q$ is a simple one, but other ones are possible.
\end{remark}

As mentioned in Section \ref{sec:limiting}, instead of directly imposing the correct average over the cell, the value $q_\text{C}$ of the reconstruction at the cell center is imposed. Later, given the cell average, $q_\text{C}$ is found as the solution of a linear equation. This is a purely algorithmic detour not affecting the results, but it makes it easier to combine all the different cases. Below, these linear equations linking $\bar q$ and $q_\text{C}$ are given together with the reconstructions.

If the reconstruction happens on the unit square, then $\Delta x = \Delta y = 1$ should be used in the formulas below. The sketches of the interpolation problem are encoded as follows: \includegraphics[width=0.02\textwidth]{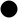} denotes the central value $q_\text{C}$, \includegraphics[width=0.02\textwidth]{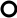} / \includegraphics[width=0.02\textwidth]{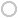} denotes a value that is not on the W edge and thus zero (gray if it is not used in the interpolation), \includegraphics[width=0.02\textwidth]{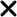} / \includegraphics[width=0.02\textwidth]{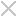} denotes one of the values $q_\text{NW}$, $q_\text{W}$, $q_\text{SW}$ (gray if it is not used in the interpolation). Values marked with an arrow do not, in principle, need to be included in the interpolation stencil, but are included here in order to achieve continuity. The colored area denotes the support of the different functions that make up the piecewise defined reconstruction. \includegraphics[width=0.02\textwidth]{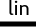} denotes an edge that is reconstructed linearly, in other words, as part of the interpolation procedure, we impose that the restriction of the reconstruction onto that edge is linear (the quadratic term vanishing).

In many cases, the reconstruction is (piecewise) biparabolic, i.e. of the form
\begin{align}
 (a_0 + a_1 x + a_2 x^2) +  (a_3  + a_4 x  + a_5 x^2 )y +(a_6  + a_7 x + a_8 x^2) y^2
\end{align}
In the following, biparabolic reconstructions are given by specifying the values of these 9 coefficients.

\subsubsection{Parabolic reconstruction on W edge} \label{ssec:para}

If edges W, S and N are all reconstructed parabolically, then the W-edge-basis-function is a biparabolic function. If either S or N (or both) are reconstructed as hat functions, the reconstruction in the cell is defined piecewise: the left and the rights halves of the cell have individual biparabolic reconstructions, which are joined in a continuous fashion.

\paragraph{Parabolic reconstruction on both neighbouring edges}

\begin{figure}
 \centering
  \includegraphics[width=0.12\textwidth]{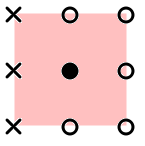}\\
 \includegraphics[width=.7\textwidth]{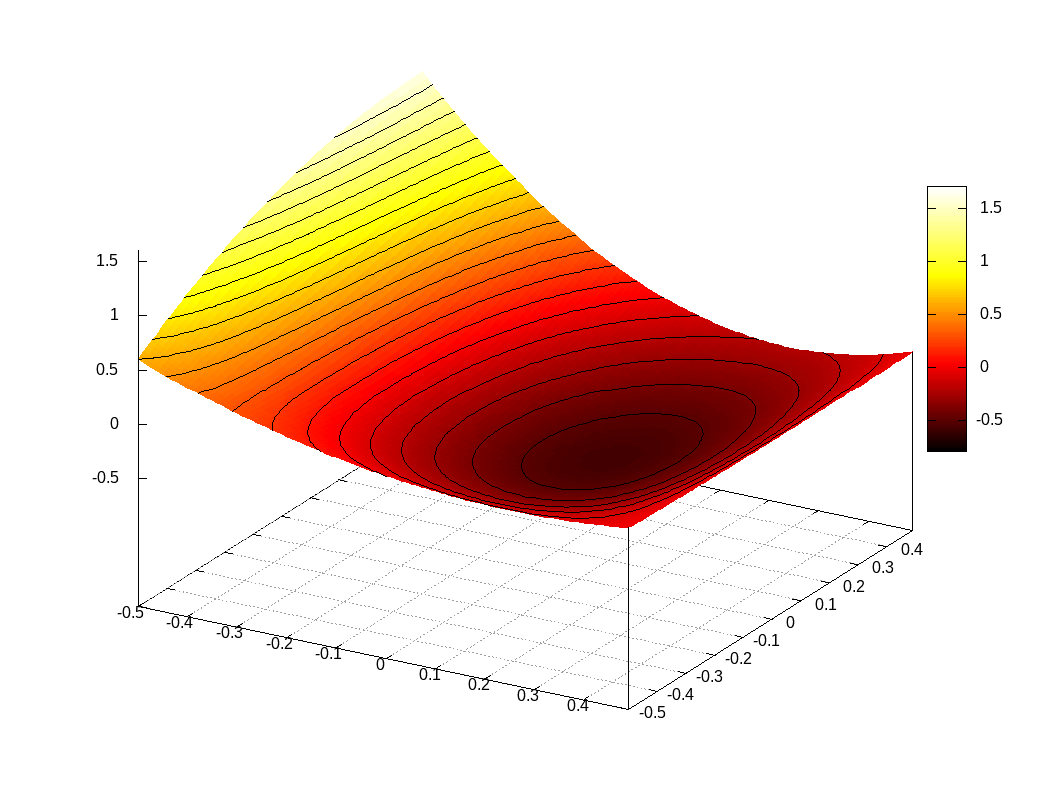}
 \caption{All edges are reconstructed parabolically, and the corresponding edge-basis-function is a simple biparabolic interpolation. $q_\text{NW} = 1.6$, $q_\text{W}=1.35$, $q_\text{SW} = 0.6$.}
 \label{fig:example-bipara}
\end{figure}

If both neighbouring edges (N and S) are reconstructed parabolically, then the reconstruction inside the cell is the trivial biparabolic reconstruction (see Figure \ref{fig:example-bipara}):
\begin{align}
q_\text{recon}^\text{W} &= \left\{a_0 = q_\text{C},a_1 = -\frac{q_\text{W}}{\Delta x},a_2 = -\frac{2 (2 q_\text{C}-q_\text{W})}{\Delta x^2},a_3 = 0,a_4 = -\frac{q_\text{NW}-q_\text{SW}}{\Delta x \Delta y},
     \right . \\ & \phantom{mmm}\left.\nonumber      
a_5 = \frac{2(q_\text{NW}-q_\text{SW})}{\Delta x^2 \Delta y},a_6 = -\frac{4 q_\text{C}}{\Delta y^2},a_7 = -\frac{2(q_\text{NW}+q_\text{SW}-2 q_\text{W})}{\Delta x \Delta y^2}, 
\right . \\ & \phantom{mmm}\left.\nonumber      
a_8 = \frac{4 (4 q_\text{C}+q_\text{NW}+q_\text{SW}-2 q_\text{W})}{\Delta x^2 \Delta y^2} \right\} \label{eq:bipararecon}\\
q_\text{C} &= \frac{1}{16}(36  \bar q -q_\text{NW}-q_\text{SW}-4 q_\text{W})
\end{align}

\newpage

\paragraph{Hat reconstruction on both neighbouring edges} 

\begin{figure}
 \centering
 \includegraphics[width=0.4\textwidth]{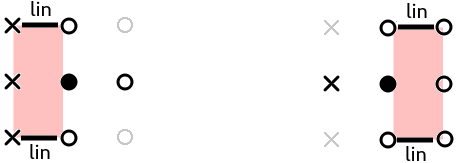} \\ \includegraphics[width=.7\textwidth]{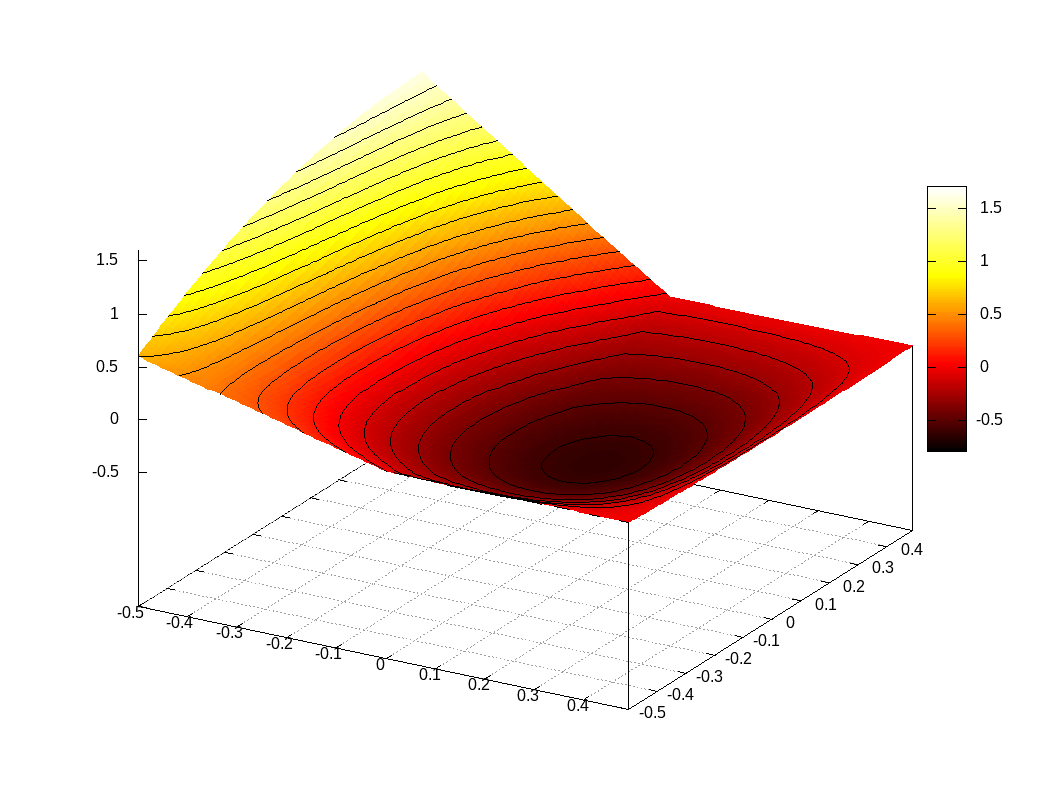}
 \caption{\emph{Top}: The case of both neighbouring edges reconstructed using hat functions, while the primary edge is reconstructed parabolically. \emph{Bottom}: The W edge is reconstructed parabolically, while the two neighbouring reconstructions are hat functions. $q_\text{NW} = 1.6$, $q_\text{W}=1.35$, $q_\text{SW} = 0.6$.}
 \label{fig:para-hat-hat}
\end{figure}
 
The interpolation problem is shown in Figure \ref{fig:para-hat-hat}.

\begin{align}
q_\text{recon}^\text{W}\Big|_{x< 0} &= \left\{a_0=q_\text{C},a_1=-\frac{q_\text{W}}{\Delta x},a_2=-\frac{2 (2 q_\text{C}-q_\text{W})}{\Delta x^2},a_3=0,a_4=
-\frac{2 (q_\text{NW}-q_\text{SW})}{\Delta x \Delta y},   \label{eq:parahathatleft}  \right . \\ & \phantom{mmm}\left.\nonumber      a_5=0,a_6=-\frac{4 q_\text{C}}{\Delta y^2},a_7=-\frac{4 (q_\text{NW}+q_\text{SW}-q_\text{W})}{\Delta x \Delta y^2},a_8=\frac{8 (2 q_\text{C}-q_\text{W})}{\Delta x^2 \Delta y^2} \right \}  \\
 q_\text{recon}^\text{W}\Big|_{x\geq 0} &= \left\{a_0=q_\text{C},a_1=-\frac{q_\text{W}}{\Delta x},a_2=-\frac{2 (2 q_\text{C}-q_\text{W})}{\Delta x^2},a_3=0,a_4=0,\label{eq:parahathatright}
 \right . \\ & \phantom{mmm}\left.\nonumber
a_5=0,a_6=-\frac{4 q_\text{C}}{\Delta y^2},a_7=\frac{4 q_\text{W}}{\Delta x \Delta y^2},a_8=\frac{8 (2 q_\text{C}-q_\text{W})}{\Delta x^2 \Delta y^2}\right\}  \\
q_\text{C} &= \frac1{32} (72  \bar q -3 (q_\text{NW}+q_\text{SW})-8 q_\text{W})
\end{align}

\newpage

\paragraph{Hat reconstruction on just one neighbouring edge} 

\begin{figure}
 \centering
 \includegraphics[width=0.4\textwidth]{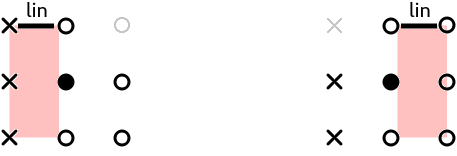} \\ \includegraphics[width=.7\textwidth]{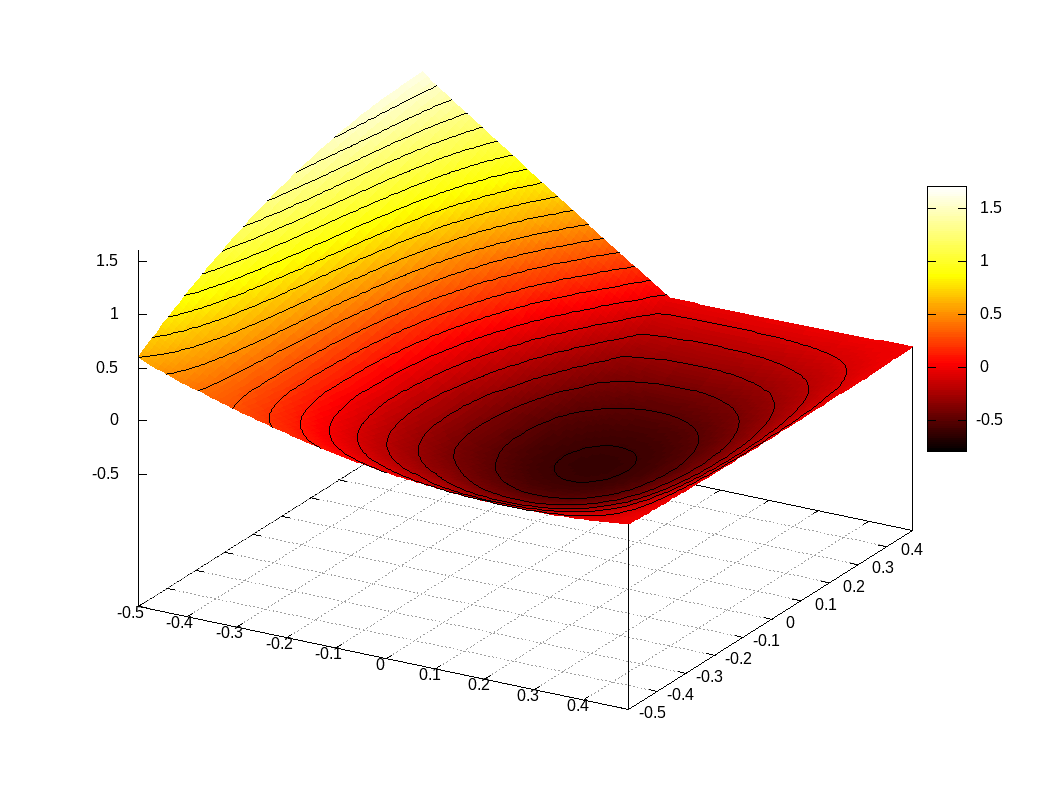}
 \caption{\emph{Top}: The case of the N edge reconstructed using hat functions, while the primary edge and the S-edge is reconstructed parabolically. \emph{Bottom}: The W and S edge is reconstructed parabolically, while the N edge is reconstructed using a hat function. $q_\text{NW} = 1.6$, $q_\text{W}=1.35$, $q_\text{SW} = 0.6$.}
 \label{fig:para-hat-para}
\end{figure}

If the N edge is reconstructed using a hat function, and both the W-edge and the S-edge parabolically, then one reconstructs the cell as follows (Figure \ref{fig:para-hat-para}):
\begin{align}
q_\text{recon}^\text{W} \Big |_{x < 0} &= \left\{a_0=q_\text{C},a_1=-\frac{q_\text{W}}{\Delta x},a_2=-\frac{2 (2 q_\text{C}-q_\text{W})}{\Delta x^2},a_3=0,a_4=-\frac{2 q_\text{NW}-q_\text{SW}}{\Delta x \Delta y},\label{eq:parahatparaleft}
     \right . \\ & \phantom{mmm}\left.\nonumber      
a_5=-\frac{2 q_\text{SW}}{\Delta x^2 \Delta y},a_6=-\frac{4 q_\text{C}}{\Delta y^2},a_7=-\frac{2 (2 q_\text{NW}+q_\text{SW}-2 q_\text{W})}{\Delta x \Delta y^2},
\right . \\ & \phantom{mmm}\left.\nonumber  
a_8=\frac{4 (4 q_\text{C}+q_\text{SW}-2 q_\text{W})}{\Delta x^2 \Delta y^2}\right\} \\
q_\text{recon}^\text{W}\Big |_{x \geq 0} &= \left\{a_0=q_\text{C},a_1=-\frac{q_\text{W}}{\Delta x},a_2=-\frac{2 (2 q_\text{C}-q_\text{W})}{\Delta x^2},a_3=0,a_4=\frac{q_\text{SW}}{\Delta x \Delta y},  \label{eq:parahatpararight}   \right . \\ & \phantom{mmm}\left.\nonumber      a_5=-\frac{2 q_\text{SW}}{\Delta x^2 \Delta y},a_6=-\frac{4 q_\text{C}}{\Delta y^2},a_7=-\frac{2 (q_\text{SW}-2 q_\text{W})}{\Delta x \Delta y^2},
\right . \\ & \phantom{mmm}\left.\nonumber  
a_8=\frac{4 (4 q_\text{C}+q_\text{SW}-2 q_\text{W})}{\Delta x^2 \Delta y^2}\right\} \\
q_\text{C} &= \frac{1}{32} (72  \bar q -3 q_\text{NW}-2 (q_\text{SW}+4 q_\text{W}))
\end{align}

\begin{figure}
 \centering
 \includegraphics[width=0.4\textwidth]{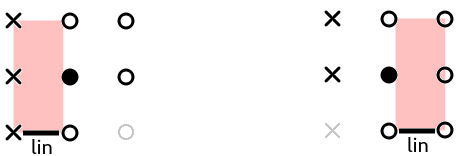} \hfill \includegraphics[width=.7\textwidth]{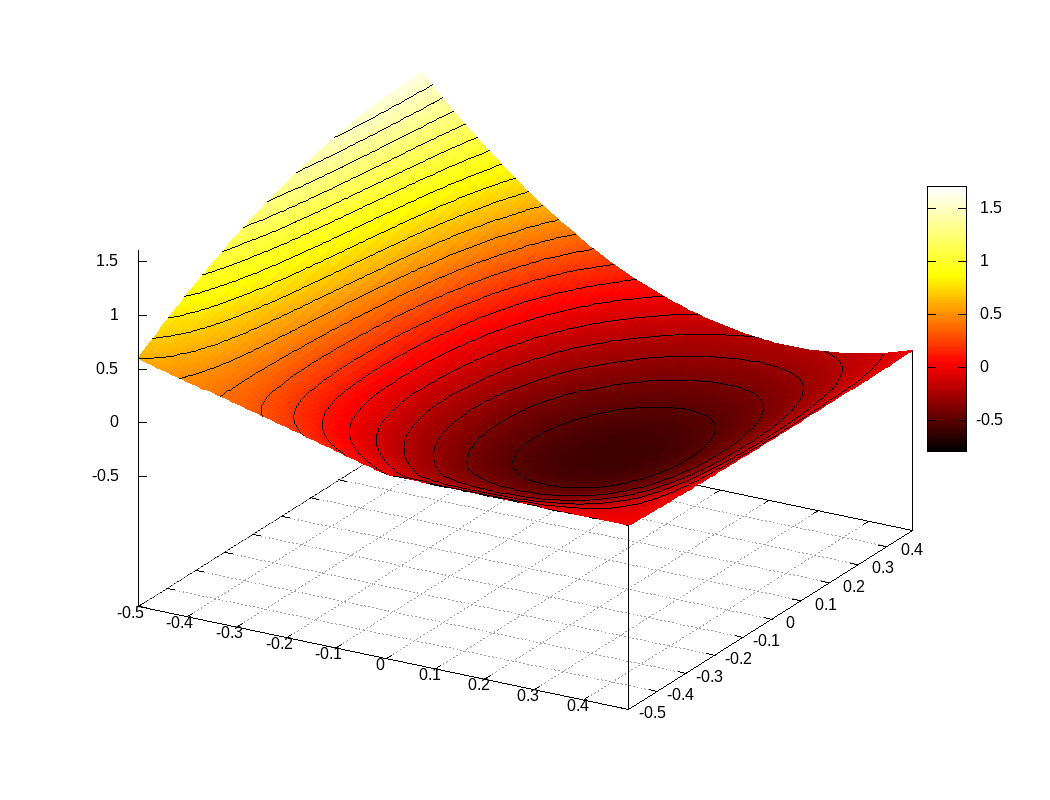}
 \caption{\emph{Top}: The case of the S edge reconstructed using hat functions, while the primary edge is reconstructed parabolically. \emph{Bottom}: The W and N edge is reconstructed parabolically, while the S edge is reconstructed using a hat function. $q_\text{NW} = 1.6$, $q_\text{W}=1.35$, $q_\text{SW} = 0.6$.}
 \label{fig:para-para-hat}
\end{figure}

If it is the S edge, then (Figure \ref{fig:para-para-hat}):
\begin{align}
q_\text{recon}^\text{W}\Big |_{x < 0} &= \left \{a_0=q_\text{C},a_1=-\frac{q_\text{W}}{\Delta x},a_2=-\frac{2 (2 q_\text{C}-q_\text{W})}{\Delta x^2},a_3=0,a_4=-\frac{q_\text{NW}-2 q_\text{SW}}{\Delta x \Delta y},  \label{eq:paraparahatleft}   \right . \\ & \phantom{mmm}\left.\nonumber      a_5=\frac{2 q_\text{NW}}{\Delta x^2 \Delta y},a_6=-\frac{4 q_\text{C}}{\Delta y^2},a_7=-\frac{2 (q_\text{NW}+2 q_\text{SW}-2 q_\text{W})}{\Delta x \Delta y^2},
\right . \\ & \phantom{mmm}\left.\nonumber  
a_8=\frac{4 (4 q_\text{C}+q_\text{NW}-2 q_\text{W})}{\Delta x^2 \Delta y^2}\right \} \\
q_\text{recon}^\text{W}\Big |_{x \geq 0} &= \left\{a_0=q_\text{C},a_1=-\frac{q_\text{W}}{\Delta x},a_2=-\frac{2 (2 q_\text{C}-q_\text{W})}{\Delta x^2},a_3=0,a_4=-\frac{q_\text{NW}}{\Delta x \Delta y},   \label{eq:paraparahatright}   \right . \\ & \phantom{mmm}\left.\nonumber      a_5=\frac{2 q_\text{NW}}{\Delta x^2 \Delta y},a_6=-\frac{4 q_\text{C}}{\Delta y^2},a_7=-\frac{2 (q_\text{NW}-2 q_\text{W})}{\Delta x \Delta y^2},
\right . \\ & \phantom{mmm}\left.\nonumber  
a_8=\frac{4 (4 q_\text{C}+q_\text{NW}-2 q_\text{W})}{\Delta x^2 \Delta y^2}\right\}\\
q_\text{C}&= \frac1{32} (72  \bar q -2 q_\text{NW}-3 q_\text{SW}-8 q_\text{W})
\end{align}

\paragraph{Proof of continuity} \label{ssec:continuitypara}

It is obvious from the sketches of the interpolation problem in Figures \ref{fig:example-bipara}--\ref{fig:para-para-hat} that the reconstructions interpolate the values on the cell interfaces. What remains to be shown is that the piecewise defined reconstruction is continuous:

\begin{theorem} \label{thm:continuitypara}
 The reconstructions from Section \ref{ssec:para} are continuous along the line $x=0$ where the two pieces are joined.
\end{theorem}
\begin{proof}
 As is obvious from the sketches of the interpolation problems in Figures \ref{fig:para-hat-hat}--\ref{fig:para-para-hat}, the three points along $x=0$, i.e.
 \begin{align}
  q_\text{recon}\left(0, \frac{\Delta y}{2} \right) &= 0 &
  q_\text{recon}\left(0, 0 \right) &= q_\text{C} &
  q_\text{recon}\left(0, -\frac{\Delta y}{2} \right) &= 0 
 \end{align}
 are part of the interpolation. Recall that the restriction of a biparabolic function onto the straight line $x=0$ is a parabola in $y$, and that the latter is uniquely defined by three points. Therefore, all the values of the reconstruction along $x=0$ agree for all the reconstructions presented in Section \ref{ssec:para}.
\end{proof}

\subsubsection{Hat reconstruction on W edge}\label{ssec:hat}

If the W-edge is reconstructed as a hat function, then necessarily one needs to consider a piecewise defined reconstruction with the pieces joined along $y = 0$. The reconstruction in each piece only depends on whether the other adjacent edge is reconstructed parabolically or as a hat function. One thus has less cases to consider.

Consider the top piece, i.e. the one defined on $[-\frac{\Delta x}{2}, \frac{\Delta x}{2}] \times [0, \frac{\Delta y}{2}]$. It is bordered by the N-edge. If the N-edge is reconstructed as a hat function then one needs additionally to define the reconstruction piecewise in the left and right halves (joined along $x=0$), i.e. the reconstruction is piecewise by quadrant. This is not necessary if the N-edge is reconstructed parabolically.

\paragraph{Parabolic reconstruction on at least one neighbouring edge} 

Here, the situation is considered in which either the N-edge or the S-edge are reconstructed as parabolae. Then it is possible to provide a biparabolic reconstruction of, respectively, the top or bottom half of the cell. 

These cases can occur individually or simultaneously. If both the N-edge and the S-edge are reconstructed parabolically, then the entire reconstruction of the cell is given by the two pieces given in \eqref{eq:topparthatpara}--\eqref{eq:bottomparthatpara}. If, for example, the N-edge is reconstructed parabolically, and the S-edge as a hat function, then the top piece of the reconstruction in the cell is to be taken from \eqref{eq:topparthatpara}, while the bottom piece used should be the one from \eqref{eq:hathatbottomleft}--\eqref{eq:hathatbottomright}.

\begin{figure}
 \centering
 \includegraphics[width=0.4\textwidth]{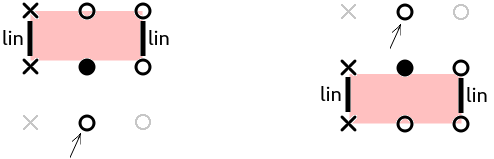} \\ \includegraphics[width=.7\textwidth]{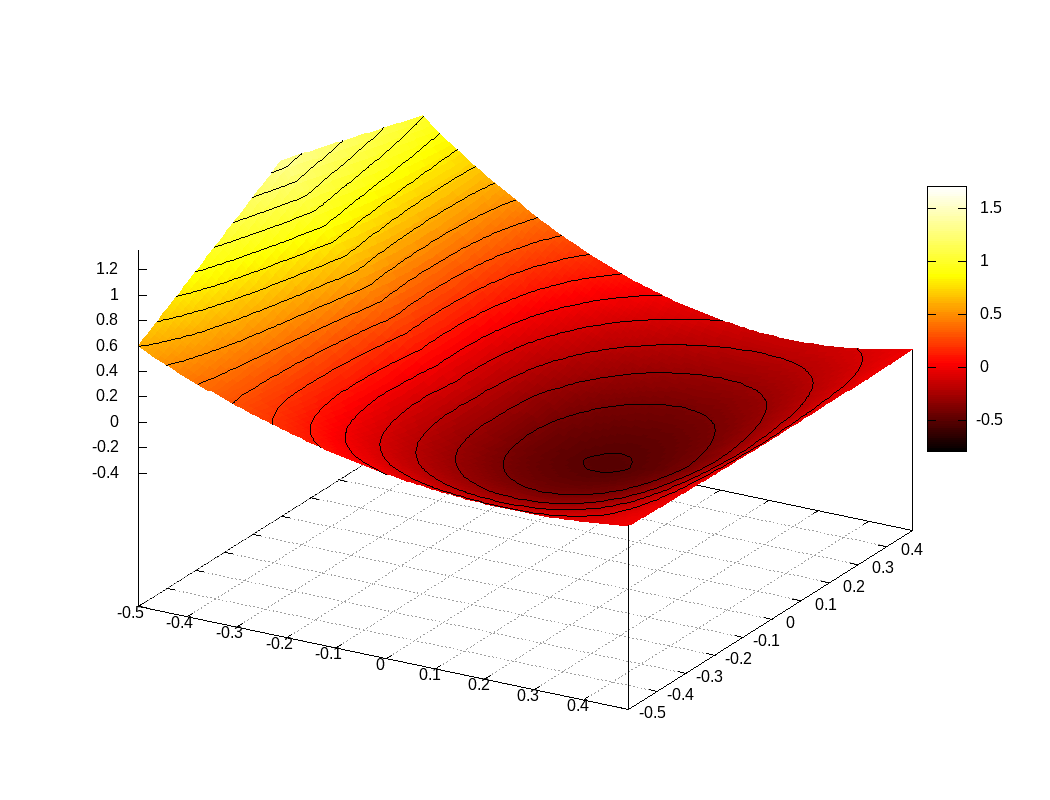}
 \caption{\emph{Top}: The case of the neighbouring edges reconstructed using parabolas, while the primary edge is reconstructed using the hat function. \emph{Bottom}: The W edge is reconstructed as a hat function, the other edges are reconstructed parabolically. $q_\text{NW} = 1$, $q_\text{W}=1.5$, $q_\text{SW} = 0$.}
 \label{fig:hat-para}
\end{figure}

See Figure \ref{fig:hat-para} for the setup of the interpolation problem.

\begin{align}
q_\text{recon}^{W} \Big|_{y\geq 0} &= \left\{a_0=q_\text{C},a_1=-\frac{q_\text{W}}{\Delta x},a_2=-\frac{2 (2 q_\text{C}-q_\text{W})}{\Delta x^2},a_3=0,a_4=-\frac{2 (q_\text{NW}-q_\text{W})}{\Delta x \Delta y},     \right .\label{eq:topparthatpara}  \\ & \phantom{mmm}\left.\nonumber      a_5=\frac{4 (q_\text{NW}-q_\text{W})}{\Delta x^2 \Delta y},a_6=-\frac{4 q_\text{C}}{\Delta y^2},a_7=0,a_8=\frac{16 q_\text{C}}{\Delta x^2 \Delta y^2}\right\} \label{eq:hatWparaN}\\
\frac1{\Delta x \Delta y} & \int_{y\geq 0}  q_\text{recon} \,\dd x \dd y = \frac{2 q_\text{C}}9+\frac{q_\text{NW}+q_\text{W}}{24}\\
q_\text{recon}^{W} \Big|_{y < 0} &= \left\{a_0=q_\text{C},a_1=-\frac{q_\text{W}}{\Delta x},a_2=-\frac{2 (2 q_\text{C}-q_\text{W})}{\Delta x^2},a_3=0,a_4=\frac{2 (q_\text{SW}-q_\text{W})}{\Delta x \Delta y},     \right . \label{eq:bottomparthatpara}\\ & \phantom{mmm}\left.\nonumber      a_5=-\frac{4 (q_\text{SW}-q_\text{W})}{\Delta x^2 \Delta y},a_6=-\frac{4 q_\text{C}}{\Delta y^2},a_7=0,a_8=\frac{16 q_\text{C}}{\Delta x^2 \Delta y^2}\right\} \label{eq:hatWparaS}\\
\frac1{\Delta x \Delta y}& \int_{y < 0} q_\text{recon} \,\dd x \dd y = \frac{2 q_\text{C}}9+ \frac{q_\text{SW}+q_\text{W}}{24}
\end{align}

\newpage

\paragraph{Hat reconstruction on at least one neighbouring edge} \label{ssec:hathat}

\begin{figure}
 \centering
 \includegraphics[width=0.4\textwidth]{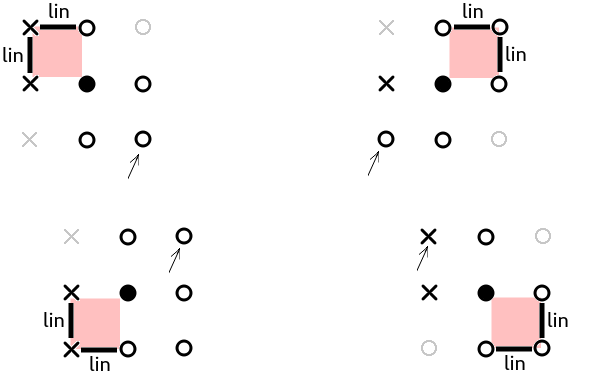} \\ \includegraphics[width=.45\textwidth]{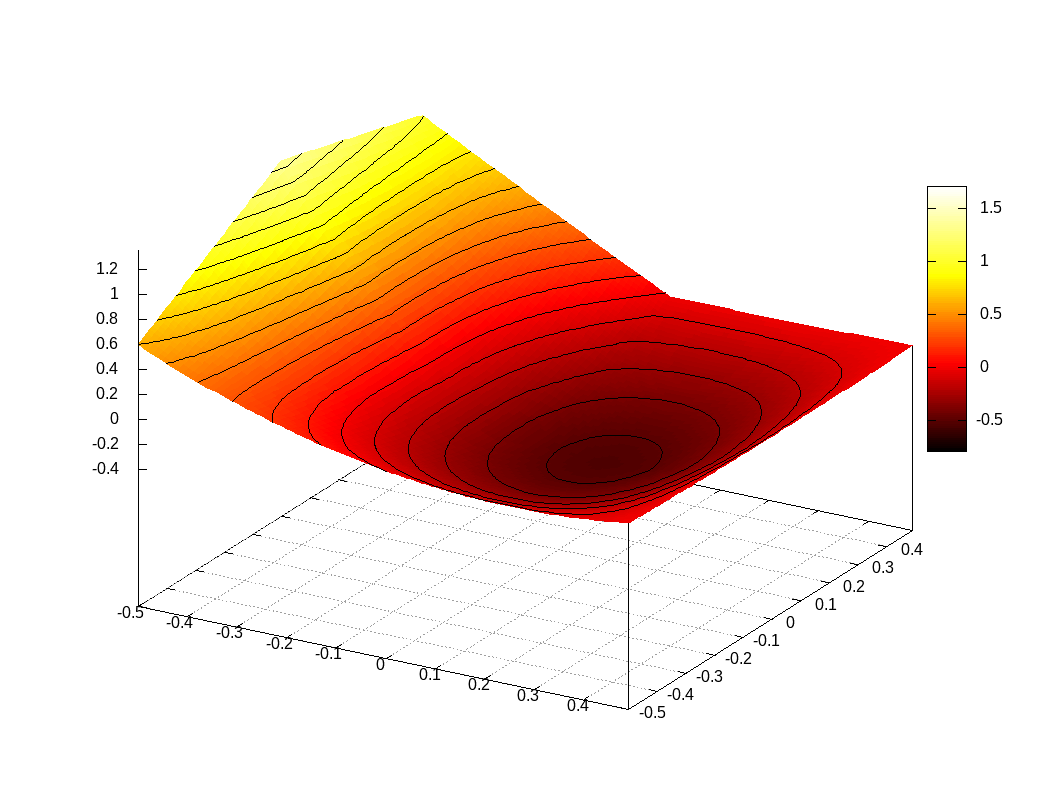} \hfill \includegraphics[width=.45\textwidth]{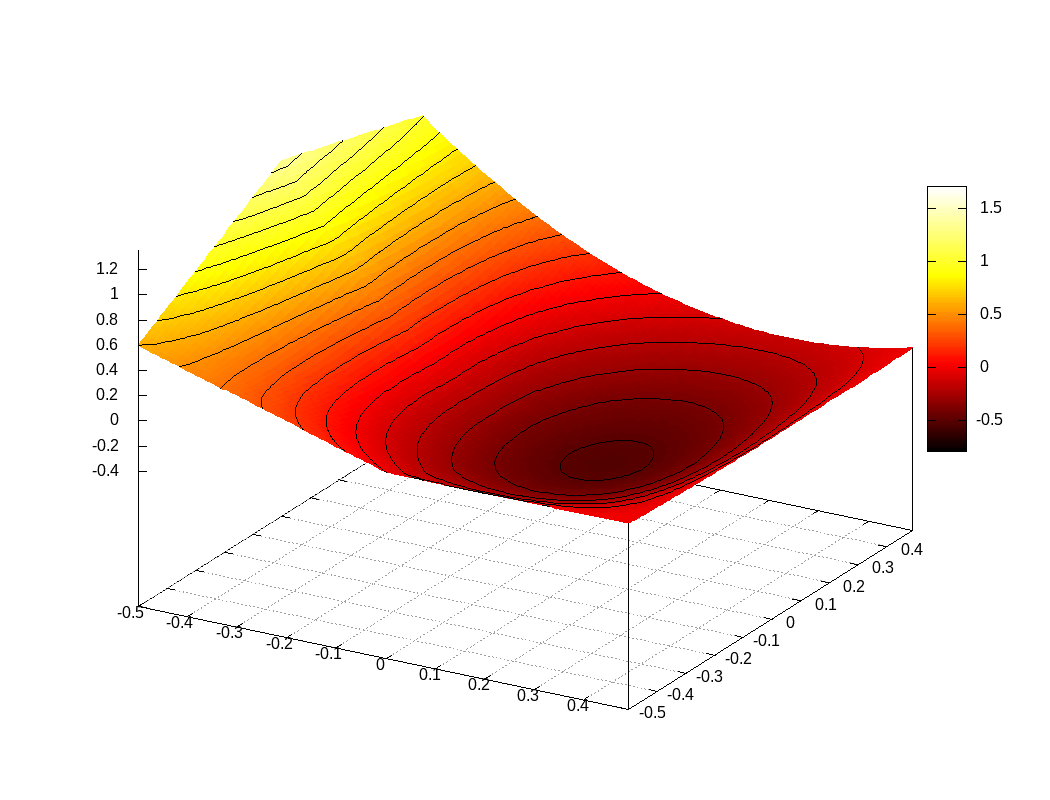} \\
 \includegraphics[width=.7\textwidth]{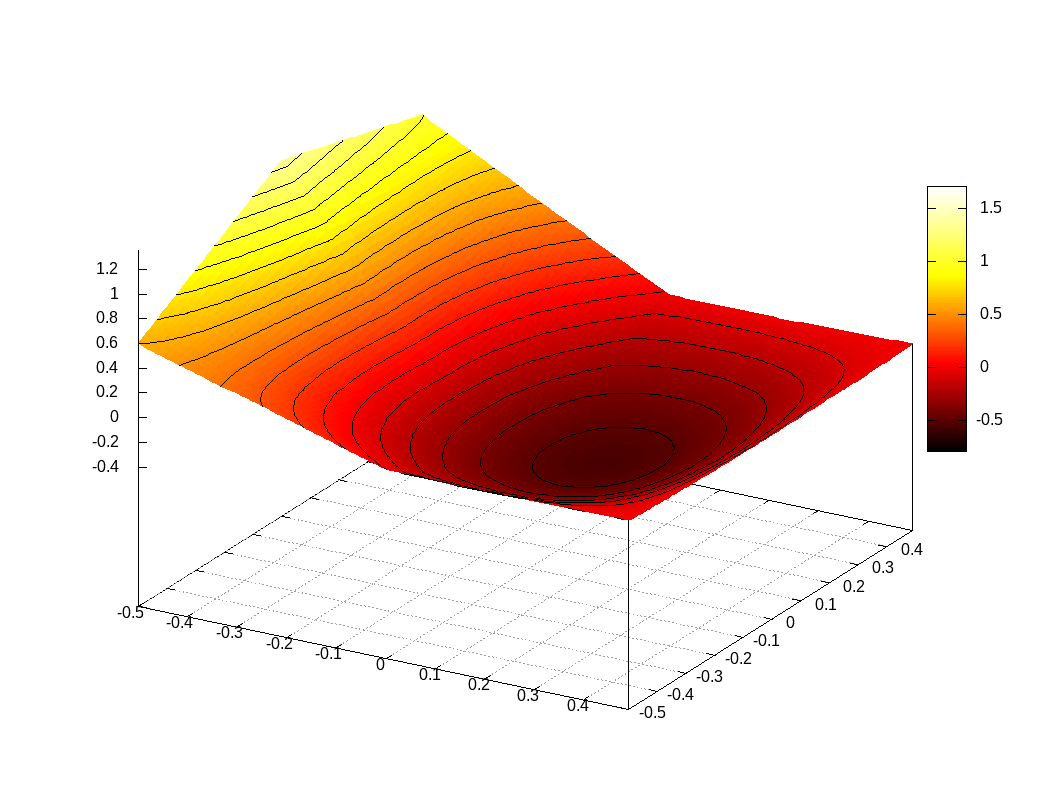}
 \caption{\emph{Top}: The case of (possibly) all three edges reconstructed using hat functions. The reconstruction inside the cell is defined on four quadrants. \emph{Middle}: The W edge is reconstructed as a hat function, and also N (\emph{left}) / S (\emph{right}). \emph{Bottom}: All edges are reconstructed as hat functions.}
 \label{fig:hat-hat}
\end{figure}

In this case the reconstruction is defined piecewise on each quadrant. The biparabolic reconstructions are obtained from interpolation problems shown in Figure \ref{fig:hat-hat}.

If the N edge is reconstructed as a hat function, then the top half $[-\frac{\Delta x}{2}, \frac{\Delta x}{2}] \times [0, \frac{\Delta y}{2}]$ of the cell is to be reconstructed as
\begin{align}
q_\text{recon}^\text{W} \Big |_{y \geq 0, x < 0}&=\left\{a_0=q_\text{C},a_1=-\frac{q_\text{W}}{\Delta x},a_2=-\frac{2 (2 q_\text{C}-q_\text{W})}{\Delta x^2},a_3=0,a_4=-\frac{3 q_\text{NW}-2 q_\text{W}}{\Delta x \Delta y},  \label{eq:hatWhatNleft}   \right . \\ & \phantom{mmm}\left.\nonumber      a_5=\frac{2 (q_\text{NW}-2 q_\text{W})}{\Delta x^2 \Delta y},a_6=-\frac{4 q_\text{C}}{\Delta y^2},a_7=-\frac{2 q_\text{NW}}{\Delta x \Delta y^2},
\right . \\ & \phantom{mmm}\left.\nonumber  
a_8=\frac{4 (4 q_\text{C}-q_\text{NW})}{\Delta x^2 \Delta y^2}\right\} \\
q_\text{recon}^\text{W} \Big |_{y \geq 0, x \geq 0}&= \left\{a_0=q_\text{C},a_1=-\frac{q_\text{W}}{\Delta x},a_2=-\frac{2 (2 q_\text{C}-q_\text{W})}{\Delta x^2},a_3=0,a_4=\frac{q_\text{SW}}{\Delta x \Delta y}, \label{eq:hatWhatNright}    \right . \\ & \phantom{mmm}\left.\nonumber      a_5=-\frac{2 q_\text{SW}}{\Delta x^2 \Delta y},a_6=-\frac{4 q_\text{C}}{\Delta y^2},a_7=-\frac{2 (q_\text{SW}-2 q_\text{W})}{\Delta x \Delta y^2},
\right . \\ & \phantom{mmm}\left.\nonumber  
a_8=\frac{4 (4 q_\text{C}+q_\text{SW}-2 q_\text{W})}{\Delta x^2 \Delta y^2}\right \} \\
\frac1{\Delta x \Delta y}&\int_{y \geq 0} q_\text{recon} \,\dd x \dd y = \frac{2 q_\text{C}}9+\frac1{576} (35 q_\text{NW}+q_\text{SW}+22 q_\text{W}) 
\end{align}

If the S edge is reconstructed as a hat function, then the reconstruction reads

\begin{align}
q_\text{recon}^\text{W} \Big |_{y < 0, x < 0}&= \left\{ a_0=q_\text{C},a_1=-\frac{q_\text{W}}{\Delta x},a_2=-\frac{2 (2 q_\text{C}-q_\text{W})}{\Delta x^2},a_3=0,  \label{eq:hathatbottomleft}   \right . \\ & \phantom{mmm}\left.\nonumber    a_4=-\frac{-3 q_\text{SW}+2 q_\text{W}}{\Delta x \Delta y},  a_5=-\frac{2 (q_\text{SW}-2 q_\text{W})}{\Delta x^2 \Delta y},a_6=-\frac{4 q_\text{C}}{\Delta y^2},
\right . \\ & \phantom{mmm}\left.\nonumber  
a_7=-\frac{2 q_\text{SW}}{\Delta x \Delta y^2},a_8=\frac{4 (4 q_\text{C}-q_\text{SW})}{\Delta x^2 \Delta y^2}\right \} \\
q_\text{recon}^\text{W}\Big |_{y < 0, x \geq 0} &= \left\{a_0=q_\text{C},a_1=-\frac{q_\text{W}}{\Delta x},a_2=-\frac{2 (2 q_\text{C}-q_\text{W})}{\Delta x^2},a_3=0,a_4=-\frac{q_\text{NW}}{\Delta x \Delta y},    \label{eq:hathatbottomright} \right . \\ & \phantom{mmm}\left.\nonumber      a_5=\frac{2 q_\text{NW}}{\Delta x^2 \Delta y},a_6=-\frac{4 q_\text{C}}{\Delta y^2},a_7=-\frac{2 (q_\text{NW}-2 q_\text{W})}{\Delta x \Delta y^2},
\right . \\ & \phantom{mmm}\left.\nonumber  
a_8=\frac{4 (4 q_\text{C}+q_\text{NW}-2 q_\text{W})}{\Delta x^2 \Delta y^2}\right\} \\
\frac1{\Delta x \Delta y}& \int_{y <0} q_\text{recon} \,\dd x \dd y =  \frac{2 q_\text{C}}9+ \frac1{576} (q_\text{NW}+35 q_\text{SW}+22 q_\text{W})
\end{align}

\paragraph{Proof of continuity} \label{ssec:continuityhat}

\begin{theorem} \label{thm:continuityhat}
 The reconstructions in Section \ref{ssec:hat} are continuous along $x = 0$ and along $y = 0$.
\end{theorem}
\begin{proof}
 In complete analogy to the proof of Theorem \ref{thm:continuitypara} one observes from the sketches of the interpolation problem in Figures \ref{fig:hat-para}--\ref{fig:hat-hat} that the points along $x = 0$ and $y = 0$ are always included. The three points along $x = 0$ and the three points along $y = 0$ each define a unique parabola.
\end{proof}

\subsection{Plateau-limiting} \label{sec:plateau}

\begin{figure}
 \centering
\includegraphics[width=0.5\textwidth]{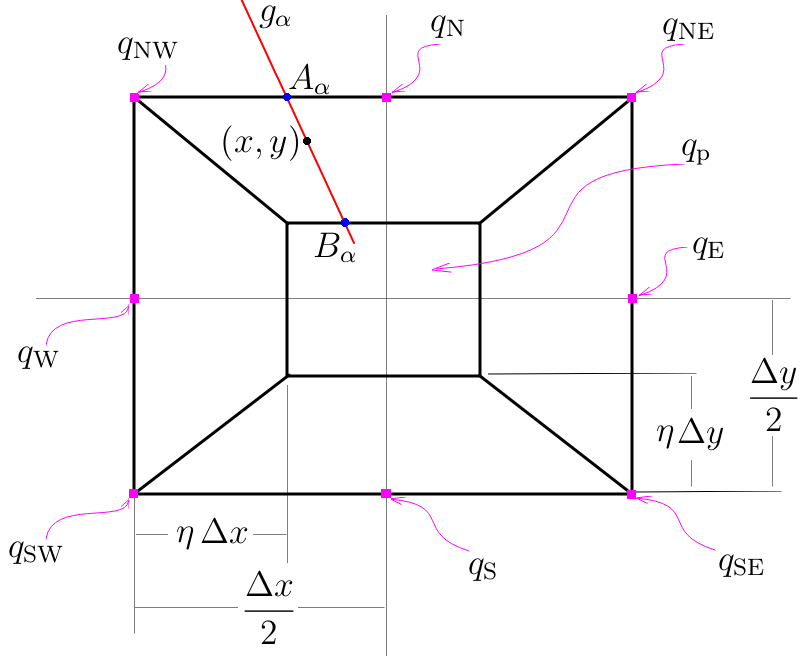} 
\caption{Sketch of the interpolation between the plateau and the boundary of the cell.}
\label{fig:sketch}
\end{figure}

Consider a situation in which \eqref{eq:quasimonotone} is true, while the reconstruction described above exceeds $m$ or $M$. In that case, the idea of a plateau reconstruction (see Figure \ref{fig:sketch}) is to introduce a rectangle a distance $\eta \Delta x$ (or $\eta \Delta y$) away from the cell boundary, i.e. $$\left[\Delta x \left(-\frac12 + \eta \right), \Delta x \left(\frac12 - \eta \right)\right] \times \left[\Delta y \left(-\frac12 + \eta \right), \Delta y \left(\frac12 - \eta \right)\right]$$ with $\eta \in (0,\frac12)$ where the value of the reconstruction will be constant and equal to $q_\text{p}$, a value to be determined to ensure that the average of the reconstruction equals the given average (see Figure \ref{fig:plateau} for an example). This rectangle will be referred to as \textbf{plateau}. The remaining four trapezes will be the supports of functions that continuously join the reconstruction along the edge to the plateau in the simplest possible way. Because reconstructions along edges are either parabolas or hats, every trapezoidal region is either joining the plateau to a parabola or to a hat function. $\eta$ will be chosen in such a way that the maximum principle is guaranteed. It is clear that, as \eqref{eq:quasimonotone} is true, this can always be done by choosing $\eta$ small enough.

\subsubsection{Interpolation in the trapezes} \label{ssec:plateauinterpolation}

Consider for definiteness the northern trapeze. Define a point $A_\alpha := \left(-\frac{\Delta x}{2} + \alpha \Delta x, \frac{\Delta y}{2}\right) \in \mathbb R^2$ parametrized by $\alpha \in [0,1]$. Define a point $$B_\alpha := \left(\Delta x(-\frac12 + \eta) + \alpha \Delta x(1 -2 \eta) , \Delta y(\frac12 - \eta)\right) $$ on the northern edge of the plateau. Observe that as $\alpha$ goes from 0 to 1, both points move all the way from the left to the right on their respective edges. The straight line
\begin{align}
 g_\alpha := \left\{ (x, y) : \frac{x}{\Delta x} = -\frac{1}{2} + \alpha - \left(\frac{y}{\Delta y} - \frac{1}{2}\right)  ( 1  -2 \alpha)  \right\}
\end{align}
connects them. Obviously, given $x$ and $y$ there is a unique
\begin{align}
 \alpha = \frac{\frac{x}{\Delta x} + \frac{y}{\Delta y} }{2\frac{y}{\Delta y}} = \frac{x\Delta y + y\Delta x }{2y\Delta x} \label{eq:alphaasfctofxan}
\end{align}
The idea of the reconstruction is to associate to a point $(x, y)$ the value given by a linear interpolation between the value of the reconstruction at $A_\alpha$ and the (constant) value $q_\text{p}$ at $B_\alpha$. In particular this means that the diagonal edges of the reconstruction (connections between the corners of the cell and the corners of the plateau) are reconstructed as straight lines.

The four trapezes can be reconstructed individually, because continuity along the diagonal segments where they join is already guaranteed by the above procedure. For a given trapeze, the choice of reconstruction thus merely depends on whether the adjacent edge is reconstructed parabolically (see Section \ref{ssec:plateaupara}) or as a hat function (see Section \ref{ssec:plateauhat}).

\subsubsection{Parabolic reconstruction along the edge} \label{ssec:plateaupara}

The parabolic reconstruction along the N-edge is given by
\begin{align}
 q_\text{parabolic}^\text{N}(x) = q_\text{N} + \frac{x}{\Delta x} (q_\text{NE} - q_\text{NW}) + 2\frac{x^2}{\Delta x^2}(q_\text{NE} + q_\text{NW} - 2 q_\text{N})  \qquad x \in \left[-\frac{\Delta x}{2}, \frac{\Delta x}{2}\right]
\end{align}
The value of this parabolic reconstruction is sought at the location $\xi$ of point $A_\alpha$ with $\alpha$ given by \eqref{eq:alphaasfctofxan}:
\begin{align}
 \xi = \Delta x \left( -\frac{1}{2} + \frac{\frac{x}{\Delta x}\Delta y + y }{2y} \right) = \Delta x\frac{\frac{x}{\Delta x}}{2\frac{y}{\Delta y}}
\end{align}
Finally, the reconstruction at $(x, y)$ is assigned the value
\begin{align}
 q_\text{recon}^\text{N}(x, y) &:= q_\text{parabolic}^\text{N}(\xi) + \left(y - \frac{\Delta y}{2}\right) \frac{q_\text{p} - q_\text{parabolic}^\text{N}(\xi)}{- \Delta y \eta}  \\
 &= q_\text{parabolic}^\text{N}(\xi)\left(1 +  \frac{\frac{y}{\Delta y} - \frac{1}{2}}{\eta}  \right) - \frac{\frac{y}{\Delta y} - \frac{1}{2}}{\eta}  q_\text{p} 
\end{align}
with  
\begin{align}
 q_\text{parabolic}^\text{N}(\xi) &=  q_\text{N} + \frac{\hat x}{2\hat y} (q_\text{NE} - q_\text{NW}) + 2\left( \frac{\hat x}{2\hat y} \right)^2 (q_\text{NE} + q_\text{NW} - 2 q_\text{N}) 
\end{align}
and $\hat x := \frac{x}{\Delta x}$ and $\hat y := \frac{y}{\Delta y}$. Observe that the reconstruction is not polynomial, but lies in
\begin{align}
 \mathrm{span}\left(1, \hat x, \hat y, \frac{\hat x}{\hat y}, \frac{{\hat x}^2}{\hat y}, \frac{{\hat x}^2}{{\hat y}^2}  \right)
\end{align}

For reference we give the four reconstructions: {\footnotesize
\begin{align}
 q_\text{recon}^\text{trapeze,W}(x, y) &= q_\text{p}\frac{1+2 \hat x}{2 \eta} + (-1+2 \eta-2 \hat x) \left( \frac{q_\text{W} }{2 \eta} - \frac{(q_\text{NW}-q_\text{SW})  y}{4 \eta \hat x}+\frac{(q_\text{NW}+q_\text{SW}-2 q_\text{W}) \hat  y^2}{4 \eta \hat x^2} \right) \label{eq:trapezeWparabolic}\\
 q_\text{recon}^\text{trapeze,E}(x, y) &= q_\text{p}\frac{1-2 \hat  x}{2 \eta}+ (-1+2 \eta+2 \hat x) \left( \frac{q_\text{E} }{2 \eta} +\frac{(q_\text{NE}-q_\text{SE})  \hat y}{4 \eta \hat x}-\frac{(2 q_\text{E}-q_\text{NE}-q_\text{SE})  \hat y^2}{4 \eta \hat x^2} \right ) \label{eq:trapezeEparabolic}\\
 q_\text{recon}^\text{trapeze,N}(x, y) &= q_\text{p}\frac{1-2  \hat y}{2 \eta}+   (-1+2 \eta+2 \hat y) \left( -\frac{(2 q_\text{N}-q_\text{NE}-q_\text{NW}) \hat x^2 }{4 \eta \hat y^2}+\frac{(q_\text{NE}-q_\text{NW}) \hat x }{4 \eta \hat y}+\frac{q_\text{N} }{2 \eta}\right) \label{eq:trapezeNparabolic}\\
 q_\text{recon}^\text{trapeze,S}(x, y) &= q_\text{p}\frac{1+2\hat y}{2\eta} + (1-2 \eta+2 \hat y) \left( \frac{(2 q_\text{S}-q_\text{SE}-q_\text{SW}) \hat x^2 }{4 \eta\hat  y^2}+\frac{(q_\text{SE}-q_\text{SW}) \hat x }{4 \eta \hat y}-\frac{q_\text{S} }{2\eta} \right) \label{eq:trapezeSparabolic}
\end{align}}

The integrals over the four regions are

\begin{align}
\frac{1}{\Delta x \Delta y} \int_{\text{trapeze,W}} q_\text{recon} \,\dd x \dd y &=  \frac1{36} \eta \Big(6 (3-4 \eta) q_\text{P}-(2 \eta-3) (4 q_\text{W} + q_\text{NW}+q_\text{SW})\Big)\\
 \frac{1}{\Delta x \Delta y} \int_{\text{trapeze,E}} q_\text{recon} \,\dd x \dd y &= \frac1{36}  \eta \Big( 6( 3   -4\eta)  q_\text{P}  -(2 \eta-3)( 4q_\text{E} +q_\text{NE}+q_\text{SE})    \Big)\\
 \frac{1}{\Delta x \Delta y} \int_{\text{trapeze,N}} q_\text{recon} \,\dd x \dd y &= \frac1{36} \eta \Big(6 (3-4 \eta) q_\text{P} -(2 \eta-3) (4 q_\text{N}+q_\text{NE}+q_\text{NW}) \Big)\\
 \frac{1}{\Delta x \Delta y} \int_{\text{trapeze,S}} q_\text{recon} \,\dd x \dd y &= \frac1{36}\eta \Big (6 (3-4 \eta) q_\text{P}-(2 \eta-3) (4 q_\text{S}+q_\text{SE}+q_\text{SW})\Big)
\end{align}
and the integral over the plateau obviously
\begin{align}
 \frac{1}{\Delta x \Delta y} \int_{\text{plateau}} q_\text{recon} \,\dd x \dd y = (1-2\eta)^2 q_\text{p}
\end{align}

\subsubsection{Hat-function reconstruction along the edge} \label{ssec:plateauhat}

If an edge is reconstructed using a hat-function, then the reconstruction of the trapeze follows the algorithm outlined at the beginning of Section \ref{sec:plateau}, but is naturally defined in a piecewise fashion. The reconstruction of the W-trapeze is

\begin{align}
q_\text{recon}^{\text{trapeze,W}}(x,y) \Big |_{y \geq 0} &=
q_\text{W}-\frac{\Delta x (q_\text{NW}-q_\text{W}) y}{x\Delta y}+\frac{\left(\frac{\Delta x}{2}+x\right) \left(q_\text{P}-q_\text{W}+\frac{\Delta x (q_\text{NW}-q_\text{W}) y }{x\Delta y} \right)}{\Delta x \eta} \label{eq:trapezeWhattop} \\
q_\text{recon}^{\text{trapeze,W}}(x,y) \Big |_{y < 0} &=
q_\text{W}+\frac{\Delta x (q_\text{SW}-q_\text{W}) y}{x\Delta y}+\frac{\left(\frac{\Delta x}{2}+x\right) (q_\text{P}-q_\text{W}-\frac{\Delta x (q_\text{SW}-q_\text{W}) y}{x\Delta y}}{\Delta x \eta} \label{eq:trapezeWhatbottom}
\end{align}

\begin{align}
\frac{1}{\Delta x \Delta y} \int_{\text{trapeze,W}} q_\text{recon} \,\dd x \dd y &= 
\frac{1}{6} (3-4 \eta) \eta q_\text{p} +\frac{1}{24} \eta (2 \eta-3) (q_\text{NW} +q_\text{SW}+2 q_\text{W}) 
\end{align}

The reconstructions of the other trapezes can be obtained by rotation as in Equations \eqref{eq:rotationreconS}--\eqref{eq:rotationreconE}.

\subsubsection{Choice of the plateau value and the maximum principle} \label{ssec:plateaumaximumprinciple}

\begin{theorem} \label{thm:reconplateauinterpol}
 There exists a choice of $\eta$ such that the reconstruction is conservative and $m \leq q_\text{recon}(x,y) \leq M$ for all $x,y$ inside the cell.
\end{theorem}
\begin{proof}
For any choice of $q_\text{p} \in (m, M)$, the reconstruction inside the cell fulfills $m \leq q_\text{recon} \leq M$, because the reconstructions inside the trapezes are interpolations along straight lines between $q_\text{p}$ and a maximum-preserving reconstruction along the edge. For the same reason, as $\eta \to 0$, the average of the reconstruction over the cell approaches $q_\text{p}$, because the reconstructions inside the trapezes remain bounded and their contribution to the cell average thus vanishes in the limit. Thus, for all $\epsilon > 0$ sufficiently small one can find an $\eta > 0$ such that $\frac{1}{\Delta x \Delta y} \int_c q_\text{recon}(x, y) \, \dd x \dd y = q_\text{p} + a$ with $|a| < \epsilon$. Then, choosing $q_\text{p} := \bar q - a$ ensures conservativity of the reconstruction. At the same time, as $m < \bar q < M$, one simply needs to choose $\epsilon < \min\left(M-\bar q, \bar q - m\right)$ to ensure that $m < q_\text{p} < M$. 
\end{proof}

For example, if all edges are reconstructed parabolically, then the average of the reconstruction over the entire cell is
\begin{align}
 q_\text{p}- \frac{1}{9} \eta (2 \eta-3) \Big (4 E -6 q_\text{p} +2 V \Big) \overset{!}{=} \bar q
\end{align}
(where $4V := q_\text{NE}+q_\text{NW} +q_\text{SE}+q_\text{SW}$, $4 E := q_\text{E}+ q_\text{N} +q_\text{S} +q_\text{W}$)
which gives the value of $q_\text{p}$:
\begin{align}
 q_\text{p}  = \frac{9 \bar q +\eta (2 \eta-3) (4E+2V) }{ 3 (3-6 \eta+4 \eta^2)}
\end{align}
The polynomial in the denominator does not have real zeros. 

What thus remains is the choice of $\eta$. 
The only bounds on $\eta$ originate from the condition 
\begin{align}
 m < q_\text{p} < M
\end{align}
The equation $q_\text{p} = \mu \in \{ m,M\}$ is quadratic in $\eta$ -- and this is true in general and not just in this example. It is therefore easy to identify real, positive solutions and to take their minimum. 
In practice, having established a minimum, $\eta$ is chosen to be half of it. In case no real, positive solutions are identified, $\eta$ is not subject to any conditions and we choose $\eta = \frac14$.


\begin{thebibliography}{MRKG03}

\bibitem[AB23a]{abgrall22}
Remi Abgrall and Wasilij Barsukow.
\newblock Extensions of {A}ctive {F}lux to arbitrary order of accuracy.
\newblock {\em ESAIM: Mathematical Modelling and Numerical Analysis},
  57(2):991--1027, 2023.

\bibitem[AB23b]{abgrall22proceeding}
R{\'e}mi Abgrall and Wasilij Barsukow.
\newblock A hybrid finite element--finite volume method for conservation laws.
\newblock {\em Applied Mathematics and Computation}, 447:127846, 2023.

\bibitem[Abg22]{abgrall20}
R{\'e}mi Abgrall.
\newblock A combination of {R}esidual {D}istribution and the {A}ctive {F}lux
  formulations or a new class of schemes that can combine several writings of
  the same hyperbolic problem: application to the 1d {E}uler equations.
\newblock {\em Communications on Applied Mathematics and Computation}, pages
  1--33, 2022.

\bibitem[Bar21a]{barsukow19activeflux}
Wasilij Barsukow.
\newblock The active flux scheme for nonlinear problems.
\newblock {\em Journal of Scientific Computing}, 86(1):1--34, 2021.

\bibitem[Bar21b]{barsukow20cgk}
Wasilij Barsukow.
\newblock Truly multi-dimensional all-speed schemes for the euler equations on
  cartesian grids.
\newblock {\em Journal of Computational Physics}, 435:110216, 2021.

\bibitem[BB23]{barsukow20swaf}
Wasilij Barsukow and Jonas~P Berberich.
\newblock A well-balanced {A}ctive {F}lux method for the shallow water
  equations with wetting and drying.
\newblock {\em Communications on Applied Mathematics and Computation}, pages
  1--46, 2023.

\bibitem[BEK{\etalchar{+}}17]{barsukow16}
Wasilij Barsukow, Philipp~VF Edelmann, Christian Klingenberg, Fabian Miczek,
  and Friedrich~K R{\"o}pke.
\newblock A numerical scheme for the compressible low-{M}ach number regime of
  ideal fluid dynamics.
\newblock {\em Journal of Scientific Computing}, 72(2):623--646, 2017.

\bibitem[BHKR19]{barsukow18activeflux}
Wasilij Barsukow, Jonathan Hohm, Christian Klingenberg, and Philip~L Roe.
\newblock The active flux scheme on {C}artesian grids and its low {M}ach number
  limit.
\newblock {\em Journal of Scientific Computing}, 81(1):594--622, 2019.

\bibitem[BK22]{barsukow17}
Wasilij Barsukow and Christian Klingenberg.
\newblock Exact solution and a truly multidimensional {G}odunov scheme for the
  acoustic equations.
\newblock {\em ESAIM: M2AN}, 56(1), 2022.

\bibitem[CH23]{chudzik23}
Erik Chudzik and Christiane Helzel.
\newblock A review of cartesian grid active flux methods for hyperbolic
  conservation laws.
\newblock In {\em International Conference on Finite Volumes for Complex
  Applications}, pages 93--109. Springer, 2023.

\bibitem[CHK21]{chudzik21}
Erik Chudzik, Christiane Helzel, and David Kerkmann.
\newblock The cartesian grid active flux method: Linear stability and bound
  preserving limiting.
\newblock {\em Applied Mathematics and Computation}, 393:125501, 2021.

\bibitem[ER13]{eymann13}
Timothy~A Eymann and Philip~L Roe.
\newblock Multidimensional active flux schemes.
\newblock In {\em 21st AIAA computational fluid dynamics conference}, 2013.

\bibitem[FR15]{fan15}
Doreen Fan and Philip~L Roe.
\newblock Investigations of a new scheme for wave propagation.
\newblock In {\em 22nd AIAA Computational Fluid Dynamics Conference}, page
  2449, 2015.

\bibitem[GC90]{gresho90}
Philip~M Gresho and Stevens~T Chan.
\newblock On the theory of semi-implicit projection methods for viscous
  incompressible flow and its implementation via a finite element method that
  also introduces a nearly consistent mass matrix. {P}art 2: {I}mplementation.
\newblock {\em International Journal for Numerical Methods in Fluids},
  11(5):621--659, 1990.

\bibitem[HKS19]{kerkmann18}
Christiane Helzel, David Kerkmann, and Leonardo Scandurra.
\newblock A new {ADER} method inspired by the active flux method.
\newblock {\em Journal of Scientific Computing}, 80(3):1463--1497, 2019.

\bibitem[HS99]{hu99}
Changqing Hu and Chi-Wang Shu.
\newblock Weighted essentially non-oscillatory schemes on triangular meshes.
\newblock {\em Journal of Computational Physics}, 150(1):97--127, 1999.

\bibitem[LL98]{lax98}
Peter~D Lax and Xu-Dong Liu.
\newblock Solution of two-dimensional riemann problems of gas dynamics by
  positive schemes.
\newblock {\em SIAM Journal on Scientific Computing}, 19(2):319--340, 1998.

\bibitem[Mae17]{maeng17}
Jungyeoul Maeng.
\newblock {\em On the advective component of active flux schemes for nonlinear
  hyperbolic conservation laws}.
\newblock PhD thesis, University of Michigan, Dissertation, 2017.

\bibitem[MRKG03]{munz03}
C-D Munz, Sabine Roller, Rupert Klein, and Karl~J Geratz.
\newblock The extension of incompressible flow solvers to the weakly
  compressible regime.
\newblock {\em Computers \& Fluids}, 32(2):173--196, 2003.

\bibitem[PP11]{peraire11}
Jaime Peraire and Per-Olof Persson.
\newblock High-order discontinuous {G}alerkin methods for {CFD}.
\newblock In {\em Adaptive high-order methods in computational fluid dynamics},
  pages 119--152. World Scientific, 2011.

\bibitem[RLM15]{roe15}
Philip~L Roe, Tyler Lung, and Jungyeoul Maeng.
\newblock New approaches to limiting.
\newblock In {\em 22nd AIAA Computational Fluid Dynamics Conference}, page
  2913, 2015.

\bibitem[RMF18]{roe18}
Philip~L Roe, Jungyeoul Maeng, and Doreen Fan.
\newblock Comparing active flux and {D}iscontinuous {G}alerkin methods for
  compressble flow.
\newblock In {\em 2018 AIAA Aerospace Sciences Meeting}, page 0836, 2018.

\bibitem[Roe21]{roe21}
Philip Roe.
\newblock Designing {CFD} methods for bandwidth—a physical approach.
\newblock {\em Computers \& Fluids}, 214:104774, 2021.

\bibitem[vL77]{vanleer77}
Bram van Leer.
\newblock Towards the ultimate conservative difference scheme. {IV}. {A} new
  approach to numerical convection.
\newblock {\em Journal of computational physics}, 23(3):276--299, 1977.

\bibitem[Zen14]{zeng14}
Xianyi Zeng.
\newblock A high-order hybrid finite difference--finite volume approach with
  application to inviscid compressible flow problems: a preliminary study.
\newblock {\em Computers \& Fluids}, 98:91--110, 2014.

\bibitem[Zen19]{zeng19}
Xianyi Zeng.
\newblock Linear hybrid-variable methods for advection equations.
\newblock {\em Advances in Computational Mathematics}, 45(2):929--980, 2019.

\end{thebibliography}
\end{document}